\documentclass{gtart}     

\input gtmonout
\volumenumber{2}
\volumeyear{1999}
\volumename{Proceedings of the Kirbyfest}
\pagenumbers{489}{553}
\papernumber{24}
\received{10 September 1998}\revised{8 June 1999}
\published{22 November 1999}

\nocolon
\newtheorem{thm}{Theorem}[section]  
\newtheorem{cor}[thm]{Corollary}
\newtheorem{lemma}[thm]{Lemma} 
\newtheorem{prop}[thm]{Proposition} 
\theoremstyle{remark}
\newtheorem{defn}[thm]{Definition} 
\newtheorem{example}[thm]{Example} 
\newcommand{\aaa}{\mbox{$\alpha$}}
\newcommand{\map}{\mbox{$\rightarrow$}}

\newcommand{\Aaa}{\mbox{$\mathcal A$}}
\newcommand{\Fff}{\mbox{$\mathcal F$}}  
\newcommand{\bbb}{\mbox{$\beta$}}
\newcommand{\sss}{\mbox{$\sigma$}}  
 
\newcommand{\rrr}{\mbox{$\rho$}} 
\newcommand{\Ggg}{\mbox{$\Gamma$}}

\newcommand{\ttt}{\mbox{$\tau$}} 
\newcommand{\bdd}{\mbox{$\partial$}}

\newcommand{\Ss}{\mbox{$\Sigma$}}
\newcommand{\Ddd}{\mbox{$\Delta$}}

\newcommand{\Thp}{\mbox{$\Theta_P$}}
\newcommand{\Thq}{\mbox{$\Theta_Q$}}
\newcommand{\aub} {\mbox{$A \cup_{P} B$}}
\newcommand{\awb} {\mbox{$A_- \cup_{P_-} B_-$}}
\newcommand{\xuy} {\mbox{$X \cup_{Q} Y$}}
\newcommand{\xwy} {\mbox{$X_- \cup_{Q_-} Y_-$}}
\newcommand{\avb} {\mbox{$A' \cup_{P'} B'$}}
\newcommand{\xvy} {\mbox{$X' \cup_{Q'} Y'$}}
\newcommand{\px} {\mbox{$P_X$}}
\newcommand{\py} {\mbox{$P_Y$}}
\newcommand{\qa} {\mbox{$Q_A$}}
\newcommand{\qb} {\mbox{$Q_B$}}
\newcommand{\lll}{\mbox{$\lambda$}}

\newcommand{\pf}{\mbox{\bf Proof}\qua}
\input rlepsf
\let\relabela\adjustrelabel

\begin{document}

\title{Genus two Heegaard splittings\\ of orientable
three--manifolds}                     
\asciititle{Genus two Heegaard splittings of orientable
three-manifolds}                     

\authors{Hyam Rubinstein\\Martin Scharlemann}                  
\address{Department of Mathematics, University  of Melbourne\\
Parkville, Vic 3052, Australia}
\secondaddress{Mathematics Department, University of California\\
        Santa Barbara, CA 93106, USA}
\asciiaddress{Department of Mathematics, University  of Melbourne\\
Parkville, Vic 3052, Australia\\
Mathematics Department, University of California\\
Santa Barbara, CA 93106,  USA}

\email{rubin@ms.unimelb.edu.au, mgscharl@math.ucsb.edu}

\begin{abstract}   
It was shown by Bonahon--Otal and Hodgson--Rubinstein that any two genus--one
Heegaard splittings of the same $3$--manifold (typically a lens space) are
isotopic.  On the other hand, it was shown by Boileau, Collins and
Zieschang that certain  Seifert manifolds have distinct genus--two Heegaard
splittings.  In an earlier  paper, we presented a technique for comparing
Heegaard splittings of the same  manifold and, using this technique,
derived the uniqueness theorem  for lens space splittings as a simple
corollary.  Here we use a similar technique to examine, in general, ways in
which two non-isotopic genus--two Heegard  splittings of the same
$3$-manifold compare, with a particular focus on how the  corresponding
hyperelliptic involutions are related.
\end{abstract}

\asciiabstract{It was shown by Bonahon-Otal and Hodgson-Rubinstein
that any two genus-one Heegaard splittings of the same 3-manifold
(typically a lens space) are isotopic.  On the other hand, it was
shown by Boileau, Collins and Zieschang that certain Seifert manifolds
have distinct genus-two Heegaard splittings.  In an earlier paper, we
presented a technique for comparing Heegaard splittings of the same
manifold and, using this technique, derived the uniqueness theorem for
lens space splittings as a simple corollary.  Here we use a similar
technique to examine, in general, ways in which two non-isotopic
genus-two Heegard splittings of the same 3-manifold compare, with a
particular focus on how the corresponding hyperelliptic involutions
are related.}

\makeshorttitle
%
%

\section{Introduction}

It is shown in \cite{Bo}, \cite{HR} that any two genus one Heegaard
splittings of the same manifold (typically a lens space) are
isotopic.  On the other hand, it is shown in \cite{BGM}, \cite{Mo}
that certain Seifert manifolds have distinct genus two Heegaard
splittings (see also Section
\ref{seifert} below).  In \cite{RS} we present a technique for
comparing Heegaard splittings of the same manifold and derive the
uniqueness theorem for lens space splittings as a
simple corollary.  The intent of this paper is to use the
technique of \cite{RS} to examine, in general, how two genus two Heegard
splittings of the same $3$--manifold compare.

One potential way of creating genus two Heegaard split
$3$--manifolds is to ``stabilize'' a splitting of lower genus (see
\cite[Section 3.1]{Sc}).  But since, up to isotopy, stabilization is
unique and since genus one Heegaard splittings are
known to be unique, this process cannot produce non-isotopic
splittings.  So we may as well restrict to genus two splittings that
are not stabilized.  A second way of creating a $3$--manifold
equipped with a genus two Heegaard splitting is to take the
connected sum of two $3$--manifolds, each with a genus one splitting. 
But (again since genus one splittings are unique) any two Heegaard
splittings of the same manifold that are constructed in this way can be
made to coincide outside a collar of the summing sphere.  Within
that collar there is one possible difference, a ``spin'' corresponding
to the non-trivial element of $\pi_1 (RP(2))$, where $RP(2)$
parameterizes unoriented planes in $3$--space and the spin reverses the
two sides of the plane. Put more simply, the two splittings differ only
in the choice of which side of the torus in one summand is identified
with a given side of the splitting torus in the other summand.  The
first examples of this type are given in \cite{Mon}, \cite{Vi}.

 These easier cases having been
considered,  interest will now focus on genus two splittings that are
``irreducible'' (see \cite[Section 3.2]{Sc}).  It is a consequence of
\cite{CG} that a genus two splitting which is irreducible is also
``strongly irreducible'' (see \cite[Section 3.3]{Sc}, or the proof of
Lemma \ref{str.irred} below).  That is, if $M = \aub$ is a Heegaard
splitting, then any pair of essential disks, one in $A$ and one in $B$,
have boundaries that intersect in at least two points.

The result of our program is a listing, in Sections \ref{seifert} and
\ref{examples}, of all ways in which two strongly irreducible genus two
Heegaard splittings of the same closed orientable $3$--manifold $M$
compare.  The proof that this is an exhaustive listing is the subject
of the rest of the paper.  What we do not know is when two Heegaard
splittings constructed in the ways described are authentically
different.  That is, we do not have the sort of algebraic invariants
which would allow us to assert that there is no global isotopy of $M$
that carries one splitting into another.  For the case of Seifert
manifolds (eg
\cite{BoO}) such algebraic invariants can be (non-trivially) derived
from the very explicit form of the fundamental group.

Any $3$--manifold with a genus two Heegaard splitting admits an
orientation preserving involution whose quotient space is $S^3$ and
whose branching locus is a $3$--bridge knot (cf \cite{BH}).  The
examples constructed in Section \ref{examples} are sufficiently
explicit that we can derive from them global theorems.  Here are a
few:  If $M$ is an atoroidal closed orientable $3$--manifold then the
involutions coming from distinct Heegaard splittings necessarily
commute.  More generally, the commutator of two different involutions
can be obtained by some composition of Dehn twists around essential
tori in $M$.  Finally, two genus two splittings become equivalent
after a single stabilization.

The results we obtain easily generalize to compact orientable
$3$--manifolds with boundary, essentially by substituting boundary tori
any place in which Dehn surgery circles appear.  

We expect the methods and results here may be helpful in understanding
$3$--bridge knots (which appear as branch sets, as described above) and
in understanding the mapping class groups of genus two $3$--manifolds.

The authors gratefully acknowledge the support of, respectively, the 
Australian Research Council and both MSRI and the National Science
Foundation.

\section{Cabling handlebodies}
\label{cable}

Imbed the solid torus $S^1 \times D^2$ in ${\bf C^2}$ as $\{ (z_1, z_2)
| |z_1| = 1, |z_2| \leq 1 \}$.  Define a natural orientation-preserving
involution $\Theta\co   S^1 \times D^2 \map S^1 \times D^2$ by $\Theta
(z_1, z_2) =  (\overline{z}_1, \overline{z}_2)$. Notice that the
fixed points of $\Theta$ are precisely the two arcs $\{ (z_1, z_2)
| z_1 = \pm 1, -1 \leq z_2 \leq 1 \}$ and the quotient space is $B^3$. 
On the torus $S^1 \times \bdd D^2$ the fixed points of $\Theta$ are the
four points $\{ (\pm 1, \pm 1) \}$. 

For any pair of integers $(p, q) \neq (0,0)$ we can define the {\em $(p,
q)$ torus link} $L_{p,q} \subset S^1 \times \bdd
D^2$ to be  $\{ (z_1, z_2) \in S^1 \times \bdd
D^2 | (z_1)^p = (z_2)^q \}$.  The $(1,0)$ torus link is a meridian
and the $(k, 1)$ torus link is a longitude of the solid torus. A
{\em torus knot} is a torus link of one component which is not a
meridian or longitude.  In other words, a torus knot is a torus link in
which $p$ and $q$ are relatively prime, and neither is zero. Up to
orientation preserving homeomorphism of $S^1 \times \bdd D^2$ (given by
Dehn twists) we can also assume, for a torus knot, that $1 \leq p <
q$.

\medskip
{\bf Remark}\qua  Let $L_{p,q} \subset S^1 \times \bdd D^2$ be a torus
knot,
\aaa\ an arc that spans the annulus $S^1 \times \bdd D^2 - L_{p,q}$ and
\bbb\ be a radius of the disk $\{ point \} \times D^2 \subset S^1
\times D^2$. Then the complement of a neighborhood of $L_{p,q} \cup
\aaa$ in
$S^1 \times D^2$ is isotopic to a neighborhood of $(S^1 \times \{ 0
\}) \cup \bbb$ in $S^1 \times D^2$.  This fact is useful later in
understanding how cabling is affected by stabilization.

\medskip
Clearly $\Theta$ preserves any torus link
$L_{p,q}$.  If $L_{p,q}$ is a torus knot, so $p$ and $q$ can't
both be even, the involution $\Theta | L_{p,q}$ has precisely two fixed
points: $(1, 1)$ and either $(-1, -1)$, if $p$ and $q$ are both odd;
or $(-1, 1)$ if $p$ is even; or $(1, -1)$ if $q$ is even. This has the
following consequence.  Let $N$ be an equivariant neighborhood of the
torus knot $L_{p,q}$ in  $S^1 \times D^2$.  Then $N$ is
topologically a solid torus, and the fixed points of $\Theta | N$ are
two arcs.  That is, $\Theta | N$ is topologically conjugate to
$\Theta$. (See Figure 1.)


\begin{figure}[ht!]
\centerline{
\epsfxsize=76mm \epsfbox{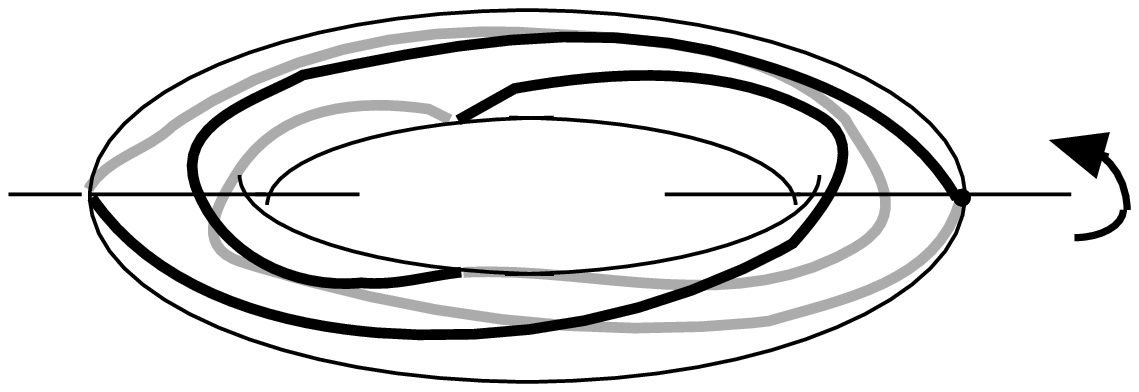}}
\caption{}
\end{figure} 

The involution $\Theta$ can be used to build an involution of a genus
two handlebody $H$ as follows.  Create $H$ by attaching together two
copies  of $S^1 \times D^2$ along an equivariant disk neighborhood
$E$ of $(1, 1) \in  S^1 \times \bdd D^2$ in each copy.  Then $\Theta$
acting simultaneously on both copies will produce an involution of $H$,
which we continue to denote $\Theta$. Again the quotient is $B^3$ but
the fixed point set consists of three arcs. (See Figure 2 for a
topologically equivalent picture.) We will call
$\Theta$ the {\em standard involution} on $H$.   It has the following
very useful properties:  it carries any simple closed curve in $\bdd H$
to an isotopic copy of the curve, and, up to isotopy, any homeomorphism
of
$\bdd H$ commutes with it. It will later be useful to distinguish
involutions of different handlebodies, and since, up to isotopy rel
boundary, this involution is determined by its action on $\bdd H$, it is
legitimate, and will later be useful, to denote the involution by 
$\Theta_{\bdd H}$.


\begin{figure}[ht!]
\centerline{\epsfxsize=72mm \epsfbox{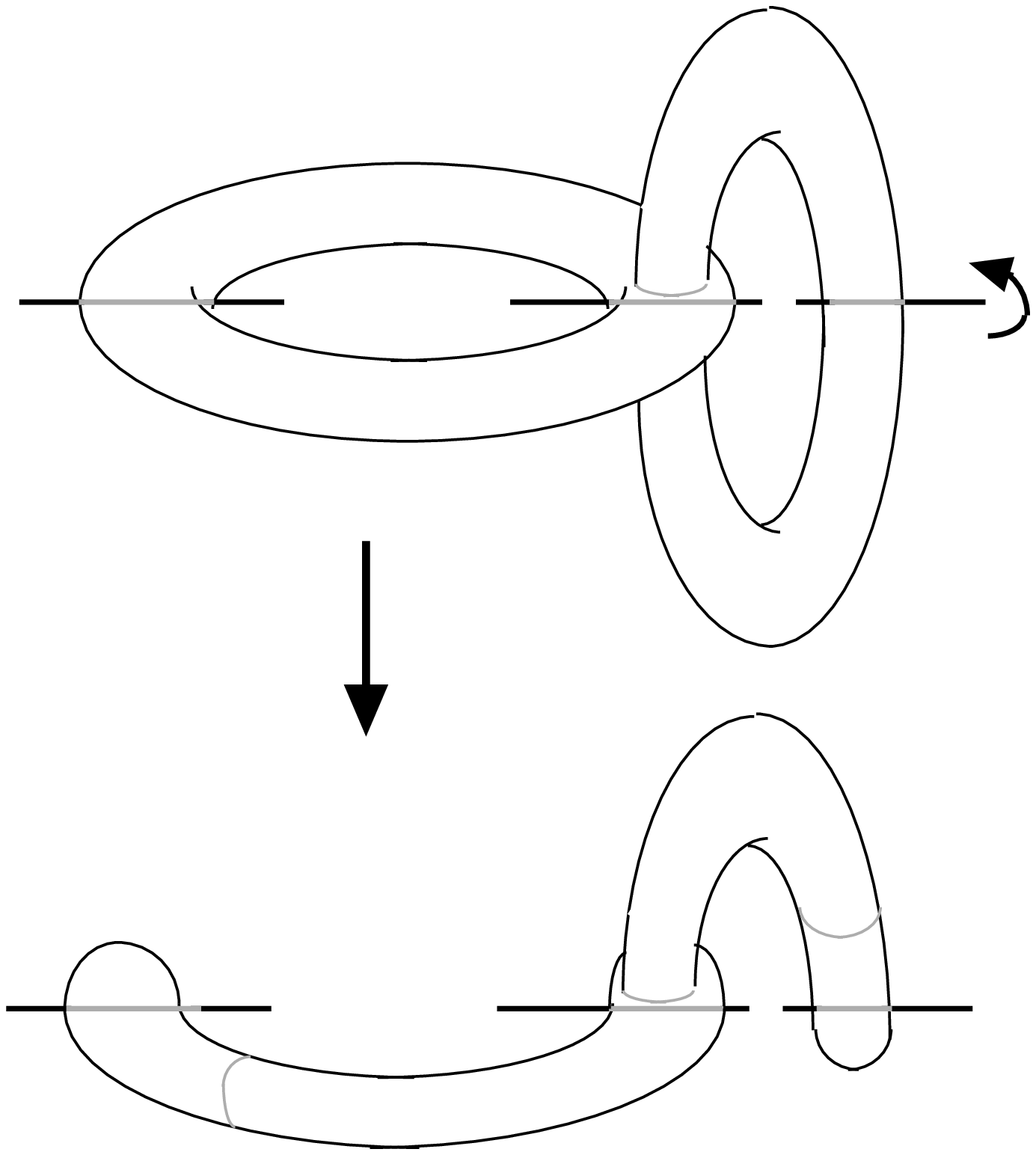}}
\caption{}
\end{figure}

Two alternative involutions of the genus two handlebody  $H = (S^1
\times D^2)_1 \cup_{E} (S^1 \times D^2)_2$ will
sometimes be useful.  Consider the involution that rotates $H$
around a diameter of $E$, exchanging $(S^1 \times D^2)_1$ and  $(S^1
\times D^2)_2$.  (See Figure 3.)  The diameter is the set of fixed
points, and the quotient space is a solid torus.  This will be called
the {\em minor involution} on $H$.  The final involution is best
understood by thinking of $H$ as a neighborhood of the union of two
circles that meet so that the planes containing them are
perpendicular, as in Figure 2.  Then $H$ is the union of two solid tori
in which a core of fixed points in one solid torus coincides in the
other solid torus with a diameter of a fiber.  Under
this identification, a full
$\pi$ rotation of one solid torus around its core coincides in the
other solid torus with the standard involution, and one of the arc of
fixed points in the second torus is a subarc of the core of the
first.  The quotient of this involution is a solid torus and the fixed
point set is the core of the first solid torus together with an
additional boundary parallel proper arc in the second solid torus. 
This involution will be called the {\em circular involution}.


\begin{figure}[ht!]
\centerline{
\epsfxsize=75mm \epsfbox{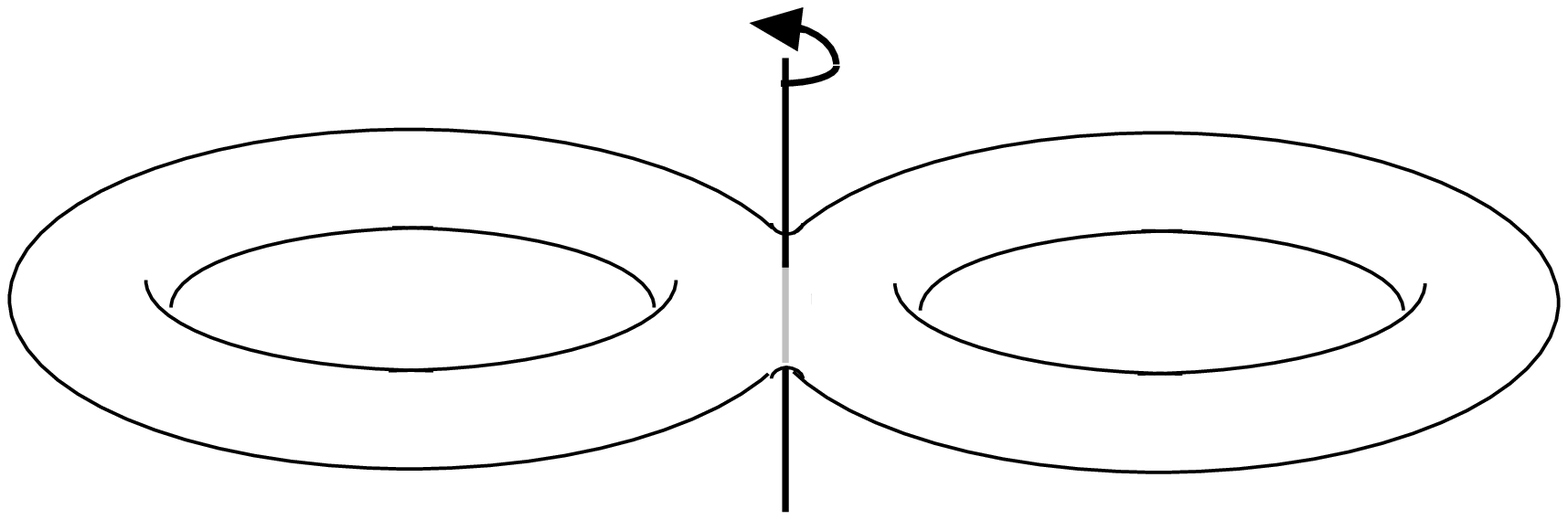}}
\caption{}
\end{figure}

In analogy to definitions in the case of a solid torus, we have:

\begin{defn}
\label{long}
A {\em meridian disk} of a handlebody $H$ is an essential disk
in $H$.  Its boundary is a {\em meridian curve}, or, more simply, a
meridian.  A {\em longitude} of $H$ is a simple closed curve in $\bdd
H$ that intersects some meridian curve exactly once.  A meridian disk
can be separating or non-separating.  Two longitudes $\lll, \lll'
\subset \bdd H$ are {\em separated longitudes} if they lie on opposite
sides of a separating meridian disk.
\end{defn}

There is a useful way of imbedding one genus two
handlebody in another.  Begin with $H = (S^1 \times D^2)_1 \cup_{E}
(S^1 \times D^2)_2$, on which $\Theta$ operates as above.  Let $N$ be
an equivariant neighborhood of a torus knot in $(S^1 \times D^2)_1$. 
Choose $N$ large enough (or $E$ small enough) that $E \subset N \cap 
\bdd (S^1 \times D^2)_1$.  Then  $H' = N \cup_{E}
(S^1 \times D^2)_2$ is a new genus two handlebody on which $\Theta$
continues to act.  In fact $\Theta_{\bdd H'} = \Theta_{\bdd H} | H'$. 
We say the handlebody $H'$ is obtained by {\em cabling into} $H$ or,
dually, $H$ is obtained by cabling out of $H'$. (See Figure 4.)


\begin{figure}[ht!]
\centerline{\relabelbox\small
\epsfxsize=82mm \epsfbox{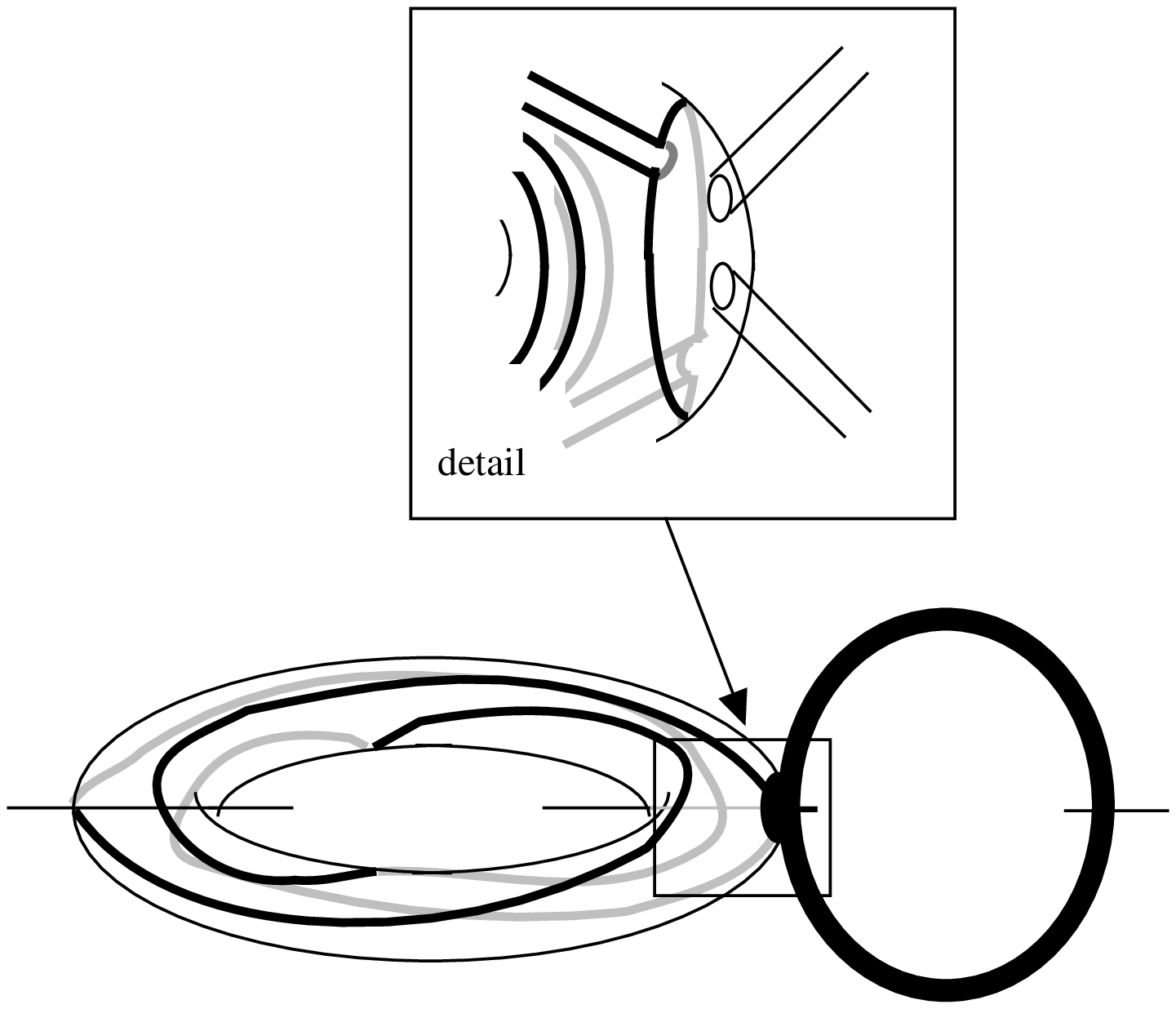}
\relabela <-2pt, 0pt> {detail}{detail}
\endrelabelbox}
\caption{}
\end{figure}

There is another useful way to view cabling into $H$. Recall the
process of Dehn surgery:  Let $q/s$ be a rational number and \aaa\ be a
simple closed curve in a $3$--manifold $M$.  Then we say a manifold $M'$
is obtained from $M$ by {\em $q/s$--surgery on \aaa} if a solid torus
neighborhood $\eta (\aaa)$ is removed from $M$ and is replaced by a
solid torus whose meridian is attached to $L_{q,s} \subset \bdd \eta
(\aaa)$.  Unless there is a natural choice of longitude for $\eta
(\aaa)$ (eg when $M = S^3$), $q/s$ is only defined modulo the
integers or, put another way, we can with no loss of generality take $0
\leq q < s$.  

If we take \aaa\ to be the core $S^1 \times \{ 0 \} \subset S^1 \times
D^2$ and perform $q/s$ surgery, then the result $M'$ is still a solid
torus, but $L_{q, s}$ becomes a meridian of $M'$.  The curve
$L_{p,r}$, with $ps - qr = 1$ becomes a longitude of $M'$, because it
intersects $L_{q,s}$ in one  point.  A longitude $L_{0,1}$ becomes the
torus knot $L_{p,q} \subset M'$ because it intersects a longitude $p$
times and a meridian $q$ times.  So another way of viewing $H' \subset
H$ is this:  Attach a neighborhood (containing $E$, but disjoint from
$\aaa \subset (S^1 \times D^2)_1$) of the longitude $(S^1 \times \{ 1
\})_1$ to $(S^1 \times D^2)_2$ to form $H'$.  Then do $q/s$ surgery to
\aaa\ to give $H$ containing $H'$. The advantage of this point of view
is that the construction is more obviously $\Theta$ equivariant (since
both the longitude and the core
\aaa\ are clearly preserved by $\Theta$) once we observe once and
for all, from the earlier viewpoint (see Figure 1), that Dehn surgery is
$\Theta$ equivariant.  

Of course it is also possible to cable into $H$ via a similar
construction in $(S^1 \times D^2)_2$, perhaps at the same time as we
cable in via $(S^1 \times D^2)_1$.

\section{Seifert examples of multiple Heegaard splittings}
\label{seifert}

A Heegaard splitting of a closed orientable $3$--manifold
$M$ is a decomposition $M = \aub$ in which $A$ and $B$ are
handlebodies, and $P = \bdd A = \bdd B$.  In other words, $M$ is
obtained by gluing two handlebodies together by some homeomorphism of
their boundaries. If the
splitting is genus two, then the splitting induces an involution on
$M$.  Indeed the standard involutions of $A$ and $B$ can be made to
coincide on $P$, since the standard involution on $A$, say, commutes
with the gluing homeomorphism
$\bdd A \map \bdd B$.  So we can regard $\Theta_P$ as an
involution of $M$ (cf \cite{BH}).   

We are interested in understanding closed orientable $3$--manifolds that
admit more than one isotopy class of genus two Heegaard splitting. That
is, splittings $M = \aub = \xuy$ in which the genus two surfaces $P$ and
$Q$ are not isotopic.  (A separate but related question is whether
there is an ambient {\em homeomorphism} which carries one to the other,
ie, whether the splittings are {\em homeomorphic}.)  In this section
we begin by discussing a class of manifolds for which the answer is
well understood.

It is a
consequence of the classification theorem of Moriah and Schultens
\cite{MS} that Heegaard splittings of closed Seifert manifolds (with
orientable base and fiberings) are either ``vertical'' or
``horizontal''.  The consequence which is relevant here is that any
such Seifert manifold which has a genus two splitting is in fact a
Seifert manifold over $S^2$ with three exceptional fibers.  Through
earlier work of Boileau and Otal \cite{BO} it was already known that
genus two splittings of these manifolds were either vertical or
horizontal and this led Boileau, Collins and Zieschang
\cite{BCZ} and, independently, Moriah \cite{Mo} to give a complete
classification of genus 2 Heegaard decompositions in this case.  In
general, there are several.

Most (the vertical splittings) can be constructed as follows: 
Take regular neighborhoods of two exceptional fibers and connect them
with an arc (transverse to the fibering) that projects to an imbedded
arc in the base space connecting the two exceptional points, which are
the projections of the exceptional fibers.  Any two such arcs are
isotopic, so the only choice involved is in the pair of exceptional
fibers.  It is shown in
\cite{BCZ} that this choice can make a difference---different choices
can result in Heegaard splittings that are not even homeomorphic.

The various vertical splittings do have one common property, however. 
They all share the same standard involution.  All that is involved in
demonstrating this is the proper construction of the involution on the
Seifert manifold $M$.  In the base space, put all three exceptional
fibers on the equator of the sphere.  Now consider the
orientation preserving involution of $M$ that simultaneously reverses
the direction of every fiber and reflects the base $S^2$ through the
equator (ie, exchanges the fiber lying over a point with the fiber
lying over its reflection).  This involution induces the natural
involution on a neighborhood of any fiber that lies over the equator,
specifically the exceptional fibers.  If we choose two of them, and
connect them via a subarc of the equator, the involution on $M$ is the
standard involution on the corresponding Heegaard splitting.

Two types of Seifert manifolds have additional splittings (see
\cite[Proposition 2.5]{BCZ}).  One, denoted $V(2,3,a)$, is the $2$--fold
branched cover of the $3$--bridge torus knot $L_{(a, 3)} \subset S^3, a
\geq 7$, and the other, denoted $W(2, 4, b)$, is a $2$--fold branched
cover over the $3$--bridge link which is the union of the torus knot
$L_{(b,2)} \subset S^3, b \geq 5$ and the core of the solid torus on
which it lies.  Since these are three--bridge links, there is a
sphere that divides them each into two families of three unknotted arcs
in $B^3$.  The $2$--fold branched cover of three unknotted arcs in
$B^3$ is just the genus two handlebody (in fact the inverse
operation to quotienting out by the standard involution).  So this
view of the links defines a Heegaard splitting on the
double--branched cover $M$.

In both cases the natural fibering of $S^3$ by torus knots of the
relevant slope lifts to the Seifert fibering on the double--branched
cover.  The torus knots lie on tori, each of which induces a genus
one Heegaard splitting of $S^3$.  The natural involution of $S^3$
defined by this splitting (rotation about an unknot \aaa\ in $S^3$ that
intersects the cores of both solid tori, see
\cite[Figures 4 and 5]{BCZ}), preserves the fibering of
$S^3$ and induces the natural involution on any fiber that
intersects its axis.  We can arrange that the exceptional fibers
(including those on which we take the $2$--fold branched cover) 
intersect \aaa.  Then the standard involution of
$S^3$ simultaneously does three things.  It induces the standard
involution \Thp\ on $M$ that comes from its vertical Heegaard
decomposition $M = \aub$; it preserves the $3$--bridge link
$L \subset S^3$ that is the image of the fixed point set of the
other involution \Thq; and it preserves the sphere which lifts to the
other Heegaard surface $Q \subset M = \xuy$, while interchanging $X$
and $Y$.  It follows easily that \Thp\ and \Thq\ commute.

The product of the two involutions is again the standard involution,
but with a different axis of symmetry. To see how this can be, note that
the involution \Thq\ is in fact just a flow of $\pi$ along each
regular fiber and also along the exceptional fibers other than the
branch fibers.  Since the branch fibers have even index, a flow of
$\pi$ on regular fibers induces the identity on the
branch fibers. So in fact \Thq\ is isotopic to the identity (just
flow along the fibers).  The fixed point set of \Thp\ intersects any
exceptional fiber in two points, $\pi$ apart; indeed it is a reflection
of the fiber across those two fixed points. Hence 
\Thq\ carries the fixed point set of
\Thp\ to itself and the involutions commute.  The composition $\Thp
\Thq$ is also a reflection in each exceptional fiber---but through a
pair of points which differ by $\pi /2$ from the points across which,
in an exceptional fiber, \Thp\ reflects.  See Figure 5.

\begin{figure}[ht!] 
\centerline{\relabelbox\small
\epsfxsize=75mm \epsfbox{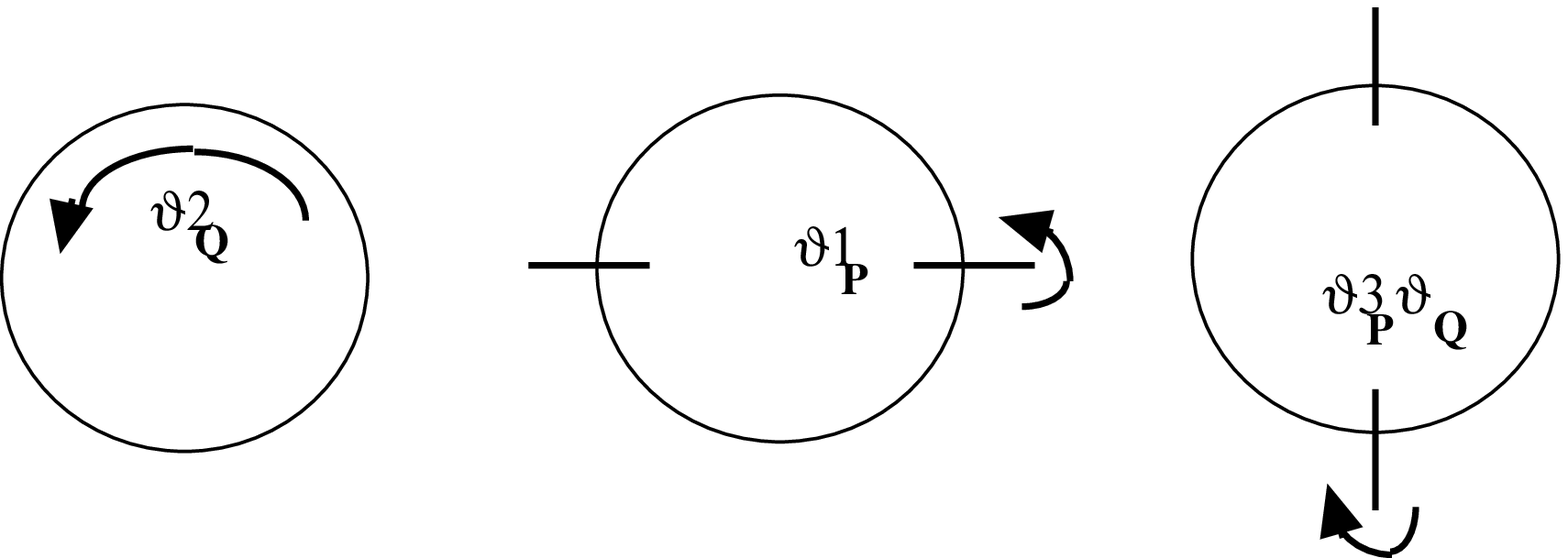}
\relabela <-6pt,-2pt> {J1}{\Thp}
\relabela <-0pt,-6pt>  {J2}{\Thq}
\relabela <-6pt,5pt>  {J3}{$\Thp\Thq$}
\endrelabelbox}
\caption{}
\end{figure}

\section{Other examples of multiple Heegaard splittings}
\label{examples}

In this section we will list a number of ways of constructing
$3$--manifolds $M$, not necessarily Seifert manifolds, which support
multiple genus two Heegaard splittings.  That is, it will follow from
the construction that $M$ has two or more Heegaard splittings which are
at least not obviously isotopic.  The constructions are elementary
enough that in all cases it will be easy to see that a single
stabilization suffices to make them isotopic. (We will only rarely
comment on this stabilization property.) They are symmetric enough that
in all cases we will be able to see directly how the corresponding
involutions of $M$ are related.  When
$M$ contains no essential separating tori then, in many cases, the
involutions from the different Heegaard splittings will be the same
and, in all cases, the involutions will at least commute.  When $M$
does contain essential separating tori, the same will be true after
possibly some Dehn twists around essential tori.

\begin{defn}

Suppose $T^2 \times I \subset M$ is a collar of an essential torus in
a compact orientable $3$--manifold $M$.  Then a homeomorphism $h\co  M
\map M$ is obtained by a {\em Dehn twist} around $T^2 \times \{ 0 \}$ if
$h$ is the identity on $M - (T^2 \times I)$.

\end{defn}

Ultimately we will show that any manifold that admits multiple
splittings will do so because the manifold, and any pair of different
splittings, appears on the list below.  This will allow us to make
conclusions about how the involutions determined by multiple splittings
are related.  What we are unable to determine is when the examples
which appear here actually do give non-isotopic splittings.  For this
one would need to demonstrate that there is no {\em global} isotopy of
$M$ carrying one splitting to another.  This requires establishing
a property of the splitting that is invariant under Nielsen moves and
showing that the property is different for two splittings.  For example,
the very rich structure of Seifert manifold fundamental groups was
exploited in \cite{BCZ} to establish that some splittings were
even globally non-isotopic.

Alternatively, as in \cite{BGM}, one could show that the associated
involutions have fixed sets which project to inequivalent knots or
links in $S^3$.  Note that non-isotopic splittings can probably arise
even when the associated involutions have fixed sets projecting tot he
same knots or links in $S^3$.  In this case, there would be
inequivalent $3$--bridge representations of these knots or links.

\subsection{Cablings}
\label{symcable}

Consider the graph $\Ggg \subset S^2 \subset S^3$ consisting of two
orthogonal polar great circles.  One polar circle will be denoted \lll\
and the other will be thought of as two edges $e_a$ and $e_b$ attached
to \lll.  The full $\pi$ rotation $\Pi_{\mu}\co  S^3 \map S^3$ about the
equator $\mu$ of $S^2$ preserves $\Ggg$.  (Here ``full $\pi$
rotation'' means this:  Regard $S^3$ as the join of $\mu$ with
another circle, and rotate this second circle half-way round.)  Without
changing notation, thicken
$\Ggg$ equivariantly, so it becomes a genus three handlebody and note
that on the two genus two sub-handlebodies
$\lll
\cup e_a$ and
$\lll \cup e_b$, $\Pi_{\mu}$ restricts to the standard involution.  

Now divide the solid torus $\lll$ in two by a longitudinal annulus
\Aaa\ perpendicular to $S^2$. The annulus $\Aaa$ splits $\lll$ into
two solid tori $\lll_a$ and $\lll_b$.  Both ends of the $1$--handle
$e_a$ are attached to $\lll_a$ and both ends of $e_b$ to $\lll_b$.
Define genus two handlebodies $A$ and $B$ by $A = \lll_a \cup e_a$  and
$B = \lll_b \cup e_b$.  Then $\Pi_{\mu}$ preserves $A$ and $B$  and on
them restricts to the standard involution.  Finally, construct a closed
$3$--manifold $M$ from $\Ggg$ by gluing $\bdd A - \Aaa$ to $\bdd B -
\Aaa$ by any homeomorphism (rel boundary).   Such a $3$--manifold $M$
and genus two Heegaard splitting $M = \aub$ is characterized by the
requirement that a longitude of one handlebody is identified with a
longitude of the other.  See Figure 6.


\begin{figure}[ht!] 
\centerline{\relabelbox\small
\epsfxsize=100mm \epsfbox{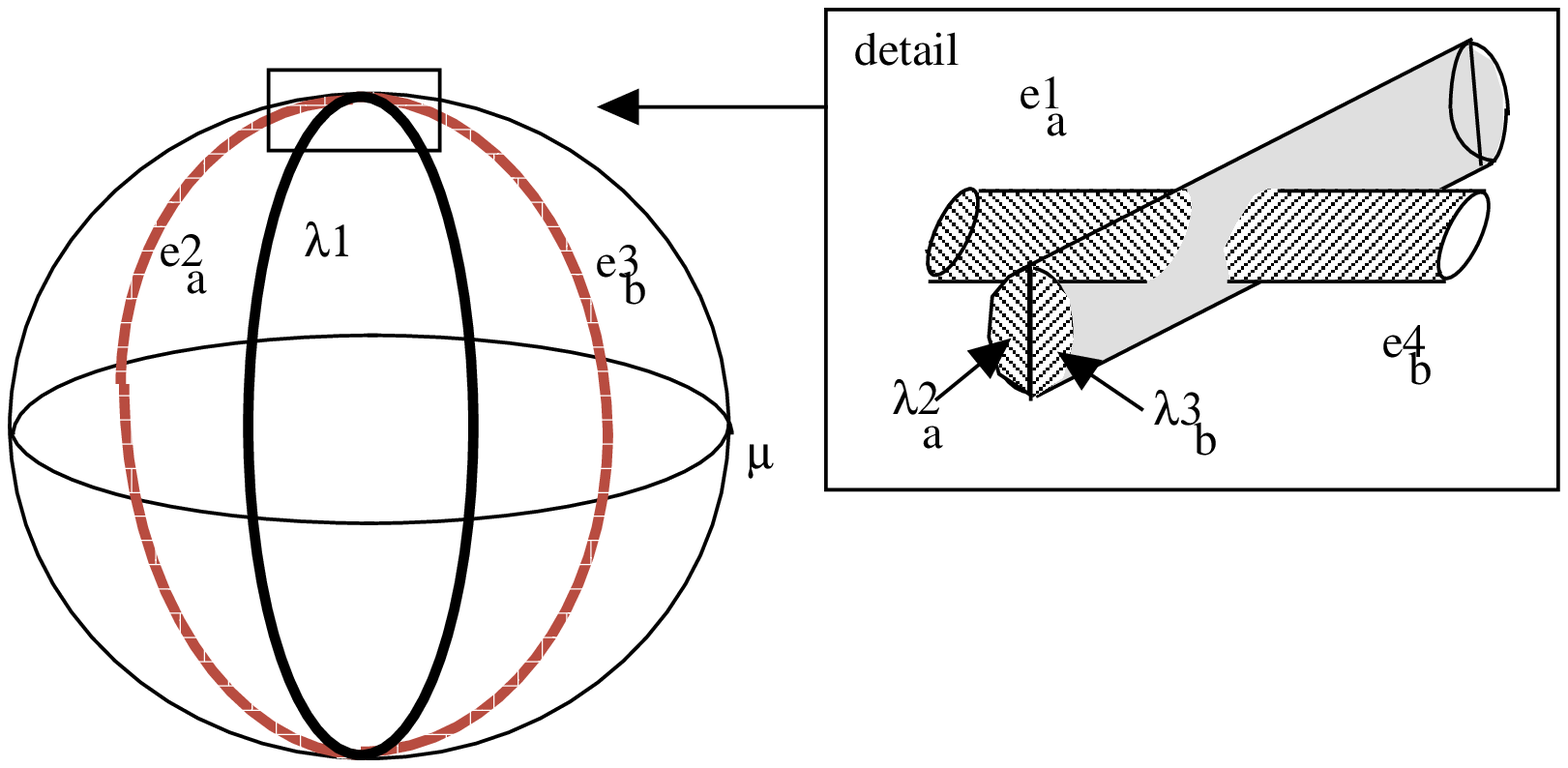}
\relabel {detail}{detail}
\relabel {l1}{$\lambda$}
\relabela <-2pt,0pt> {l2}{$\lambda_a$}
\relabela <-1pt,0pt> {l3}{$\lambda_b$}
\relabel {m}{$\mu$}
\relabela <0pt,-7pt> {e1}{$e_a$}
\relabel {e2}{$e_b$}
\relabel {e3}{$e_a$}
\relabela <0pt,2pt> {e4}{$e_b$}
\endrelabelbox}
\caption{}
\end{figure}

\medskip
{\bf Question}\qua  Which $3$--manifolds have such Heegaard
splittings?
\medskip

So far we have described a certain kind of Heegaard splitting, but have
not exhibited multiple splittings of the same $3$--manifold.  But such
examples can easily be built from this construction:  Let $\aaa_a$ and
$\aaa_b$ be the core curves of $\lll_a$ and $\lll_b$ respectively.  

\medskip
{\bf Variation 1}\qua   Alter $M$ by Dehn surgery on $\aaa_a$, and
call the result $M'$.  The splitting surface $P$ remains a Heegaard
splitting surface for $M'$, but now a longitude of $B' = B$ is attached
to a twisted curve in $\bdd A'$.  Since $\aaa_a$ and $\aaa_b$ are
parallel in $M$, we could also have gotten $M'$ by the same Dehn surgery
on $\aaa_b$.  But the isotopy from $\aaa_a$ to $\aaa_b$ crosses $P$,
so the splitting surface is apparently different in the two
splittings.  In fact, one splitting surface is obtained from the other
by cabling out of $B'$ and into $A'$.  It follows from the Remark in
Section \ref{cable} that the two become equivalent after a single
stabilization.

\medskip
{\bf Variation 2}\qua   Alter $M$ by Dehn surgery on both $\aaa_a$ and
$\aaa_b$ and call the result $M'$.  (Note that $M'$ then
contains a Seifert submanifold.)  In $\lll$ the annulus \Aaa\
separates the two singular fibers $\aaa_a$ and $\aaa_b$.  New
splitting surfaces for $M'$ can be created by replacing \Aaa\ by any
other annulus in $\lll$ that separates the
singular fibers and has the same boundary .  There are an integer's
worth of choices, basically because the braid group $B_2 \cong {\bf
Z}$.  Equivalently, alter $P$ by Dehn twisting around the separating
torus
$\bdd \lll$.

\subsection{Double cablings}
\label{doubcable}

Just as the previous example of symmetric cabling is a special case of
Heegaard splittings, so the example here of double cablings is
a special case of the symmetric cabling above, with additional
parts of the boundaries of $A$ and $B$ identified.  

Consider the graph
in $\Ggg \subset S^2 \subset S^3$ consisting of two circles $\mu_n$ and
$\mu_s$ of constant latitude, together with two edges  $e_a$ and $e_b$
spanning the annulus between them in $S^2$.  Both
$e_a$ and $e_b$ are segments of a polar great circle $\lll$.  The full
$\pi$ rotation $\Pi_{\lll}\co  S^3 \map S^3$ about $\lll$ preserves
$\Ggg$.  Without changing notation, thicken $\Ggg$ equivariantly, so it
becomes a genus three handlebody and note that on the two genus two
sub-handlebodies $\mu_n \cup e_a \cup \mu_s$ and $\mu_n \cup e_b
\cup \mu_s$, $\Pi$ restricts to the standard involution.

Now remove from both $\mu_n$ and $\mu_s$ annuli $\Aaa_n$ and $\Aaa_s$
respectively, chosen so that the boundary of each of the annuli is the
$(2,2)$ torus link in the solid torus.  That is, each boundary
component is the $(1,1)$ torus knot, where a preferred longitude of
the solid torus $\mu_n$ or $\mu_s$ is that determined by intersection
with $S^2$.  Place $\Aaa_n$ and $\Aaa_s$ so that they are
perpendicular to $S^2$ at the points where the edges $e_a$ and $e_b$
are attached.  Then $\Aaa_n$ divides $\mu_n$ into two solid tori, one
of them $\mu_{na}$ attached to one end of $e_a$ and the other
$\mu_{nb}$ attached to an end of $e_b$.  The annulus $\Aaa_s$ similarly
divides the solid torus $\mu_s$.  

Define genus two handlebodies $A$ and $B$ by $A = \mu_{na} \cup e_a
\cup \mu_{sa}$  and $B = \mu_{nb} \cup e_b \cup \mu_{sb}$.  Then
$\Pi_{\lll}$ preserves $A$ and $B$  and on them restricts to the
standard involution.  Finally, construct a closed $3$--manifold $M$ from
$\Ggg$ by gluing $\bdd A - (\Aaa_n \cup \Aaa_s)$ to  $\bdd B - (\Aaa_n
\cup \Aaa_s)$ by any homeomorphism (rel boundary).   Such a
$3$--manifold $M$ and genus two Heegaard splitting $M = \aub$ is
characterized by the requirement that two separated longitudes of one
handlebody are identified with two separated longitudes of the other. 
See Figure 7.


\begin{figure}[ht!] 
\centerline{\relabelbox\small
\epsfxsize=105mm \epsfbox{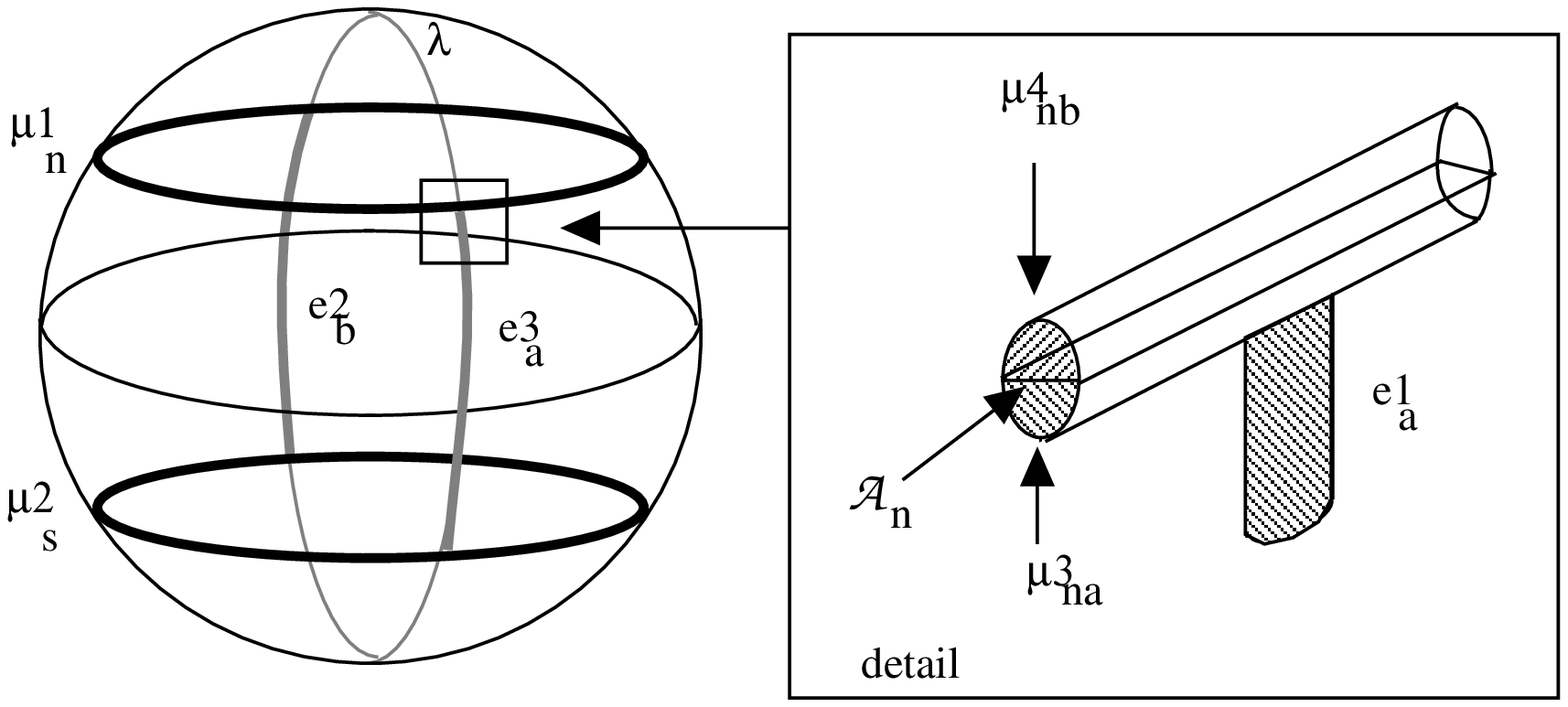}
\relabel {detail}{detail}
\relabela <0pt,-2pt> {l}{$\lambda$}
\relabel {m1}{$\mu_n$}
\relabel {m2}{$\mu_s$}
\relabel {m3}{$\mu_{na}$}
\relabela <0pt,-4pt> {m4}{$\mu_{nb}$}
\relabel {e1}{$e_a$}
\relabel {e2}{$e_b$}
\relabel {e3}{$e_a$}
\relabel {A}{$\mathcal A_n$}
\endrelabelbox}
\caption{}
\end{figure}

\medskip
{\bf Question}\qua  Which $3$--manifolds have such Heegaard
splittings?
\medskip

Just as in example \ref{symcable}, manifolds with multiple Heegaard
splittings can easily be built from this construction:  

\medskip
{\bf Variation 1}\qua  Let $\aaa_{na}, \aaa_{nb}, \aaa_{sa}, \aaa_{sb}$
be the core curves of $\mu_{na}, \mu_{nb}, \mu_{sa}$ and $\mu_{sb}$
respectively. Do Dehn surgery on one or more of these curves,
changing $M$ to $M'$.  If a single Dehn surgery is done in $\mu_n$
and/or $\mu_s$ then there is a choice on which of the possible core
curves it is done.  If two Dehn surgeries are done in $\mu_n$ and/or
$\mu_s$ then there is an integer's worth of choices of replacements for
$\Aaa_n$ and/or $\Aaa_s$, corresponding to Dehn twists around $\bdd
\mu_n$ and/or $\bdd \mu_s$. Up to such Dehn twists, all these Heegaard
splittings induce the same natural involution on $M'$. 

\medskip
{\bf Variation 2}\qua   Let $\rrr_a$ be a simple closed curve in the
$4$ punctured sphere $\bdd A \cap \bdd \Ggg$ with the property that 
$\rrr_a$ intersects the separating meridian disk orthogonal to $e_a$
exactly twice and a meridian disk of each of
$\mu_{na}$ and $\mu_{sa}$ in a single point.  Similarly define
$\rrr_b$.  Suppose the gluing homeomorphism $h\co  \bdd A \cap
\bdd \Ggg \map \bdd B \cap \bdd \Ggg$ has $h(\rrr_a) = \rrr_b$, and
call the resulting curve \rrr.  

Push $\rrr$ into $A$ and do any Dehn surgery on the curve.  Since
$\rrr$ is a longitude of $A$ the result is a handlebody. 
Similarly, if the curve were pushed into $B$
before doing surgery, then $B$ remains a handlebody. So this gives
two alternative splittings. But this is not new, since this
construction is obviously just a special case of a single cabling
(Example \ref{symcable}).  However, if we do surgery on
the curve \rrr\ after pushing into $A$ and simultaneously do surgery
on one or both of $\aaa_{nb}$ and $\aaa_{sb}$ we still get a Heegaard
splitting.  Now push \rrr\ into $B$ and simultaneously move the other
surgeries to $\aaa_{na}$ and/or $\aaa_{sa}$ and get an alternative
splitting.

\subsection{Non-separating tori}
\label{nonsep}

Let $\Ggg, \Aaa_n, \Aaa_s$ and the four core curves $\aaa_{na},
\aaa_{nb}, \aaa_{sa}, \aaa_{sb}$ be defined as they were in the
previous case, Example \ref{doubcable}, but now consider the
$\pi$--rotation $\Pi_{\mu}$ that rotates $S^3$ around the equator $\mu$
of $S^2$.  This involution preserves $\Ggg$ and the
$1$--handles $e_a$ and $e_b$, but it exchanges north and south, so
$\mu_n$ is exchanged with $\mu_s$, and $\Aaa_n$ with $\Aaa_s$.  Remove
small tubular neighborhoods of core curves $\aaa_n$ and $\aaa_s$ of the 
solid tori $\mu_n$ and $\mu_s$, and with them small core sub-annuli
of $\Aaa_n$ and $\Aaa_s$. Choose these neighborhoods so that they are
exchanged by $\Pi_{\mu}$ and call their boundary tori $T_n$ and $T_s$. 
Attach $T_n$ to $T_s$ by an orientation reversing homeomorphism $h$
that identifies the annulus $T_n \cap \mu_{n a}$ with $T_s \cap
\mu_{s a}$ and the annulus $T_n \cap \mu_{n b}$ with $T_s \cap
\mu_{s b}$. Choose $h$ so that the orientation reversing
composition $\Pi_{\mu} h\co  T_n \map T_n$ fixes two meridian circles,
$\ttt_+$ and $\ttt_-$, lying respectively on the meridian disks of
$\mu_n$ at which $e_a$ and $e_b$ are attached. The resulting manifold
$\Ggg_T$ is orientable and, in fact, is homeomorphic to $T^2 \times I$
with two $1$--handles attached.  Let $T \subset \Ggg_T$ be the
non-separating torus which is the image of $T_s$ (and so also $T_n$).  
Also denote by $\ttt_{\pm a}$ the two arcs of intersection
of $\ttt_+$  and $\ttt_-$ with $A$; these arcs lie on the longitudinal
annulus $A \cap T$.  Similarly denote the two arcs $\ttt_{\pm} \cap B$
by $\ttt_{\pm b}$.  See Figure 8.


\begin{figure}[ht!] 
\centerline{\relabelbox\small
\epsfxsize=105mm \epsfbox{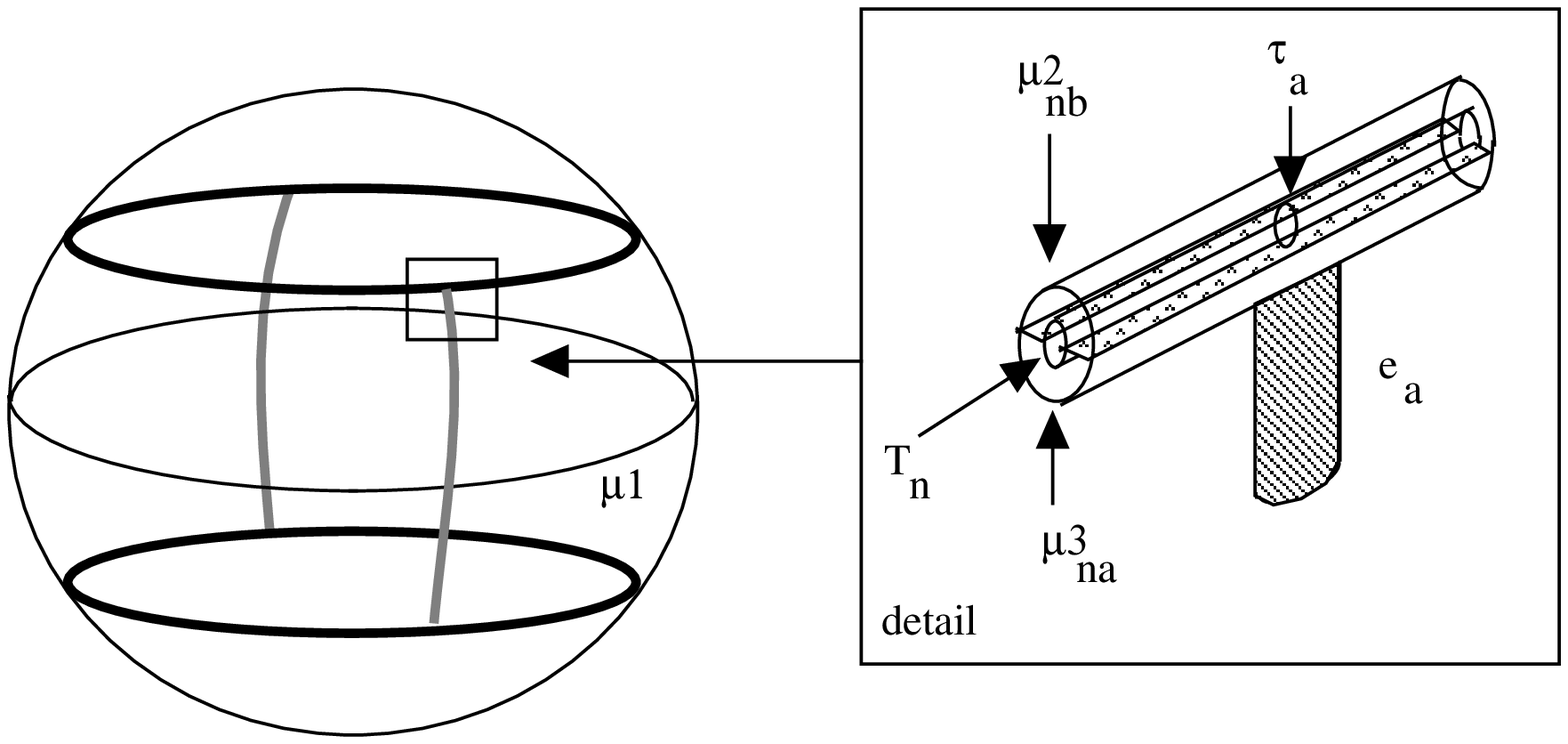}
\relabel {detail}{detail}
\relabel {m1}{$\mu$}
\relabel {m3}{$\mu_{na}$}
\relabela <0pt, -2pt> {m2}{$\mu_{nb}$}
\relabel {e}{$e_a$}
\relabela <0pt, -2pt> {t}{$\tau_a$}
\relabela <-1pt, -1pt> {T}{$T_n$}
\endrelabelbox}
\caption{}
\end{figure}

Note that in $\Ggg_T$ the union of $e_a, \mu_{n a}$, and $\mu_{sa}$ is
a genus two handlebody $A$ that intersects $T$ in a longitudinal
annulus. Similarly, the remainder is a genus two handlebody $B$ that
also intersects $T$ in a longitudinal
annulus.  The involution
$\Pi_{\mu}$ acts on $\Ggg_T$, preserves $T$ (exchanging its two
sides and fixing the two meridians $\ttt_{\pm}$), and
preserves both $A$ and $B$. The fixed points of the involution on $A$
consist of the arc $\mu \cap e_a$ and also the two arcs $\ttt_{\pm
a}$.  It is easy to see that this is the standard involution on $A$,
and, similarly, $\Pi_{\mu} | B$ is the standard involution.  Now glue
together the $4$--punctured spheres $\bdd A \cap \bdd \Ggg$ and $\bdd B
\cap \bdd \Ggg$ by any homeomorphism rel boundary.   The resulting
$3$--manifold $M$ and genus two Heegaard splitting $M = \aub$ has
standard involution $\Theta_P$ induced by $\Pi_{\mu}$.  The splitting is
characterized by the requirement that two distinct longitudes of one
handlebody, coannular within the handlebody,  are identified with two
similar longitudes of the other. 

\medskip
{\bf Question}\qua  Which $3$--manifolds have such Heegaard
splittings?
\medskip

Much as in the previous examples, manifolds with multiple Heegaard
splittings can be built from this construction:  

\medskip
{\bf Variation 1}\qua    We can assume that the deleted
neighborhoods of $\aaa_n$ and $\aaa_s$ in the
construction of $M$ above were small enough to leave the
parallel core curves $\aaa_{na}, \aaa_{nb}, \aaa_{sa},
\aaa_{sb}$ intact. Do Dehn surgery on 
$\aaa_{n a}$ (or, equivalently, $\aaa_{s a}$), changing $M$ to $M'$. 
The same manifold $M'$ can be obtained by doing the same Dehn surgery
to   $\aaa_{n b}$ (or, equivalently, $\aaa_{s b}$), but the Heegaard
splittings are not obviously isotopic, for they differ  by cabling into
$A$ and out of $B$.  

\medskip
{\bf Variation 2}\qua Do Dehn surgery on both 
$\aaa_{n a}$ and $\aaa_{n b}$ (or, equivalently, both $\aaa_{s a}$ and
$\aaa_{s b}$), changing $M$ to $M'$.  This inserts two singular fibers
in the collar $T^2 \times I$ between $\bdd \mu_s$ and $\bdd \mu_n$ and
these are separated by two spanning annuli, the remains of the annuli
$\Aaa_n$ and $\Aaa_s$ glued together. View this region as a Seifert
manifold, with two exceptional fibers, over the annulus $S^1 \times
I$. Let $p_a, p_b$ denote the projections of the two exceptional 
fibers to the annulus $S^1 \times I$.  The spanning annuli
project to two spanning arcs in
$(S^1
\times I) -
\{ p_a, p_b
\}$. There is a choice of such spanning arcs, and so of spanning annuli
between
$\bdd \mu_s$ and $\bdd \mu_n$ that still produce a Heegaard splitting. 
The choices of arcs all differ by braid moves in $(S^1 \times I) - \{
p_a, p_b \}$, and these correspond to Dehn twists around essential tori
in $M'$.

\medskip
{\bf Variation 3}\qua  This variation does not involve Dehn surgery. Let
$R_A$ be the long rectangle that cuts the $1$--handle $e_a$ down the
middle, intersecting every disk fiber of $e_a$ in a single diameter,
always perpendicular to $S^2$.  Extend $R_A$ by attaching meridian
disks of $\mu_{na}$ and $\mu_{nb}$ so the ends of $R_A$ become
identified to $\ttt_{+ a}$. Since the identification is orientation
reversing, $R_A$ becomes a M{\"o}bius band in $A$, corresponding to
the M{\"o}bius band spanned by $L_{1,2}$ in one of the solid tori
summands of $A$. Define $R_B$ similarly, but add a half-twist, so that
$R_B$ becomes a non-separating longitudinal annulus in $B$.  

Now construct $M$ as
above, choosing a gluing homeomorphism $\bdd A \cap \bdd \Ggg \map
\bdd B
\cap \bdd \Ggg$ so that $R_A \cap \bdd \Ggg$ ends up disjoint from
$R_B \cap \bdd \Ggg$. There are an integer's worth of possibilities
for this gluing, corresponding to Dehn twists around the annulus
complement of the two spanning arcs of $R_A$ in the $4$--punctured
sphere $\bdd A \cap \bdd \Ggg$.  The four arcs of $R_A$ and $R_B$
divide the
$4$--punctured sphere into two disks.

Let $Y$ be the genus $2$ handlebody obtained from $B$ in two
steps:  First remove a collar neighborhood of the annulus $R_B$,
cutting $B$ open along a longitudinal annulus.  At this point
$\Pi_{\mu}$ is the minor involution on $Y$, since the half-twist in 
$R_B$ means that it contains the arc $\mu \cap e_b$ as well as the arc
$\ttt_{+ a}$.  To get the standard involution on $Y$,  $\pi$
rotate around an axis in $S^3$ perpendicular to $S^2$ and
passing through the points where $\mu$ intersects the cores of $e_a$ and
$e_b$.  Call this rotation $\Pi_{\perp}$.  Two arcs of fixed points lie
in the disk fiber (now split in two) where $e_b$ crosses $\mu$.  A
third arc of fixed points, more difficult to see, is what remains of a
core of the annulus $T \cap A$, once a neighborhood of $\ttt_{+ a}$ is
removed. 

Next attach a neighborhood of the M{\"o}bius band $R_A$ to $Y$.  One
can see that it is attached along a longitude of $Y$, so the effect is
to cable out of $Y$ into its complement---$Y$ remains a handlebody. 
Moreover, $\Pi_{\perp}$ still induces the standard involution on $Y$.

Similarly, if a neighborhood of $R_A$ is removed from $A$ the effect
is to cable into $A$ and if a neighborhood of $R_B$ is then attached
the result is still a handlebody $X$.  The Heegaard decomposition $M =
\xuy$ has standard involution $\Theta_Q$ induced by $\Pi_{\perp}$, since
$Y$ did.  Notice that $\Pi_{\perp}$ and $\Pi_{\mu}$ commute, with
product $\Pi_{\lll}$, so $\Theta_P$ and $\Theta_Q$ commute.  The
product involution $\Theta_P \Theta_Q$ has fixed point set in $B$
(resp.\ $X$) the core circle of $R_B$ and an additional arc which crosses
$\ttt_{- b}$ in a single point.  That is, it is the ``circular
involution'' on both handlebodies (and also on $A$ and $Y$).

\subsection{$K_4$ examples}
\label{K4} 

Let $K_4$ denote the complete graph on $4$ vertices.  Construct a
complex $\Ggg$, isomorphic to $K_4$, in $S^2$ as follows.  Let $\mu$
denote the equator and
$\lll_a, \lll_b$ two orthogonal polar great circles.  Let the edge
$e_a$ be the part of
$\lll_a$ lying in the upper hemisphere and the edge $e_b$ be the part
of $\lll_b$ that lies in the lower hemisphere. Then take $\Ggg = \mu
\cup e_a \cup e_b$. Without changing notation, thicken $\Ggg$
equivariantly, so it becomes a genus three handlebody. See Figure 9.


\begin{figure}[ht!] 
\centerline{\relabelbox\small
\epsfxsize=118mm \epsfbox{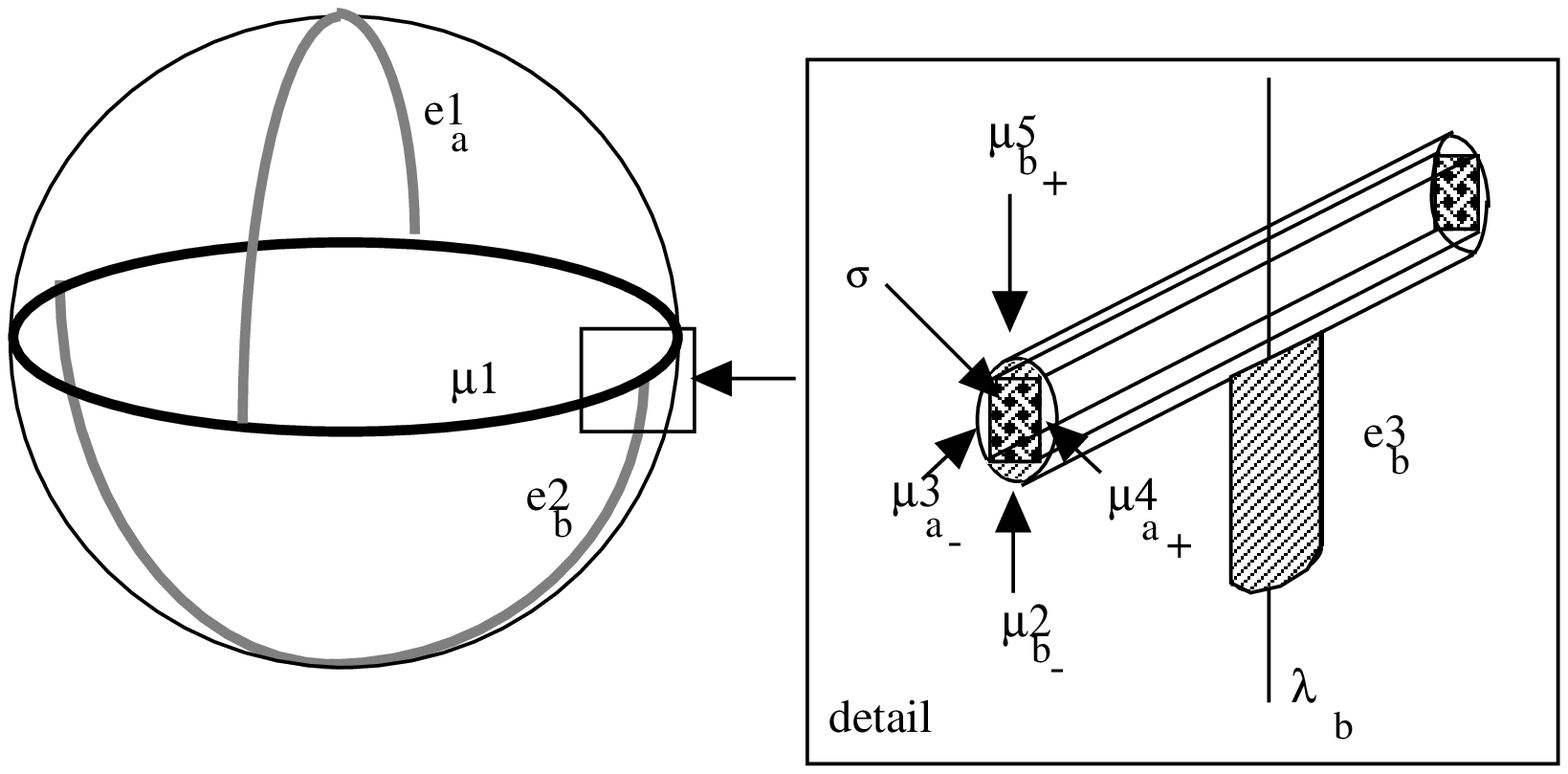}
\relabel {detail}{detail}
\relabel {m1}{$\mu$}
\relabela <0pt, -2pt> {m5}{$\mu_{b_+}$}
\relabela <-1pt, 0pt> {m2}{$\mu_{b_-}$}
\relabela <-2pt, 0pt> {m3}{$\mu_{a_-}$}
\relabela <-2pt, 0pt> {m4}{$\mu_{a_+}$}
\relabel {e1}{$e_a$}
\relabel {e2}{$e_b$}
\relabel {e3}{$e_b$}
\relabel {s}{$\sigma$}
\relabel {l}{$\lambda_b$}
\endrelabelbox}
\caption{}
\end{figure}

Consider the two $\pi$--rotations $\Pi_{a}, \Pi_{b}$ that rotate $S^3$
around, respectively, the curves $\lll_a, \lll_b$.  Both involutions
preserve $\Ggg$ and preserve also the individual parts $\mu, e_a$ and 
$e_b$.  Notice that $\Pi_{a}$ induces the minor involution on the
genus two handlebody $\mu \cup e_a$ and the standard involution on the
genus two handlebody $\mu \cup e_b$.  The symmetric statement is
true for $\Pi_b$.

Consider the link $L_{4,4} \subset \mu$.  The link intersects any
meridian disk of $\mu$ in four points.  Let \sss\ denote the
inscribed ``square'' torus $(S^1 \times square) \subset \mu$ which, in
each meridian disk of $\mu$, is the convex hull of those four points. 
The complementary closure of \sss\ in $\mu$ consists of four solid
tori. Isotope $L_{4,4}$ so that two of the complementary solid tori,
$\mu_{a_{\pm}}$, lying on opposite sides of
\sss\, are attached to $e_a$, each to one end of $e_a$.  The other two,
$\mu_{b_{\pm}}$ are then similarly attached to $e_b$.  Notice that,
paradoxically,  $\Pi_{a}$ now induces the standard involution on the
genus two handlebody $A_- = \mu_{a_{+}} \cup \mu_{a_{-}}
\cup e_a$ and the minor involution on the genus two handlebody 
$B_- = \mu_{b_{+}} \cup \mu_{b_{-}}
\cup e_b$. The latter is because $\lll_a$ is disjoint from both 
of $\mu_{b_{\pm}}$ and so only intersects the handlebody in a
diameter of a meridian disk of $e_b$.   The symmetric statements are
of course  true for $\Pi_b$.  Finally, let $M_-$ be a $3$--manifold 
obtained by gluing together the $4$--punctured spheres $\bdd A_- \cap
\bdd \Ggg$ and $\bdd B_- \cap \bdd \Ggg$ by any homeomorphism rel
boundary.  Note that so far we have not identified any Heegaard
splitting of $M$, since \sss\ is in neither $A_-$ nor $B_-$.

\medskip
{\bf Variation 1}\qua  Let $A = A_-$ and $B = B_- \cup \sss$. Then $M =
\aub$ is a Heegaard splitting, on which $\Pi_a$ is the standard
involution.  Indeed, we've already seen that $\Pi_a$ is standard on $A
= A_-$ and it is standard on $B$ since $\lll_a$ passes twice through
$\sss\ \subset B$, as well as once through $e_b$.  Alternatively, let
$X = A_- \cup \sss$ and $Y = B_-$.  Then, for exactly the same
reasons, $M = \xuy$ is a Heegaard splitting, on which $\Pi_b$ acts as
the standard involution.  Notice that $\Pi_a$ and $\Pi_b$ commute. 
Their product is $\pi$ rotation about the circle perpendicular to
$S^2$ through the poles.  This is the minor involution on {\em both}
$A$ and $B$.  It follows that $\Theta_P$ and $\Theta_Q$ commute and
their product operates as the minor involution on all four of $A, B, X,
Y$.  

\medskip
{\bf Variation 2}\qua
Let $M'$ be obtained by a Dehn surgery on the
core of \sss.  The constructions of Variation 1 give two Heegaard
splittings of $M'$ as well, with commuting standard involutions.  But
more splittings are available as well:  $A$ could be cabled into $B$ in
two different ways, essentially by moving the Dehn surgered circle into
either of $\mu_{a_{\pm}}$.  Similarly $Y$ could be cabled into $X$. 
Since such cablings have the same standard involution, the various
alternatives give involutions which either coincide or commute.

\medskip
{\bf Variation 3}\qua Let $M'$ be obtained by Dehn surgery on two
parallel circles in $\sss$.  These can be placed in a variety of
locations and still we would have Heegaard splittings:  If at most one
is placed as a core of $\mu_{a_+}$ or $\mu_{a_-}$ and the other is left
in \sss, then still $M' = \avb$ is a Heegaard splitting.  Similarly if
one is put in $\mu_{a_+}$ and the other in $\mu_{a_-}$.  In both cases
the splittings can additionally be altered by Dehn twists around the
now essential torus $\bdd \mu$.  We could similarly move one or both of
the two surgery curves into $\mu_{b_{\pm}}$ to alter the splitting $M' =
\xvy$.  Finally, we could move one into $\mu_{a_{\pm}}$ and the other
into $\mu_{b_{\pm}}$.  Then respectively $\avb$ and $\xvy$ are
alternative splittings.

\medskip
{\bf Variation 4}\qua Let $M'$ be obtained by Dehn surgery on three
parallel circles in $\sss$.  If one is placed in each of $\mu_{a_+}$
and $\mu_{a_-}$ and the third is left in \sss\ we still have a
Heegaard splitting $M' = \avb$.  Moreover, there is then a choice of
how the pair of annuli $P' \cap \mu$ lie in $\mu$.  The surgeries
change $\mu$ into a Seifert manifold over a disk, with three exceptional
fibers lying over singular points $p_{a_+}, p_{a_-}$ and $p_{\sss}$. 
The annuli $P' \cap \mu$ lie over proper arcs in the disk, which can be
altered by braid moves on the singular points.  These braid moves
translate to Dehn twists about essential tori in $M'$.   We could
similarly arrange the three surgery curves with respect to
$\mu_{b_{\pm}}$ to alter the splitting $M' = \xvy$.  

\medskip
{\bf Variation 5}\qua Let $\rrr_a$ be a simple closed curve in the
$4$ punctured sphere $\bdd A_- \cap \bdd \Ggg$ with the property that 
$\rrr_a$ intersects the separating meridian disk orthogonal to $e_a$
exactly twice and a meridian disk of each of
$\mu_{a_{\pm}}$ in a single point.  Similarly define $\rrr_b$. 
Suppose the gluing homeomorphism $h\co  \bdd A_- \cap
\bdd \Ggg \map \bdd B_- \cap \bdd \Ggg$ has $h(\rrr_a) = \rrr_b$.  

Push $\rrr_a$ into $A_-$ and do any Dehn surgery on the curve.  Since
$\rrr_a$ is a longitude of $A_-$ the result is a handlebody $A'$.  The
complement is the handlebody $B$ of Variation 1.  On the other hand,
if the curve (identified with $\rrr_b$) were pushed into $B_-$ before
doing surgery, then $B_-$ remains a handlebody $Y'$ and its complement
is the handlebody $X$ of Variation 1.  So this pair of alternative
splittings, $M = A' \cup_{P} B = X \cup_{Q} Y'$,  is in some sense a
variation of variation 1.

\medskip
{\bf Variation 6}\qua  Just as Variation 5 is a modified Variation 1,
here we modify Variations 2 and 3.  Suppose curves
$\rrr_a$ and $\rrr_b$ are identified as in Variation 5, and do Dehn
surgery on this curve \rrr. But also do another Dehn surgery on one
or two curves parallel to the core of \sss, as in Variation 2.  If \rrr\
is pushed into
$A_-$ and at most one of the other Dehn surgery curves is put in
each of $A$ and $B$ then
$\avb$ is a Heegaard splitting.  If \rrr\ is pushed into $B_-$ and
at most one of the other Dehn surgery curves is put in each of $X$
and $Y$, then
$\xvy$ is a Heegaard splitting.

\medskip
{\bf Variation 7}\qua  Topologically, $\sss \cong S^1 \times D^2$;
choose a framing so that $L_{1,1} \subset \bdd \sss$ is identified
with $S^1 \times \{ point \}$.  Remove the interior of \sss\ from
\Ggg\ and identify $\bdd \sss \cong S^1 \times \bdd D^2$ to itself by an
orientation reversing involution $\iota$ that is a reflection in the
$S^1$ factor and a $\pi$ rotation in $\bdd D^2$. In particular
$\iota$ identifies the two longitudinal annuli $A_- \cap \sss$ (resp.\
$B_- \cap \sss$).  Hence, after the identification given by $\iota$,
$A_-$ (resp.\ $B_-$) becomes a genus two handlebody $A$ (resp.\ $B$).
A closed  $3$--manifold  can then be obtained by
gluing together the $4$--punctured spheres $\bdd A_- \cap
\bdd \Ggg$ and $\bdd B_- \cap \bdd \Ggg$ by any homeomorphism rel
boundary.  Equivalently, the closed manifold $M$ is obtained from an
$M_-$ (with boundary a torus) constructed as in the initial discussion
above by identifying the torus $\bdd M_-$ to itself by an orientation
reversing involution.  The quotient of the torus is a Klein bottle $K
\subset M$, whose neighborhood typically is bounded by the
canonical torus of $M$.  

To create from this variation examples of a single manifold with
multiple splittings, apply the same trick as in earlier variations: Do
Dehn surgery on the core curve of $\mu_{b_{+}}$ (equivalently
$\mu_{b_{-}}$) and/or the core curve of $\mu_{a_{+}}$ (equivalently
$\mu_{a_{-}}$).  If we do the surgery on one curve (so the set of 
canonical tori becomes a torus cutting off a Seifert piece,
fibering over the M{\"o}bius band with one exceptional fiber) then
there is a choice of whether the curve lies in $A_-$ or $B_-$.  If we
do surgery on two curves (so the Seifert piece fibers over the
M{\"o}bius band with two exceptional fibers) then there is a choice of
which vertical annulus in the Seifert piece becomes the intersection
with the splitting surface.  In the former case the standard
involutions of the two splittings are the same and in the latter they
differ by Dehn twists about an essential torus.  

\section{Essential annuli in genus two handlebodies}

It's a consequence of the classification of surfaces that on an 
orientable surface of genus $g$ there is, up to homeomorphism, 
exactly one  non-separating simple closed curve and $[g/2]$ separating
simple closed curves.   For the genus two surface $F$, this means
that each collection
$\Ggg$ of  disjoint simple closed curves is determined up to
homeomorphism by a
$4$--tuple of non-negative  integers:
$(a, b, c, d)$ where $a \geq b \geq c$ and $c \cdot d = 0$ (see Figure 
10). Denote this $4$--tuple by $I(\Ggg)$.


\begin{figure}[ht!] 
\centerline{\relabelbox\small
\epsfxsize=75mm \epsfbox{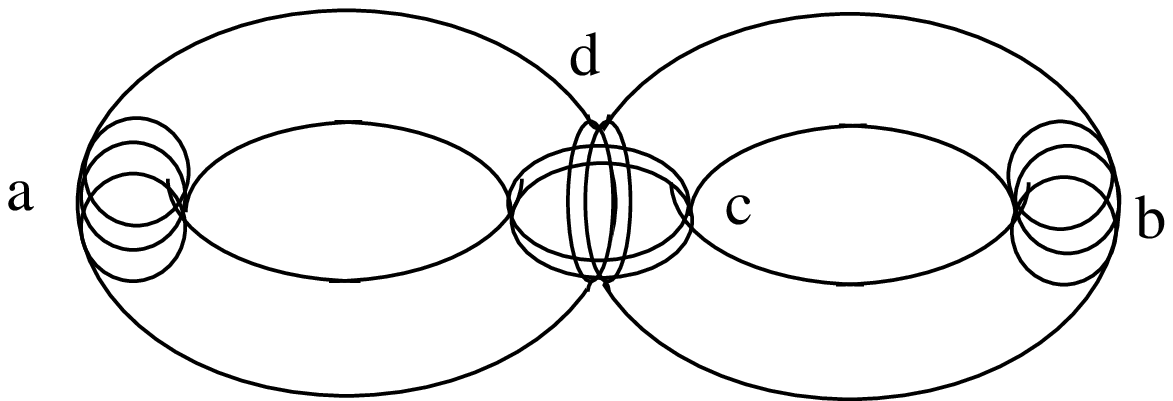}
\relabel {a}{$a$}
\relabel {b}{$b$}
\relabel {c}{$c$}
\relabela <2pt,0pt> {d}{$d$}
\endrelabelbox}
\caption{}
\end{figure}

Any collection of simple closed curves might occur as the boundary of 
some disks in a genus two handlebody and any collection of an even
number of curves might also occur as the boundary of some annuli in a
genus two handlebody, just by taking $\bdd$--parallel annuli or
tubing together disks.  To avoid such trivial constructions define:  

\begin{defn} 

A properly imbedded surface $S$ in a compact orientable $3$--manifold
$M$ is {\em essential} if $S$ is incompressible and no component of $S$ 
is $\bdd$--parallel.

\end{defn}

\begin{lemma}
\label{annuli}

Suppose $\Aaa \subset H$ is a collection of disjoint essential annuli in 
a genus
$2$ handlebody $H$.  Then $I(\bdd \Aaa) = (k, l, 0, 0)$ where $l \geq 0$ 
and $k + l$ is even.

\end{lemma}

\pf  Since $\Aaa \subset H$ is incompressible, it is
 $\bdd$--compressible.  Let $D$ be the disk obtained by a single
$\bdd$--compression.   Note that the effect of the $\bdd$--compression on
$\bdd \Aaa$ is to band sum two  distinct curves together. The band
cannot lie in an annulus in $\bdd H$ between the curves, since 
$\Aaa$ is not $\bdd$--parallel. So if $I(\bdd \Aaa) = (k, l, m, 0), m >
0$ or $(k,  l, 0, n), n > 0$, the band must lie in a pair of pants
component of $\bdd H -
\bdd  \Aaa$.  In that case  $\bdd D$ would be parallel to a component of
$\bdd A$, contradicting the assumption that $\Aaa$ is incompressible.

Finally, $k + l$ is even since each component of $\Aaa$ has two boundary
components. \qed

\begin{lemma}

Suppose $S \subset H$ is an essential oriented properly imbedded surface 
in a genus $2$ handlebody $H$ and  $\chi (S) = -1$.  Suppose that $[S]$
is trivial in $H_2(H, \bdd  H)$, and that no component of $S$ is a
disk.  Then  $I(\bdd A) = (k, l, 1, 0)$ or $(k, l, 0, 1)$.  

\end{lemma}

\pf  $S$ is $\bdd$--compressible, but the first $\bdd$--compression
can't be of an annulus component.  Indeed, the result of such a
$\bdd$--compression would be an essential disk in $H$ disjoint from $S$.
If we cut open along this disk, it would change $H$ into either one or
two solid tori.  But the only  incompressible surfaces that can be
imbedded in a solid torus are the disk and the  annulus, so $\chi (S) =
0$, a contradiction.  We conclude that the first $\bdd$--compression is
along a component $S_0$ with $\chi (S_0) = -1$.

After $\bdd$--compression $S_0$ becomes an annulus $A$.  If $A$ were
$\bdd$--parallel then the part of $S_0$ which was $\bdd$--compressed either 
lies in the region of parallelism or outside it.  In the former case,
$S_0$  would have been compressible and in the latter case it would
have been $\bdd$--parallel.  Since neither is allowed, we conclude that
$A$ is not $\bdd$--parallel.  So after the $\bdd$--compression the surface
becomes an essential collection of disjoint annuli, and Lemma
\ref{annuli} applies.

We now examine the possibilities other than those in the conclusion and 
deduce
a contradiction in each case.

\medskip
{\bf Case 1}\qua  $I(\bdd S) = (k, l, 0, n), n > 1$.

The $\bdd$--compression is into one of the complementary components and 
can
reduce $n$ by at most $1$.  So after the $\bdd$--compression the last
coordinate is still non-trivial, contradicting Lemma \ref{annuli}.

\medskip
{\bf Case 2}\qua  $I(\bdd S) = (k, l, m, 0), m > 1$.

Since $k \geq l \geq m$ the complementary components are annuli and two 
pairs
of pants.  The $\bdd$--compression then reduces $m$ by at most one, 
yielding
the same contradiction to Lemma \ref{annuli}.  

\medskip
{\bf Case 3}\qua  $I(\bdd S) = (k, l, 0, 0)$.

Since $\chi(S) = -1$, $k + l$ is odd, hence either $k$ or $l$ is odd.  
Then
there is a simple closed curve in $\bdd H$ intersecting $S$ an odd number 
of
times, contradicting the triviality of $[S]$ in $H_2(H, \bdd H)$. \qed

\medskip
{\bf  Remark}\qua  It is only a little harder to prove the same result, 
without the
assumption that $[S] = 0$, but then there is the additional
possibility that $I(S) = (1, 0, 0, 0)$. 

\begin{defn}
\label{twist}

Suppose $H$ is a handlebody and $c \subset \bdd H$ is a simple closed 
curve. Then $c$ is {\em twisted} if there is a properly imbedded disk
in $H$ which is disjoint from $c$ and, in the solid torus
complementary component
$S^1 \times D^2 \subset H$ in which $c$ lies, $c$ is a torus knot 
$L_{(p, q)}, p
\geq 2$ on $\bdd (S^1 \times D^2)$.

\end{defn}

\begin{defn}

A collection of annuli, all of whose boundary components are longitudes
is called {\em longitudinal}.  If all are twisted, then the collection 
is called {\em twisted}.

\end{defn}
 
Figures 11--13 show annuli which are respectively longitudinal, twisted
and non-separating, and twisted and separating.  Displayed in the
figure is an ``icon'' meant to schematically present the particular
annulus.  The icon is inspired by imagining taking a cross-section of
the handlebody near where the two solid tori are joined.  The
cross-section is of a meridian of the horizontal torus in the
handlebody figure together with part of the vertical torus.  Such
icons will be useful in presenting rough pictures of how families of
annuli combine to give tori in $3$--manifolds.


\begin{figure}[ht!] 
\centerline{\relabelbox\small
\epsfxsize=98mm \epsfbox{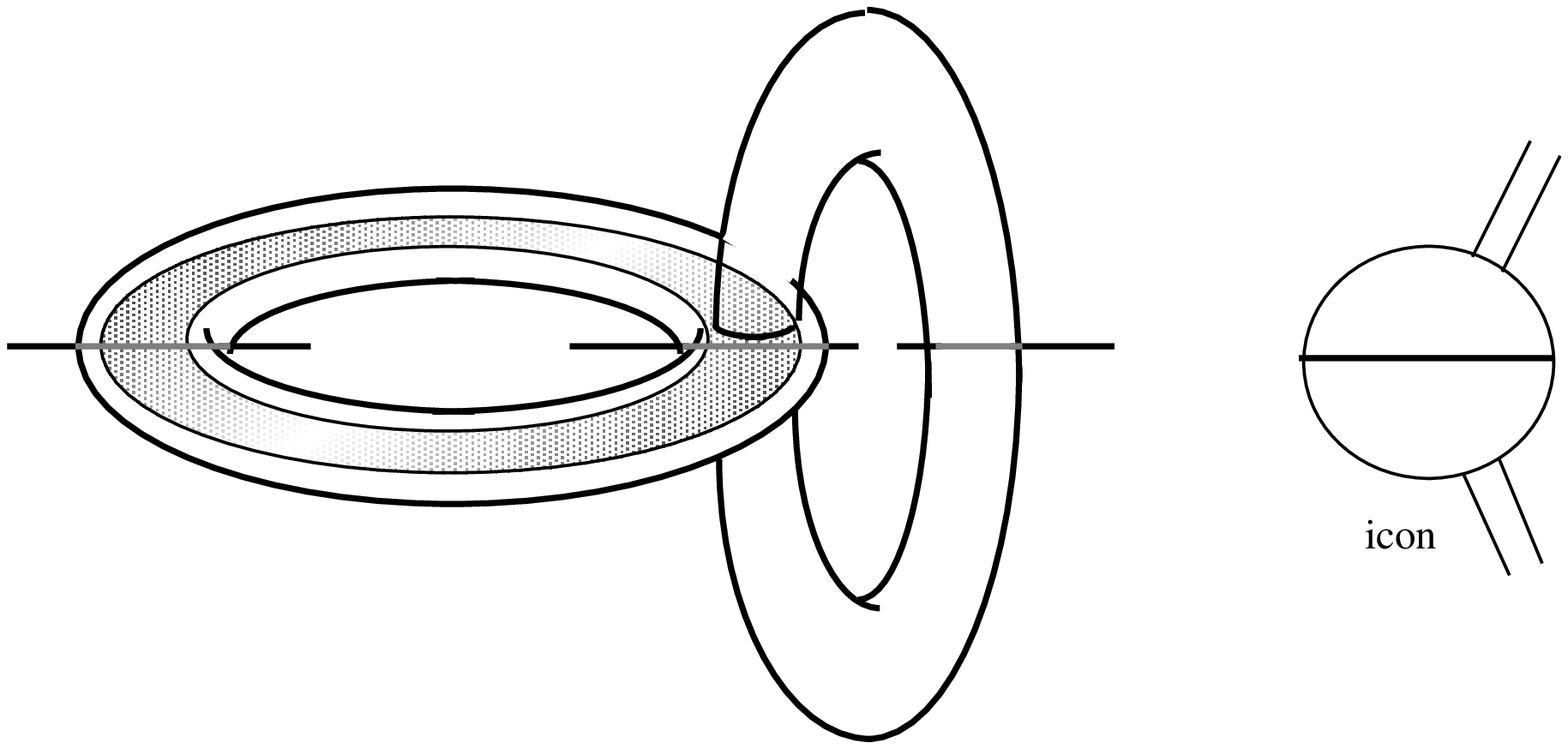}
\relabela <-5pt,0pt> {icon}{icon}
\endrelabelbox}
\caption{}
\end{figure}


\begin{figure}[ht!] 
\centerline{\relabelbox\small
\epsfxsize=84mm \epsfbox{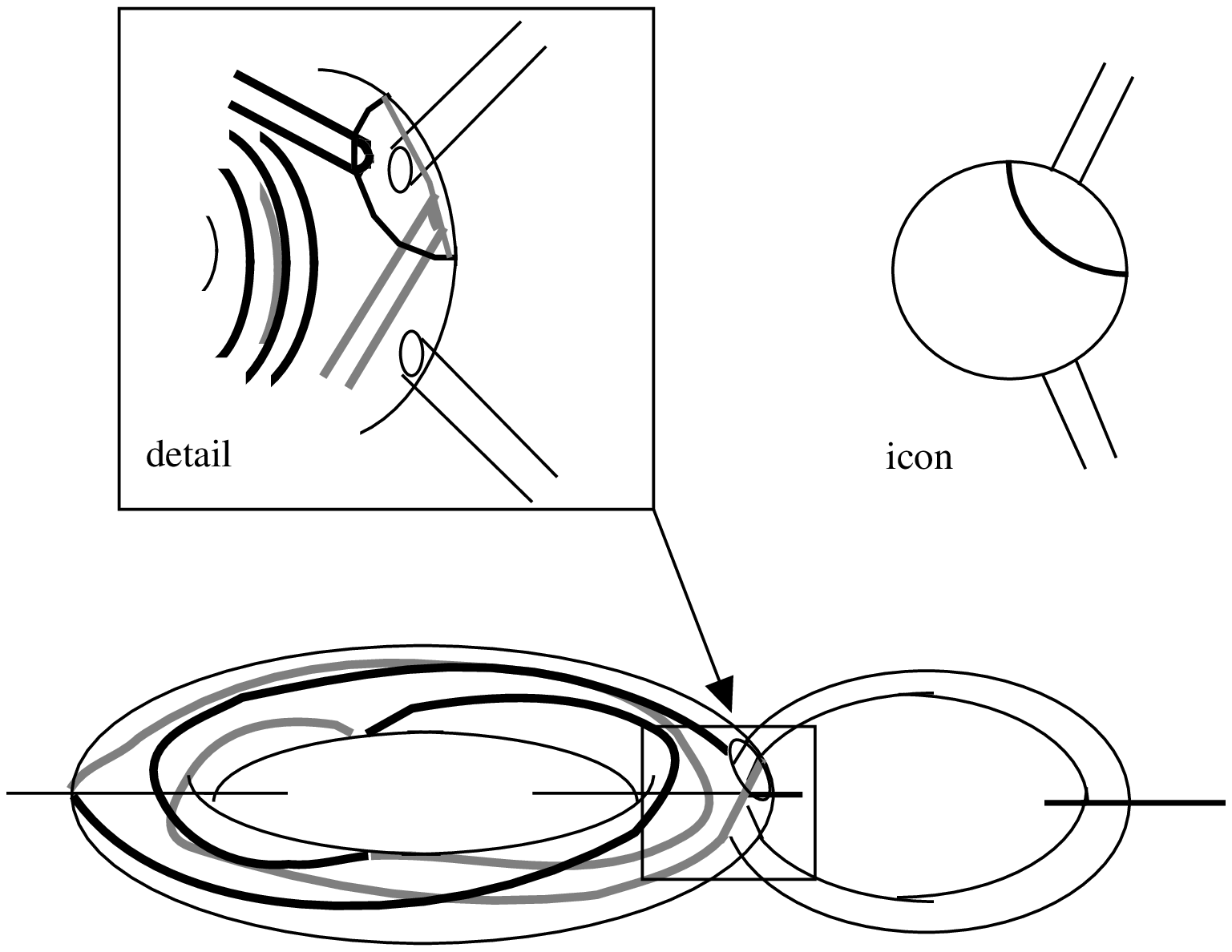}
\relabel {icon}{icon}
\relabel {detail}{detail}
\endrelabelbox}
\caption{}
\end{figure}


\begin{figure}[ht!] 
\centerline{\relabelbox\small
\epsfxsize=84mm \epsfbox{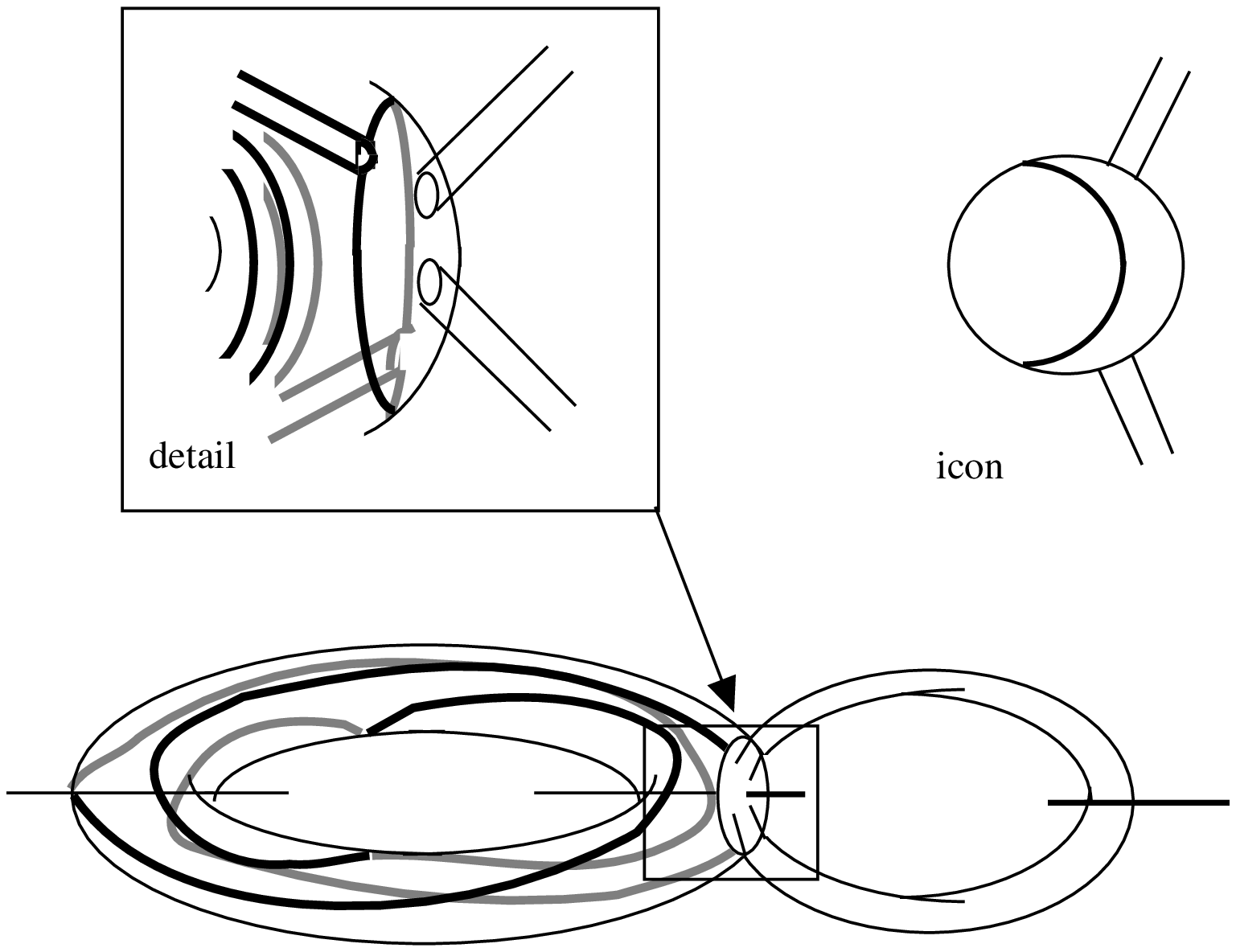}
\relabel {icon}{icon}
\relabela <-2pt,-2pt> {detail}{detail}
\endrelabelbox}
\caption{}
\end{figure}

\begin{lemma}
\label{annuli2}

Suppose $\Aaa$ is a properly imbedded essential collection of annuli in a 
genus two handlebody $H$.  Then the components of $\bdd \Aaa$ are either
all  twisted or all longitudes.  If they are all longitudes, then
the components of $\Aaa$ are all parallel and each is
non-separating in $H$.  If
they are all twisted  and $I(\bdd \Aaa) = (k, l, 0, 0)$ then one of
these two descriptions applies:

\begin{itemize}

\item  $\Aaa$ consists of two families of $k/2$ and $l/2$ parallel 
annuli, each annulus separates $H$ or 

\item $\Aaa$ consists of at most three families of parallel annuli, 
numbering respectively $e, f, g \geq 0$, with each annulus in the first
two  families non-separating, each annulus in the last family separating
and $e + f =  l, e + f + 2g = k$.  

\end{itemize}

\end{lemma}

\pf  By Lemma \ref{annuli} $I(\bdd \Aaa) = (k, l, 0, 0)$.  Let $A'$ 
denote the surface obtained from a $\bdd$--compression of $\Aaa$,
necessarily  into the unique complementary component of $\bdd H - \Aaa$
that is a $4$--punctured sphere (or twice punctured torus if $l = 0$). 
Then $A'$ contains an essential disk $D$, and $\bdd D$ is disjoint from
$\bdd  A$. 

If $D$ is a separating disk in $H$ then the complementary solid tori
contain $\Aaa$.  Any proper annulus in a solid torus is either
compressible or $\bdd$--parallel, so in each solid torus component $T$
of  $H - D$, $\Aaa$ is a collection of annuli all parallel to the
component of $\bdd T - \Aaa$ that contains $D$, and to no other
component of
$\bdd T - \Aaa$.   It follows that $\bdd \Aaa \cap T$ consists of a
collection of torus knots $L_{(p, q)}, p \geq 2$ in $\bdd T$. 

If $D$ is a non-separating disk, then $H - D$ is a single solid torus $T$ 
and all curves of $\Aaa \cap \bdd T$ are parallel in $\bdd T$. 
Each  annulus is $\bdd$--parallel to an annular component of $\bdd T -
\Aaa$ that  contains one of the two copies of $D$ lying in $\bdd T$.  If
the curves are all longitudes in $T$ (so each annulus in \Aaa\ is
$\bdd$--parallel to both annuli of $\bdd T - \Aaa$) then the annuli must
all be parallel, with a copy of
$D$ in each of the two components of
$\bdd T -
\Aaa$ to which they are boundary parallel.  If $\bdd T \cap \Aaa$
consists of  $(p, q), p \geq 2$ curves then each annulus in $\Aaa$  is
boundary parallel to exactly one annulus in $\bdd T$.  Since $\Aaa$ is
essential, such an annulus in $\bdd T$ must contain either one copy of 
$D$ or the other, or both copies of $D$.  This accounts for the three
families,  as described.  (See Figure 14.)
\qed


\begin{figure}[ht!] 
\centerline{\relabelbox\small
\epsfxsize=30mm \epsfbox{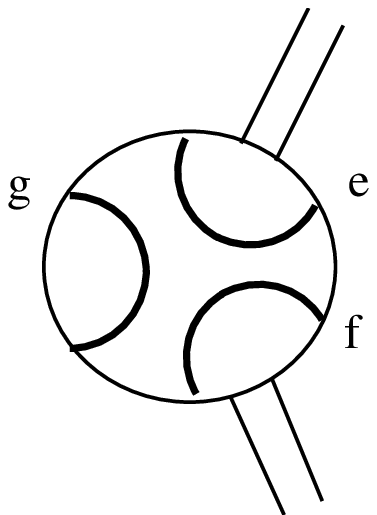}
\relabel {g}{$g$}
\relabel {e}{$e$}
\relabel {f}{$f$}
\endrelabelbox}
\caption{}
\end{figure}

\section{Canonical tori in Heegaard genus two manifolds}
\label{character}

For $M$ a closed orientable irreducible $3$--manifold there is a
(possibly empty) collection of tori, each of whose complementary
components is either a Seifert manifold or contains no essential tori
or annuli.  A minimal such collection \Fff\ of tori is called the
set of {\em canonical tori} for $M$ and is unique up to isotopy
\cite[Chapter IX]{Ja}.  
 
Suppose $M$ is of Heegaard genus two and contains an essential
torus.  Let $M = \aub$ be a (strongly irreducible) genus two Heegaard
splitting.  Using the sweep-out of $M$ by
$P$ determined by the Heegaard splitting, we can isotope \Fff\ so that
it intersects $A$ and $B$ in a collection of essential annuli. Indeed,
it is easy to arrange that all curves of $P \cap \Fff$ are essential
in both surfaces, so each component of $\Fff - P$ is an
incompressible annulus (cf \cite{RS}).  Inessential annuli in $A$ or
$B$ can be removed by an isotopy.  In the end, since no component of
\Fff\ can lie in a handlebody, $\Fff \cap A$ and $\Fff \cap B$ are
non-empty collections of essential annuli.

Note that if $T$ is a torus in \Fff\ and \aaa\ is an essential curve in
$T$, then on at least one side of $T$, \aaa\ cannot be the end of an
essential annulus. This is obvious if on one side of $T$ the component
of $M - T$ is acylindrical.  If, on the other hand, both sides are
Seifert manifolds, then the annuli must both be vertical, so the
fiberings of the Seifert manifolds agree on $T$.  This contradicts the
minimality of
\Fff. These remarks show that in $P$, if $I(P \cap \Fff) = (k, l,
0, 0)$ then $0 \leq l \leq k \leq 2$ (and of course $k + l$ is even).
With this in mind, we now examine how the tori \Fff\ can intersect $A$
and $B$.

Note that most of this section is covered by results in \cite{Ko}. 
Our perspective here is somewhat different though, as we are
interested in multiple splittings of the same manifold.  We include a
complete list of cases for future reference in later sections.  

\medskip
{\bf Case 1}\qua{\rm (Single annulus)}\qua   From \ref{annuli2} we see
that if $k = 2, l = 0$ then $\Fff$ is a separating annulus in each of
$A$ and $B$ and the Seifert manifold $V \subset M$ has base space a disk
and two singular fibers.  Since $P \cap V$ is a single annulus, call
this the {\em single annulus} case. 
Example \ref{symcable}, Variation 2 describes all splittings of this
type.  A special case is \ref{K4} Variation 3, when one of the Dehn
surgery curves is placed in $\mu_{a_{\pm}}$ and the other in either
$\mu_{b_{\pm}}$ or \sss. When one is in $\mu_{a_{\pm}}$ and the other
in  $\mu_{b_{\pm}}$, the annulus $P \cap V$ is in the part of $\bdd A
\cap \bdd \Ggg$ that's identified with  $\bdd B \cap \bdd \Ggg$. 
See Figure 15.


\begin{figure}[ht!] 
\centerline{\relabelbox\small
\epsfxsize=47mm \epsfbox{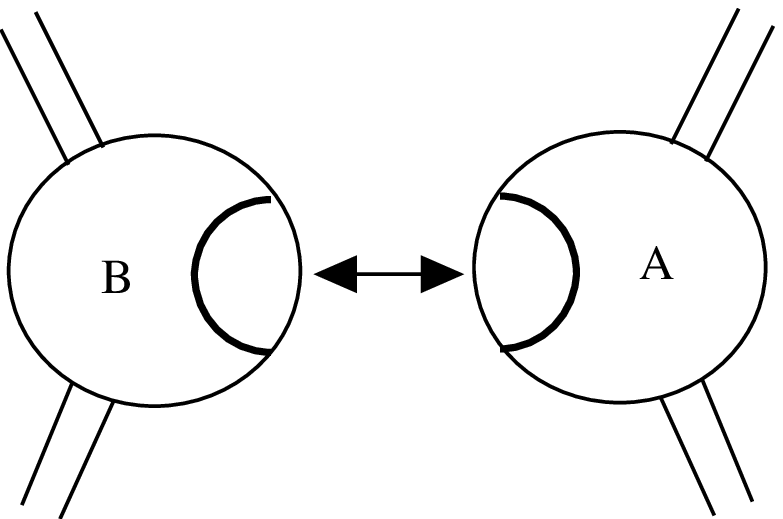}
\relabel {A}{$A$}
\relabel {B}{$B$}
\endrelabelbox}
\caption{}
\end{figure}

\medskip
{\bf Case 2}\qua (Non-separating torus)\qua If $k = l = 1$ then \Fff\
intersects both $A$ and $B$ in a single non-separating annulus.  Since
no properly imbedded annulus in $M - \Fff$ (with ends on the same side
of \Fff) is essential, the involution \Thp\ takes each annulus $\Fff
\cap A$ and $\Fff \cap B$ to itself.  This means that in each of $A$
and $B$ the curves $\Fff \cap P$ are longitudes.  Call this the {\em
single non-separating torus} case. Example \ref{nonsep} describes all
splittings of this type.  See Figure 16.


\begin{figure}[ht!] 
\centerline{\relabelbox\small
\epsfxsize=60mm \epsfbox{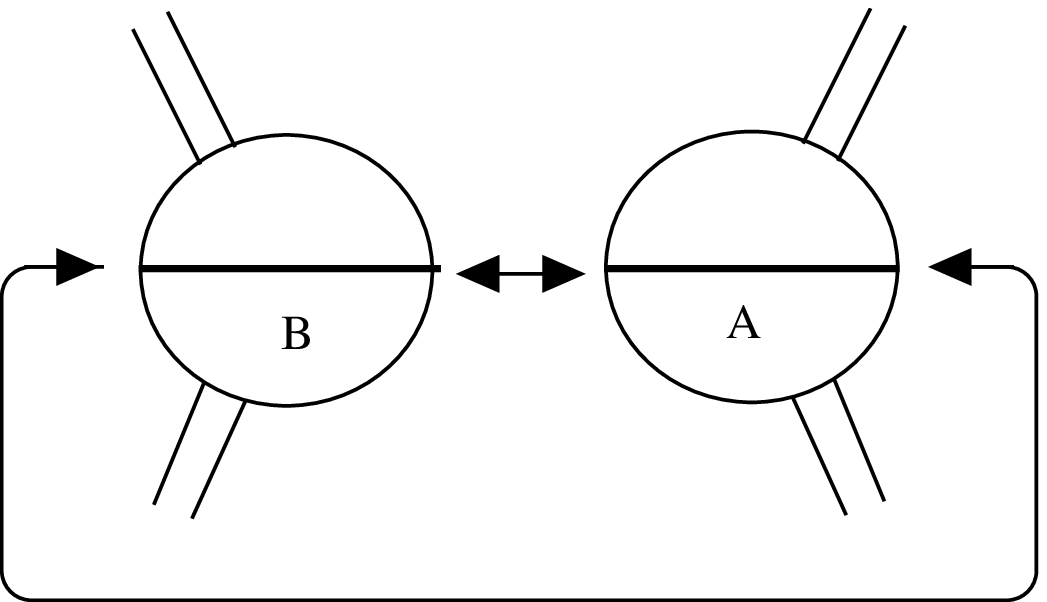}
\relabel {A}{$A$}
\relabel {B}{$B$}
\endrelabelbox}
\caption{}
\end{figure}

The case $k = l = 2$ admits a number of possibilities, depending on
whether the annuli in $A$ and/or $B$ are separating or non-separating
and, if non-separating, whether they are parallel or not.

\medskip
{\bf Case 3}\qua (Double torus)\qua  If $k = l = 2$ and annuli on both
sides are non-separating, then either $\Fff$ is a single
separating torus, for example cutting off the neighborhood of a
one-sided Klein bottle (discussed as Case 7 below), or $\Fff$ is a pair
of non-separating tori.  Between the tori lies a Seifert manifold with
base the annulus and one or two singular fibers (at most one in each of
$A$ and $B$).  Whether there are one or two singular fibers depends on
whether the annuli on one side or both sides are non-parallel.  Call
the latter the {\em double torus} case.  Example \ref{nonsep},
Variations 1 and 2 describe all splittings of this type.  See Figure
17.


\begin{figure}[ht!] 
\centerline{\relabelbox\small
\epsfxsize=60mm \epsfbox{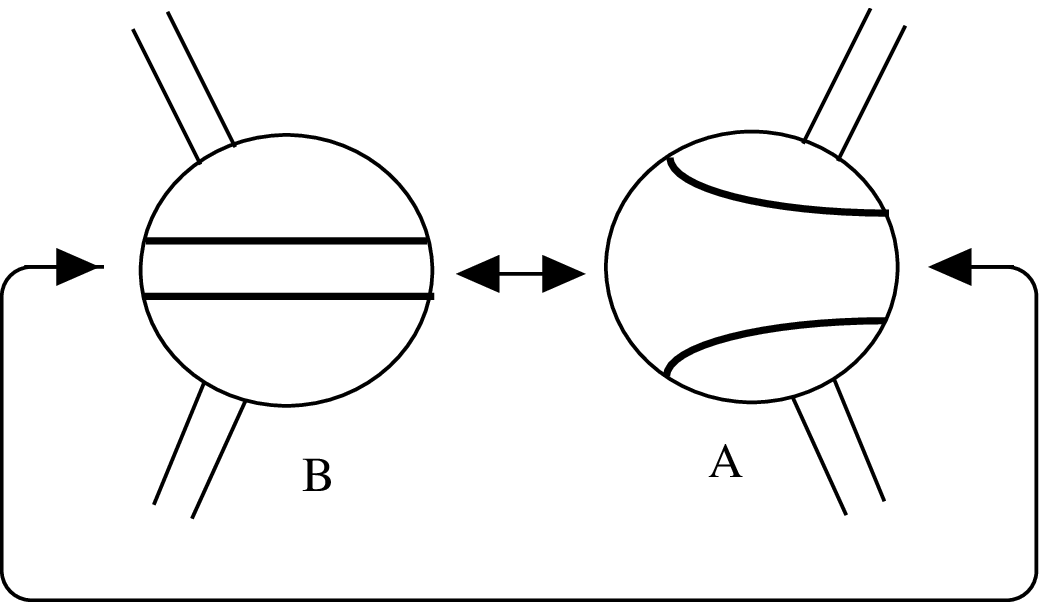}
\relabel {A}{$A$}
\relabel {B}{$B$}
\endrelabelbox}
\caption{}
\end{figure}

\medskip
{\bf Case 4}\qua (Double annulus)\qua Suppose $k = l = 2$ and in one of $A$ or
$B$, say $A$, the annuli are separating and in the other they are
non-separating and non-parallel.  Then $\Fff$ is a single separating
torus.  On one side of the torus is a Seifert manifold $V$ fibering
over the disk with three exceptional fibers, two in $A$ and
the third in $B$ lying between the pair of non-separating annuli $\Fff
\cap B$.  Call this the {\em double annulus} case. Example
\ref{K4}, Variation 4 describes all splittings of this type.  See
Figure 18.  The dotted half-circle indicates schematically that there
is an additional twisted annulus, not visible in this cross-section and
separated from the visible one by a separating disk in the
handlebody.  

\begin{figure}[ht!] 
\centerline{\relabelbox\small
\epsfxsize=60mm \epsfbox{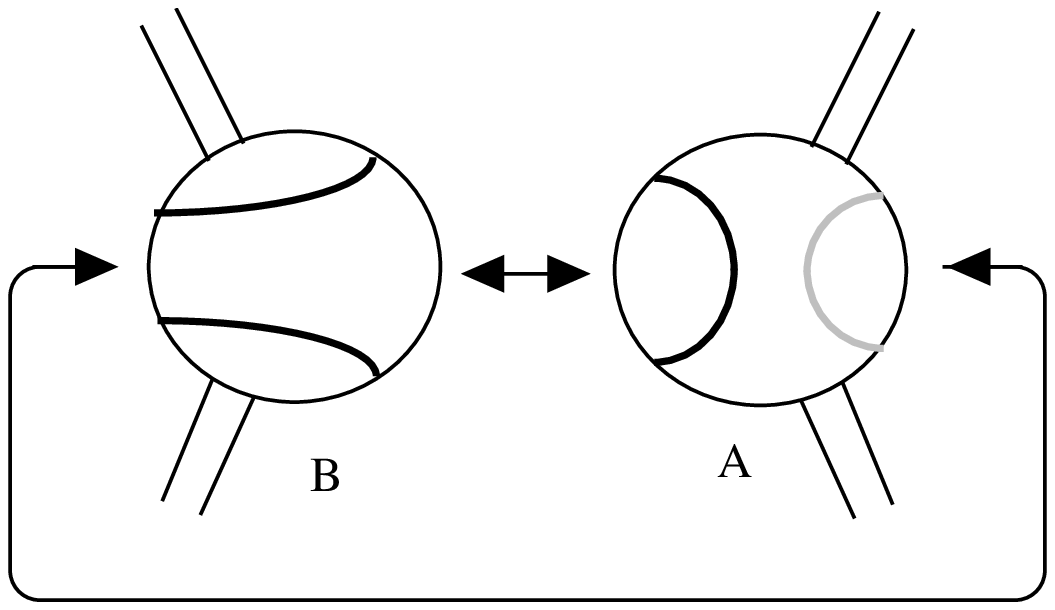}
\relabel {A}{$A$}
\relabel {B}{$B$}
\endrelabelbox}
\caption{}
\end{figure}

\medskip
{\bf Case 5}\qua (Parallel annuli)\qua Suppose $k = l = 2$ and in one of $A$ or
$B$, say $A$, the annuli are separating and in the other they are
non-separating and parallel. Then again $\Fff$ is a single separating
torus.  On one side of the torus is a Seifert manifold $V$ fibering
over the disk with two exceptional fibers, both in $A$. $\Fff
\cap B$ and $P \cap V$ are each a pair of parallel
annuli, so call it the {\em parallel annuli} case.  Notice that the
annuli $\Fff \cap B$ are longitudinal by the same argument
as in the single non-separating torus case.   Example \ref{K4},
Variation 3, with surgeries in $\mu_{a_+}$ and
$\mu_{a_-}$, describes all splittings of this type. See Figure 19.


\begin{figure}[ht!] 
\centerline{\relabelbox\small
\epsfxsize=60mm \epsfbox{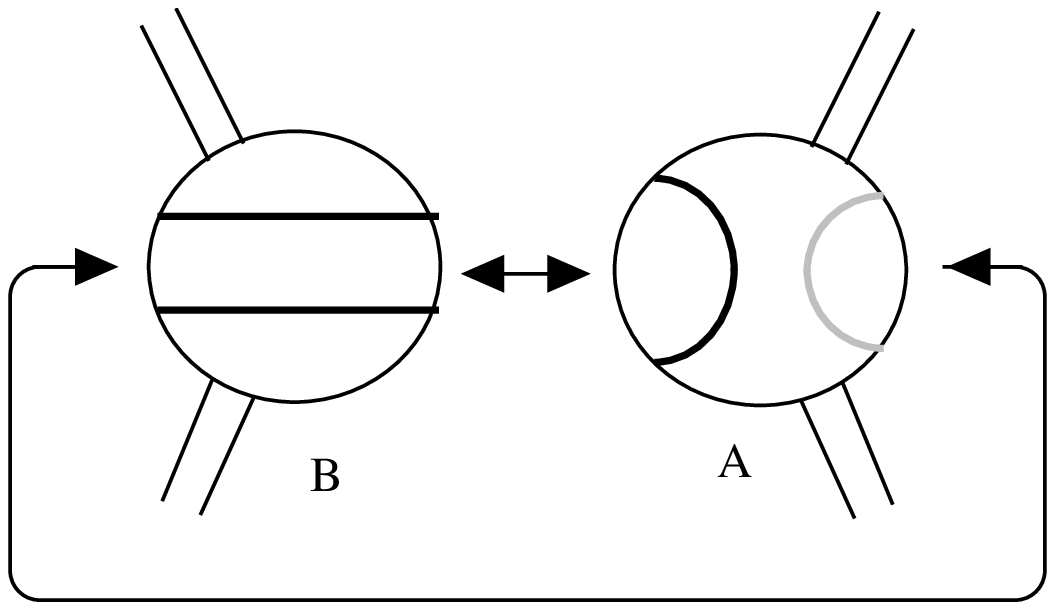}
\relabel {A}{$A$}
\relabel {B}{$B$}
\endrelabelbox}
\caption{}
\end{figure}

\medskip
{\bf Case 6}\qua (Non-parallel tori)\qua Suppose $k = l = 2$ and in both of
$A$ and $B$ the annuli are separating.  Then \Fff\ consists of two
separating tori, each bounding Seifert manifolds which fiber over the
disk with two exceptional fibers.  Example \ref{doubcable},
Variation 1, with Dehn surgery performed on all four of $\aaa_{na},
\aaa_{nb}, \aaa_{sa}, \aaa_{sb}$ describes all examples of this type.
See Figure 20.


\begin{figure}[ht!] 
\centerline{\relabelbox\small
\epsfxsize=60mm \epsfbox{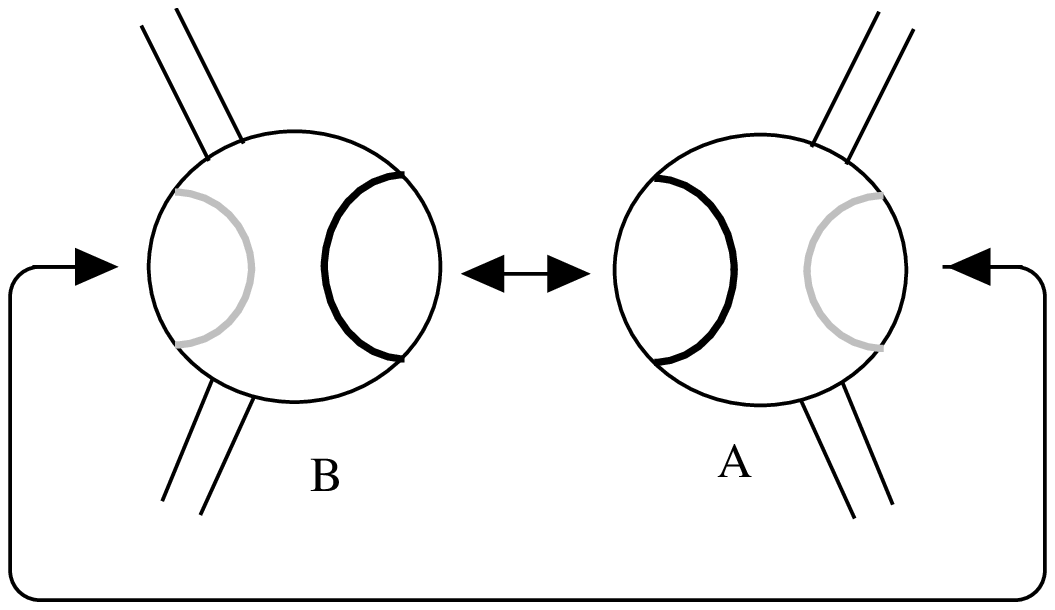}
\relabel {A}{$A$}
\relabel {B}{$B$}
\endrelabelbox}
\caption{}
\end{figure}

\medskip
{\bf Case 7}\qua (Klein bottle)\qua If $k = l = 2$ and annuli on both
sides are non-separating, then it could be that, when the pairs of
annuli are attached along their ends, the result is a single
separating torus.  The torus cuts off a Seifert piece $V$ that is the
union of the two parts of $A$ and $B$ that lie between the annuli. 
For example, if the pairs of annuli $\Fff \cap A$ and $\Fff \cap B$ are
both parallel in $A$ and $B$ respectively, then $V$ is
the neighborhood of a one-sided Klein bottle.  More generally, $V$ 
fibers over a M{\"o}bius band with zero, one, or two singular
fibers (at most one in each of $A$ and $B$).  Note that when there are
no singular fibers, so $V$ is the neighborhood of a one-sided Klein
bottle, then $V$ can also be fibered over the disk with two singular
fibers.  The fibering circles are orthogonal in $\bdd V$; in the
fibering over the M{\"o}bius band the fiber projects to a
curve in the Klein bottle whose complement is a cylinder and in the
fibering over a disk the fiber projects to a curve whose complement is
two M{\"o}bius bands. These cases correspond to Example \ref{K4}. See
Figure 21.


\begin{figure}[ht!] 
\centerline{\relabelbox\small
\epsfxsize=60mm \epsfbox{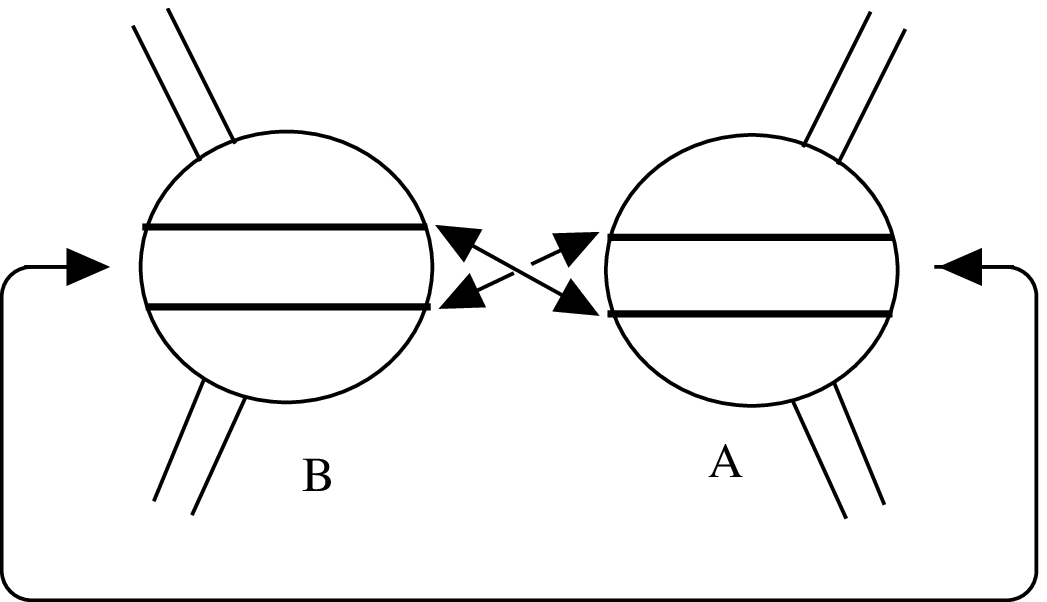}
\relabel {A}{$A$}
\relabel {B}{$B$}
\endrelabelbox}
\caption{}
\end{figure}

As can be seen from the above descriptions, each case is
determined, with one exception, by the Seifert piece $V$.  If $V$ is
just the neighborhood of a single non-separating torus, this is the
non-separating torus case.  If $V$ is a Seifert manifold over
the annulus with one or two exceptional fibers then it is the
double torus case.  If $V$ has two components (each fibering over the
disk with two exceptional fibers) then it is the non-parallel tori
case.  If $V$ fibers over
the disk with three exceptional fibers then it is the double annulus
case. If $V$ fibers over the M{\"o}bius band with one or two
exceptional fibers, then it is the Klein bottle case.  Only when
$V$ fibers over the disk with two exceptional fibers, could the
splitting be either the single annulus or the parallel annuli case
or (if both singular fibers have slope $ 1/2 $) the Klein
bottle case. 

In some situations the splittings described by the single annulus and
the parallel annuli case are closely related.  For example, begin with
Example \ref{K4} Variation 3, with one Dehn surgery circle in
each of $\mu_{a_{\pm}}$.  This is the parallel annulus case, with
canonical tori $\mu_{a_+} \cup \sss \cup \mu_{a_-} $.  Now move
the surgery circle in $\mu_{a_+}$ into $\sss$.  This is now the single
annulus case, with canonical tori $\sss \cup \mu_{a_-} $.  In
fact, if we cut along the annulus $\bdd \mu_{a_-} \cap \sss$, no
longer identifying boundaries of $A$ and $B$ there, the result is a
splitting as in Example \ref{symcable}, Variation 2.  See Figure 22.


\begin{figure}[ht!] 
\centerline{\relabelbox\small
\epsfxsize=70mm \epsfbox{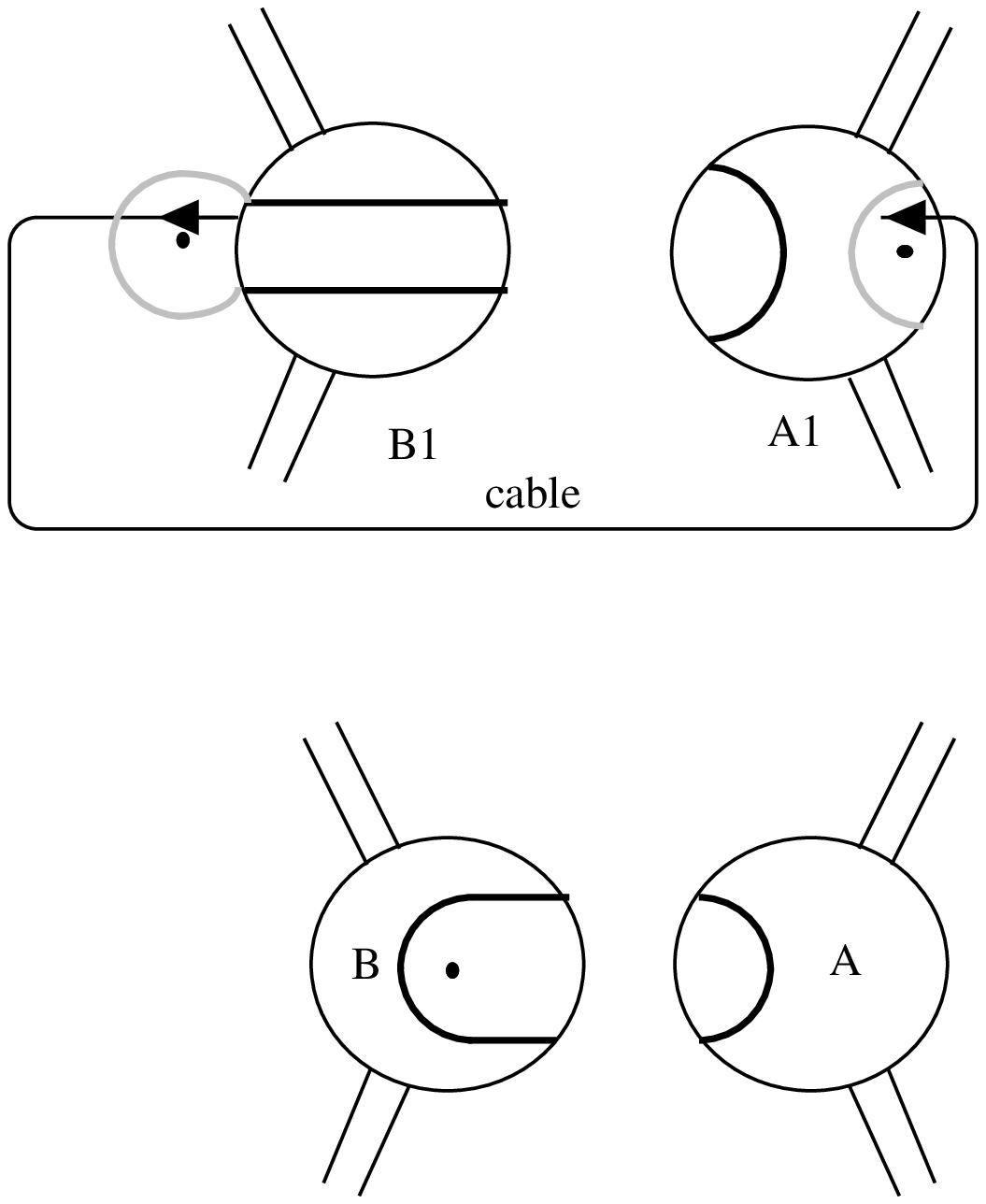}
\relabel {A}{$A$}
\relabela <-2pt,0pt> {B}{$B$}
\relabel {A1}{$A$}
\relabela <-2pt,0pt> {B1}{$B$}
\relabel {cable}{cable}
\endrelabelbox}
\caption{}
\end{figure}

We can formalize this example as follows:

\begin{lemma}
\label{switch}

Suppose $M = \aub$ has Seifert part $V$, fibering over the disk with
two exceptional fibers.  Suppose $P$ intersects $V$ as in the
parallel annuli case (that is, $P \cap V$ consists of two essential
parallel annuli) and the region between the annuli lies in $B$, say. 
Then $P$ can be cabled into $A$ to get a splitting surface $P'$
intersecting $V$ as in the single annulus case.  Moreover in $B' - V$
there is an annulus with one end a core of the annulus $\Fff \cap B'$
and other end a curve on $P'$ which is longitudinal in $A'$.

Dually, suppose $P$ intersects $V$ as in the single annulus case and in
$B - V$ there is an annulus with one end a core of $\Fff \cap B$ and 
other end a curve on $P$ which is longitudinal in $B$. Then
$P$ can be cabled into $B$ to get a splitting surface $P'$ intersecting
$V$ as in the parallel annulus case.

\end{lemma}

\pf  The first part is obtained by replacing one of the annuli in $P
\cap V$ with the incident annulus component of $\Fff \cap B$.  The
second part follows from the first by reversing the construction. 

In the example preceding the lemma, a spanning annulus as called for in
the lemma is one in \sss\ parallel to $\bdd \sss \cap
\mu_{b_{\pm}}$.   \qed

\medskip
In each of the seven cases listed above, there is a Seifert part $V$
(possibly just a thickened torus in the non-separating torus case)
which $P$ intersects in annuli and a single complementary component
$W$ which it intersects in a more complicated surface.  Since $\Fff$
lies in $V$ we know that $W$ is atoroidal.  It is also acylindrical
except possibly for an annulus whose ends in $\bdd V$ are non-fibered
curves and whose complement in $W$ is one or two solid tori. That is,
$W$ could itself be a Seifert manifold over a disk with two exceptional
fibers or over an annulus or M{\"o}bius band with one exceptional
fiber, as long as the fibering doesn't match the fibering of $V$.

In any case, $W$ has the following structure: $W = \awb$, where
$P_{-}$ is a properly imbedded surface (either a $4$--punctured sphere
or, exactly in the single annulus case, a twice punctured torus) and
$A_-, B_-$ are each genus two handlebodies.  $P_{-}$ lies in $\bdd A_-$
and $\bdd B_-$ as the complement of one or two longitudinal curves.  In
each case where this makes sense (ie, except in the single
non-separating torus case), $\bdd P_{-}$ is a fiber of the Seifert
manifold $V$ on the other side of \Fff.  $A_-$ (resp.\ $B_-$) can be
viewed as the mapping cylinders of maps from $P_{-}$ to a
$2$--complex $\Ss_A$ (resp.\ $\Ss_B$) consisting of one or two annuli in
\Fff\ and a single arc in $A_-$ (resp.\ $B_-$) with ends on the
annuli. Hence $W - \eta (\Ss_A \cup \Ss_B)$ is a product, restricting to
a product structure on the annuli $\bdd W - \eta (\Ss_A \cup \Ss_B)$. 
(Here $\eta$ denotes regular neighborhood.)  Hence it can be swept out
by $P_{-}$.

This sweep-out gives us some information about what sort of annuli
might be present in $W$.

\begin{lemma}
\label{squar}

Suppose $W$ contains an essential annulus with neither end parallel
to $\bdd P_-$.  Then $W$ contains an essential annulus or one-sided
M{\"o}bius band which intersects $P_-$ precisely in two parallel
spanning arcs.

\end{lemma}

\pf Consider how $P_-$ intersects the annulus \Aaa\ during the
sweep-out of $W$.   At the beginning it inevitably intersects \Aaa\ in
$\bdd$--compressing disks lying in $A_-$.  At the end it intersects
\Aaa\ in $\bdd$--compressing disks lying in $B_-$.  Nowhere can it
intersect it in both, so somewhere it intersects it in neither.  (The
details are standard and are suppressed.)  This means that the
intersection of $P_-$ with \Aaa\ consists entirely of spanning arcs of
\Aaa.  The squares into which \Aaa\ are cut by these arcs lie
alternately in $A_-$ and $B_-$.  It's easy to see that all these arcs
are parallel in $P_-$ so, we can assemble two of the squares into
which \Aaa\ is cut, one in $A$ and one in $B$, to produce an
annulus or one-sided M{\"o}bius band which intersects
$P_-$ precisely in two arcs.  \qed

\medskip
Let \Aaa\ be the annulus or one-sided M{\"o}bius band given by the
preceding lemma \ref{squar} and let $R_A, R_B$ be the squares in which
\Aaa\ intersects $A_-$ and $B_-$ respectively.  The complement of \Aaa\
in
$W$ is one or two solid tori, depending on whether $\Aaa \subset W$ is
separating or not.  Moreover the complement of $R_A$ in $A_-$ is also
one or two solid tori depending on whether $R_A \subset A_-$ is
separating or not, and similarly for $B_-$.  Similarly $P_- - \Aaa$ is
one or two annuli.  Since these annuli divide each solid torus of $W -
\Aaa$ into two solid tori, they are longitudinal annuli in the solid
torus. These facts give useful information about, for example, the index of
the singular fibers, but the crucial point here is that the
description is now sufficiently detailed that we have explicitly:

\begin{prop}
\label{splitchar}

Suppose $W$ contains an essential spanning  annulus with neither end
parallel to $\bdd P_-$.  Suppose the annulus is unique up to proper
isotopy and  \Aaa\ is the annulus or one-sided M{\"o}bius band given by
Lemma \ref{squar}.  Then there is an \Aaa--preserving involution
$\Theta_W$ of $W$, defined independently of  $P$ and a proper isotopy
of $P_-$ in $W$ so that after the isotopy
$\Theta_W | P_- = \Thp | P_-$.

\end{prop}

\pf  The proof is left as an exercise.  The fixed point set of 
$\Theta_W$ intersects \Aaa\ either  

\begin{itemize}

\item in two points, the centers of each of $R_A$ and $R_B$ or

\item  in two proper arcs orthogonal to the core of \Aaa\ or 

\item in the core of \Aaa,

\end{itemize}

depending on the structure of
$W$. \qed

\medskip
Proposition \ref{splitchar} is phrased to require a possible isotopy
of $P_-$ rather than of $\Aaa$, since in it application we will be
isotoping two different Heegaard splittings, using $\Aaa$ as a
reference annulus.  Also, the proper isotopy of $P_-$ in $W$ is not
necessarily fixed on $\bdd P_-$, so in fact $\Thp$ should be regarded
as $\Theta_W$ composed with some Dehn twist along a component (or two)
of $\bdd W$.

\section{Longitudes in genus $2$ handlebodies---some
technical lemmas}

We will need some technical lemmas which detect and place longitudes
in a genus two handlebody.

\begin{defn}

Two curves $\lll, \lll' \subset \bdd H$ on a genus two
handlebody $H$ are {\em separated} if they lie on opposite sides of a
separating disk in $H$. Two curves $\lll, \lll' \subset \bdd H$ are
{\em coannular} if they constitute the boundary of a properly imbedded
annulus in $H$.

\end{defn}

\begin{lemma}
\label{tech1}

Suppose $H$ is a genus two handlebody and the disjoint curves $c_1, c_2,
c_3 \subset \bdd H$ divide $\bdd H$ into two pairs of pants. Suppose
that
$c_1, c_2 \subset \bdd H$ are nonmeridinal curves which are coannular
in $H$.  Then 
$c_3$ is either meridinal or it intersects every meridian disk.  

\end{lemma}

\pf  Suppose there were a meridian disk $D$ disjoint from $c_3$ and
consider how $D$ intersects the annulus $\Aaa \subset H$ whose
boundary is $c_1 \cup c_2$.  Assume $|D \cap \Aaa|$ has been
minimized.  If
$D \cap A = \emptyset$ then $D$ is a separating disk in the handlebody
$H' = H - \eta(A)$.  Then $\bdd D$ divides $\bdd H - \bdd \Aaa$ into two
pairs of pants, so any essential curve in the complement of $\bdd D
\cup \bdd \Aaa$, eg $c_3$ is parallel either to $\bdd D$ or a
component of $\bdd \Aaa$.  But the latter violates the hypothesis.

If $D \cap A \neq \emptyset$ then consider an outermost arc of
intersection in $D$.  It cuts off a meridian disk $E$ of $H'$ that
is disjoint from $c_3$.  Two copies of $E$ banded together along the
core of $\Aaa$ in $\bdd H'$ gives a separating disk disjoint from
$c_3$.  This reduces the proof to the previous case. \qed 

\begin{lemma}
\label{tech2}

Suppose $H$ is a genus two handlebody and the disjoint curves $c_1, c_2,
c_3 \subset \bdd H$ divide $\bdd H$ into two pairs of pants. Suppose
that $c_1, c_2 \subset \bdd H$ are nonmeridinal curves which
are coannular in $H$.  Let $A$ be an
annulus, with ends denoted $\bdd_{\pm} A$, and attach $A \times I$ to
$H$ by identifying $\bdd_+ A \times I$ to a collar of $c_2$ and $\bdd_-
A \times I$ to a collar of $c_3$.  Then the resulting manifold $H'$ is
not a genus two handlebody.  

\end{lemma}

\pf  If $H'$ were a genus two handlebody, then the dual annulus
$A' = (core(A) \times I) \subset (A \times I)$ would be a non-separating
annulus in $H'$. This means that in the handlebody ($H$ again) obtained
by cutting open along $A'$, both $c_2$ and $c_3$ would be twisted or
longitudinal, but in any case each would be disjoint from some meridian
disk.  But in the case of $c_3$ this would violate Lemma \ref{tech1}.
\qed 

\begin{lemma}
\label{tech0}

Suppose $H$ is a genus two handlebody and the disjoint curves $c_1, c_2,
c_3 \subset \bdd H$ are nonmeridinal curves that divide $\bdd H$ into
two pairs of pants,
$V$ and
$V'$, with $\bdd V = \bdd V' = c_1 \cup c_2 \cup c_3$. Suppose that
$c_1, c_3 \subset \bdd H$ are separated curves
and that
$c_2$ is disjoint from some meridian disk.  Then one of
$c_1$ or $c_3$ is a longitude, and there is a disk $D$ which separates
$c_1$ and $c_3$ so that $|\bdd D \cap c_2| = 2$.

In particular, if both $c_1$ and $c_3$ are longitudes, then $H \cong V
\times I$.

\end{lemma}

\pf  Let $\Ddd$ be the union of three disjoint disks: a disk $D$ that
separates $c_1$ and $c_3$, and disjoint meridian disks $D_1$ and
$D_3$ which intersect $c_1$ and $c_3$ respectively.  Choose this
collection and a meridian disk
$D_2 \subset H$ whose boundary is disjoint from $c_2$,
so that, among all such disk collections, 
$|\Ddd \cap D_2|$ is minimal.  We can assume that $c_2$ intersects
each disk of $\Ddd$, since $c_2$ is not parallel to either $c_1$ or
$c_3$.  Hence $D_2$ is not parallel to any disk in $\Ddd$, so in fact
 $\Ddd \cap D_2 \neq \emptyset$. 

 By minimality of
$|\Ddd \cap D_2|$ all components of intersection are proper arcs in
$\Ddd$.  Consider an arc $\bbb$ of $\Ddd \cap D_2$ which is outermost in
$D_2$.  Simple counting arguments show that $\bbb \subset D$, that the
subdisk of $D_2$ cut off by \bbb\ intersects $c_1$ or $c_3$ (say
$c_1$) in a single point (for the arc is disjoint from $c_2$).  In
particular, $c_1$ is a longitude.  Even more, it follows then that
as many points of intersection with $c_2$ lie on one side of
\bbb\ in $D$ as on the other. Since this is true for any outermost
arc, it follows that all outermost arcs of $\Ddd \cap D_2$ in
$D_2$ are parallel to \bbb\ in $D - c_2$.  Furthermore we may as
well assume that all outermost disks of $D_2$ cut off by these arcs lie
on the same side of
$D$, the side containing $c_1$,
since otherwise two could be assembled to give a third meridian disk
$D_4$ which would be disjoint from the disks $D_1, D_3$ and from
the longitude $c_2$ and which would intersect $c_1$ and $c_3$ exactly
once.  The proof would then follow immediately.  (See Figure 23.)


\begin{figure}[ht!] 
\centerline{\relabelbox\small
\epsfxsize=80mm \epsfbox{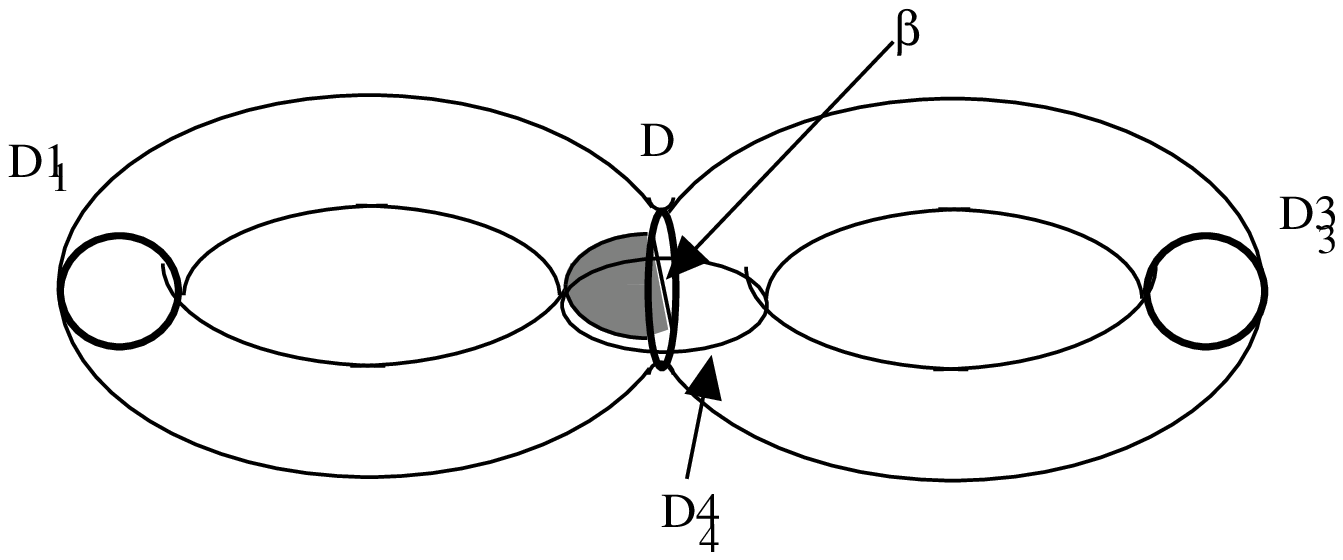}
\relabel {b}{$\beta$}
\relabela <-2pt,0pt> {D}{$D$}
\relabela <-2pt,0pt> {D1}{$D_1$}
\relabela <-2pt,0pt> {D3}{$D_3$}
\relabela <-2pt,0pt> {D4}{$D_4$}
\endrelabelbox}
\caption{}
\end{figure}

Now consider a disk component $E$ of $D_2 - \Ddd$ which is ``second to
outermost''.  That is, all but at most one arc of $\bdd E - \bdd H$ is
an outermost arc of intersection with $\Ddd$ in $D_2$.  To put it
another way, $\bdd E$ is a $2n$--gon, where every other side lies in
$\bdd H$, and of the $n$ remaining sides, at least $n-1$ are parallel
to \bbb\ in $D$.  The last side $s_n$ is perhaps an arc of $\Ddd \cap
D_2$.  (See Figure 24.)


\begin{figure}[ht!] 
\centerline{\relabelbox\small
\epsfxsize=35mm \epsfbox{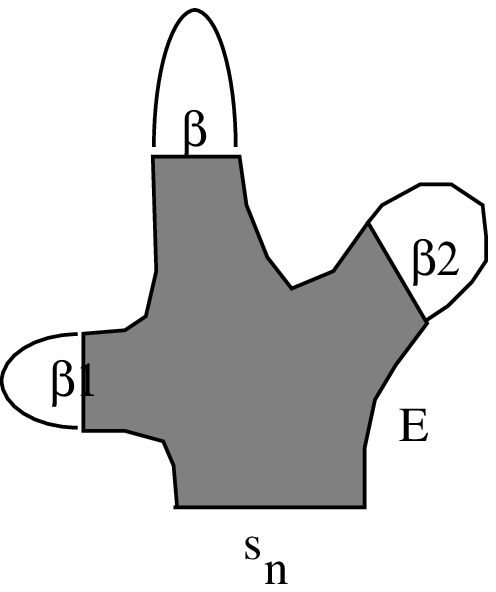}
\relabela <-1pt,0pt> {b}{$\beta$}
\relabela <-1pt,0pt> {b1}{$\beta$}
\relabela <-1pt,0pt> {b2}{$\beta$}
\relabel {s}{$s_n$}
\relabel {E}{$E$}
\endrelabelbox}
\caption{}
\end{figure}

The sides of $E$ that lie in $\bdd H$ and that have both ends on
$\bdd \bbb$ are easy to describe:  Since they are disjoint from $D_3$
and are essential in the pair of pants component of $\bdd H - \bdd
\Ddd$ on which they lie, each must cross $c_3$.  Moreover, since
they are disjoint from $c_2$, they can't cross $c_3$ more than once,
hence they cross $c_3$ exactly once.   Moreover, each must have its ends
at opposite ends of
\bbb, since if any had both ends at the same end of \bbb\ it would
follow that $c_2
\cap D_3 = \emptyset$ and that would force $c_2$ to be parallel to
$c_1$. But even one such arc of $\bdd E \cap \bdd H$, disjoint
from $c_2$ and $D_3$, crossing $c_3$ once and having ends at opposite
ends of $\bbb$, could be combined with an outermost disk of $D_2$ with
side at $\bbb$ to give a meridian disk $D_4$ as described before.  


\begin{figure}[ht!] 
\centerline{\relabelbox\small
\epsfxsize=90mm \epsfbox{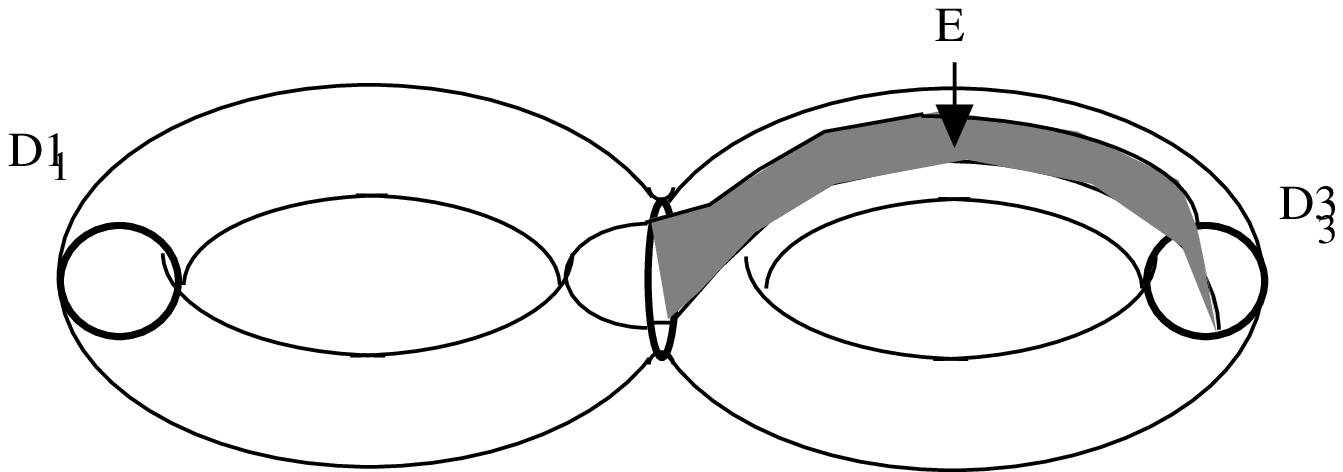}
\relabela <-1pt,0pt> {D1}{$D_1$}
\relabela <-1pt,0pt> {D3}{$D_3$}
\relabela <-2pt,0pt> {E}{$E$}
\endrelabelbox}
\caption{}
\end{figure}

So the only
remaining case to consider is
$n = 2$, with $s_n$ not parallel to $\bbb$ in $D$. So $E$ is a square,
with one side parallel to $\bbb$ and the opposite side, $s_2$, an arc
lying in
$D_3$ or $D$. (See Figure 25.)  Then simple
combinatorial arguments in the pair of pants bounded by $D$ and $D_3$
show that $s_2 \subset D_3$, since otherwise $s_2$ and $\bbb$ would
cross in $D$.  With $s_2 \subset D_3$, a simple counting argument shows
that $s_2$ cuts off from  $D_3$ a disk which can be made disjoint from
$c_3$ and intersecting $c_2$ in a single point.  The union of this disk
and $E$ along $s_2$ gives a disk, parallel to a subdisk of $D$ cut off
by
\bbb\, that is disjoint from $c_1$ and $c_3$ and intersects $c_2$
exactly once.  It follows that $D$ intersects $c_2$ twice, as
required.   \qed

\begin{lemma}

Suppose $H$ is a genus two handlebody and the curves $c_1, c_2, c_3
\subset \bdd H$ divide $\bdd H$ into two pairs of pants, $V$ and
$V'$ with $\bdd V = \bdd V' = c_1 \cup c_2 \cup c_3$. Suppose that
$c_1, c_3 \subset \bdd H$ are separated non-meridinal curves and there
is a properly imbedded disk in $H$  which intersects $c_1 \cup c_2$ in a
single point. Then $c_3$ is a longitude, and there is a disk $D$ which
separates $c_1$ and $c_3$ so that $|\bdd D \cap c_2| = 2$.

In particular, if $c_1$ is also a longitude, then $H \cong V \times I$.

\end{lemma}

\pf  Suppose there is a disk $D'$ that is disjoint from $c_1$
and intersects
$c_2$ in a single point.  Then $c_2$ is a longitude and it follows
from Lemma \ref{tech0} that some  disk $D$ separating $c_1$ from
$c_3$ intersects
$c_2$ twice. Using outermost arcs of intersection in $D$, it's then
easy to modify the disk $D'$ so that is disjoint from $D$.  Then
$D'$ must be a meridian curve for the solid torus on the side of $D$
that contains $c_3$.  Since $D'$ intersects $c_2$ in one point, it
follows that it also intersects $c_3$ in one point, completing the
proof in this case.

Suppose there is a disk $D'$ that is disjoint from $c_2$ but intersects
$c_1$ in one point (so, in fact, $c_1$ is a
longitude). By Lemma
\ref{tech0} there is a disk
$D$ that separates $c_1$ and
$c_3$ and intersects $c_2$ twice.  Choose the pair of disks $D$
and  $D'$ so that, among all such disks, $|D \cap D'|$ is minimal.

Consider an outermost
disk $E$ cut off by $D$ in $D'$, so $E$ is disjoint from both $c_1$
and $c_2$. Then $E$ lies on the side of
$D$ containing $c_3$ (since it is disjoint from $c_1$) and must
intersect $c_3$ at most (hence exactly) once, since it is disjoint
from $c_2$. Thus $c_3$ is also a longitude.   \qed

\begin{cor}
\label{tech}

Suppose $H$ is a genus two handlebody and the curves $c_1, c_2, c_3$\break
$\subset \bdd H$ divide $\bdd H$ into two pairs of pants, $V$ and
$V'$ with $\bdd V = \bdd V' = c_1 \cup c_2 \cup c_3$. Suppose that
$c_1, c_3 \subset \bdd H$ are separated curves.  Let $A$ be an
annulus with ends $\bdd_{\pm} A$.  Attach $A \times I$ to $H$ by
identifying $\bdd_+ A \times I$ to a collar of $c_1$ and $\bdd_- A
\times I$ to a collar of $c_2$.  Suppose the resulting manifold is also
a genus two handlebody. Then $c_3$ is a longitude of $H$, and there 
is a disk $D \subset H$ which separates 
$c_1$ and
$c_3$ in $H$ and which intersects $c_2$ exactly twice.  Moreover, if
$c_1$ is also a longitude of $H$, then $H \cong V \times I$.

\end{cor}

\pf  Let 
$H'$ be the new handlebody, and consider the properly imbedded dual
annulus $core(A) \times I \subset H'$.  Since it's
$\bdd$--compressible in $H'$, it follows that there is a disk $D'$ in $H$ 
which intersects $c_1 \cup c_2$ in a single point.  The result
follows from the previous lemma.  \qed

\begin{cor}
\label{tech+}

Lemma \ref{tech} remains true if $A \times I$ is replaced by any solid
torus $S^1 \times D^2$, attached at  $c_1$ and $c_2$ along parallel,
essential, non-meridinal annuli in $\bdd (S^1 \times D^2)$.

\end{cor}

The proof is the same, using either attaching annulus at $c_1$ or
$c_2$ in place of 
$core(A) \times I$. \qed

\section{Positioning a pair of splittings---the hyperbolike case}
\label{positioning}

A closed, orientable, irreducible $3$--manifold $M$ is called {\em
hyperbolike} if it has infinite fundamental group and contains no
immersed essential torus.  In the next two sections we will show that
any two Heegaard splittings of the same hyperbolike $3$--manifold can be
described by some variation of one of the examples in Section
\ref{examples}. As a consequence, the standard involutions of the
manifold induced by the two splittings commute.

In this section we will isotope the splitting surfaces $P$ and $Q$ so
that they are transverse and so that the curves of intersection and the
pieces of the surfaces cut out by them are particularly informative.
In the next section, we will move the surfaces so that they are no
longer transverse, but rather coincide as completely as possible.

Especially in the latter context, it will be useful to be able to
refer easily to the pieces of one splitting surface that lie in the
interior of one of the other handlebodies.  

\begin{defn}

Suppose $M = \aub = \xuy$ are two Heegaard splittings of $M$.  Let
$P_X$ denote the closure of $P \cap (interior X)$.  (So if $P$ and
$Q$ are transverse, as will not often be true in later discussion, then
$P_X$ is just $P \cap X$.).  Similarly define
$P_Y$, $Q_A$ and $Q_B$.

\end{defn}

 We begin with a useful lemma.

\begin{lemma}
\label{str.irred} 

Suppose \xuy\ is a genus two Heegaard splitting
of a closed hyperbolike manifold $M$ and $D_X, D_Y$ are essential
properly imbedded disks in $X$ and $Y$ respectively.  Then $|\bdd D_X
\cap \bdd D_Y| \geq 3$.

\end{lemma}

\pf  If $|D_X \cap D_Y| = 1$ then \xuy\ is stabilized and so $M$ is
either a lens space, or $S^2 \times S^1$ or $S^3$, but in any case is
not hyperbolike.  Suppose $|D_X \cap D_Y| = 0$ (so \xuy\ is {\em weakly
reducible}).  If the boundaries of $D_X$ and $D_Y$ are parallel in
$Q$, or one of the boundaries is separating, then \xuy\ is
reducible.  This means that either $M$ is reducible (hence not
hyperbolike) or \xuy\ is stabilized, and we have just shown that
this is impossible.   If the boundaries of $D_X$ and $D_Y$ are
non-separating and non-parallel then the surface
 $S$ obtained from $Q$ by doing both compressions simultaneously
is a sphere.  Moreover $S$ contains a separating essential circle
of $Q$ which is compressible on both sides, so again \xuy\ is
reducible.  

Finally, suppose $|D_X \cap D_Y| = 2$.  Then the union of
collar neighborhoods $\eta (D_X)$ and $\eta (D_Y)$ of $D_X$ and $D_Y$
along their two squares of intersection is a solid torus $W$.  Denote by
$X_-$ (resp.\ $Y_-$) the solid torus or pair of tori obtained by
compressing $X$ along $D_X$ (resp.\ $Y$ along $D_Y$).  Then $M$ is the
union of $X_-, Y_-$ and $W$, and the annuli of attachment of $W$ to
$X_-$ and $Y_-$ are either longitudinal (if the two points of
intersection of $\bdd D_X$ and $\bdd D_Y$ have opposite orientation)
or of slope $(1, 2)$ in $W$  (if the two points of
intersection of $\bdd D_X$ and $\bdd D_Y$ have the same orientation).  

So $M$ is the union of solid tori along essential annuli in their
boundary.  It is therefore either reducible or a Seifert
manifold.  In  any case it is not hyperbolike. \qed

\medskip
Suppose \xuy\ and \aub\ are two genus two Heegaard splitting
of a closed hyperbolike manifold $M$.   The two splittings define
generic sweep-outs of $M$, as described in \cite{RS}.  The pair of
sweep-outs is paramerized by points in $I \times I$.  Points in $I
\times I$ corresponding to positions where $P$ and $Q$ are not
transverse constitute a subcomplex of $I \times I$ called the {\em
graphic}.  Complementary components are called {\em regions}. 

Since the
surfaces involved have low genus, we can obtain useful information
about their relative positioning even if we allow a more liberal rule
than in \cite{RS} for labelling regions (that is, positionings in
which $P$ and $Q$ are transverse).  We label a region $A$ (resp.\ $B$)
if there is a meridian disk $D$ for $A$ (resp.\ $B$) such that $\bdd D
\subset (P - Q)$.  Labels $A$ and $B$ will be called {\em
$P$--labels}.  Similarly, we label a region $X$ (resp.\ $Y$) if there is
a meridian disk for $X$ (resp.\ $Y$) whose boundary is disjoint from
$P$. These labels are called $Q$--labels.  

Suppose that $C$ is the collection of curves of $P \cap Q$
that are essential in $P$ (resp.\ $Q$).  Then $C$ divides $P$ (resp.\
$Q$) into two parts, one lying (except for some inessential parts) in
$X$ and one in $Y$ (resp.\ $A$ and $B$).  If the two parts of $P$ (resp.\
$Q$) have even Euler characteristic we say that the positioning is {\em
$P$--even} (resp.\ $Q$--even).  If the two have odd Euler characteristic
we say it's {\em $P$--odd} (resp.\ $Q$--odd).  

\begin{lemma}
\label{odds}

If a region is $P$--odd then its $P$--labels are a
subset of the $P$--labels of any adjacent region.  Similarly
for $Q$--odd regions and $Q$--labels.

\end{lemma}

\pf  By construction, we are ignoring curves in $P \cap Q$ which
are inessential in $P$, so no component of $P - C$ is a disk.  If the
region is $P$--odd, then it follows that both parts have Euler
characteristic $-1$.  This implies that, if there is a meridian for
$A$ disjoint from $Q$ then in fact some curve in $C$ is a meridian of
$A$. Since $C$ can be pushed into \px\ or \py, this means that both
\px\ and \py\ contain a meridian of $A$.  

The effect of moving to an
adjacent region in the complement of the graphic is to alter \px\ and
\py\ by adding a band (or a disk) to one and removing it from the
other.  Clearly adding a band (or disk) doesn't destroy a curve, such
as the meridian, so one copy of the meridian of $A$ persists in at
least one of \px\ or \py\ in the new region. \qed 

\begin{lemma}
\label{evens}

If there are adjacent regions which are both $P$--even (resp.\
$Q$--even) then the $P$--labels (resp.\ $Q$--labels) of one are a subset of
the $P$--labels (resp.\ $Q$--labels) of the other.  If there are adjacent
regions which are each both $P$--even and $Q$--even then the set of all
labels for one of the regions is a subset of the labels for the other.

\end{lemma}

\pf  Suppose two adjacent regions are both $P$--even.  Moving from one
region to the other may represent moving across a center tangency, which
clearly has no effect on labels, or moving across a saddle tangency. 
The latter changes the Euler characteristic of \px\ and \py\ by $\pm
1$, so if the parity determined by $C$ doesn't change, the saddle move
must have created or destroyed an inessential curve of $P \cap Q$. 
This means that one or both ends of the band that is exchanged from
\px\ to \py\, or vice versa, lies on an inessential curve of $P$.  If
one end lies on an inessential curve, then the move is effectively an
isotopy of $C$ and so has no effect on the labelling.  If both ends lie
on the same inessential curve the effect is to add two parallel,
possibly essential, curves to $C$.  This won't add a label $A$ or $B$,
since a meridian lying in the annulus created in \py\ previously lay
in \px, but it might destroy some other meridian in \px, so a label
might be deleted.  This is the only way in which $A$ and $B$ labels
could change.  To summarize: if there is a change in $A$ or $B$ label
it's to delete a label moving from the first region to the second, and
this only happens if the corresponding band has both ends on the same
curve of $P \cap Q$, and that curve is inessential in $P$.

Now consider the situation in $Q$ if the adjacent regions
are also both $Q$--even.  Moving from the second region to the first we
have already seen that the band that's attached will have its ends on
two different curves (the two created in moving from the first to the
second region).  So no $X$ or $Y$ label can disappear.  It follows that
the set of labels for the second region is a subset of the set of
labels for the first region. \qed

\begin{lemma}
\label{different}

Any region that is $P$--even and $Q$--odd (or vice versa) has a label
that is also a label of every adjacent region.  

\end{lemma}

\pf  It's easy to see that $\chi(\px)$ and
$\chi(\qa)$ have the same parity:  For example, the sum of their
parities is the parity of the orientable surface created by doing a
double-curve sum of the two surfaces.  Furthermore, removing
curves of $P \cap Q$ that are inessential in both $P$ and $Q$ does not
alter the parity match.   So if a region is $P$--even and $Q$--odd it
follows that at least one curve in $P \cap Q$ is essential in $P$ and
inessential in $Q$ (or vice versa), ie, is a meridian $\mu$ of $A$ or
$B$ (or $X$ or $Y$).  When passing to an adjacent region in the
complement of the graphic, a band is added to either \px\ or \py, say
the former.  Before passing to the new region, move $\mu$ slightly into
\px.  Then $\mu$ will still lie in \px\ after moving to the adjacent
region.  \qed

\begin{lemma}
\label{neighbors}

If two adjacent regions have labels $A$ and $B$ then one of them has
both labels $A$ and $B$.  

\end{lemma}

\pf  This follows from \ref{odds} if either region is $P$--odd and
from \ref{evens} if both regions are $P$--even. \qed

\begin{lemma}
\label{both}

No region has both labels $A$ and $B$.

\end{lemma}

\pf  The proof is a recapitulation of ideas in \cite{RS} and
\cite{JS}.

Suppose a region has both labels.  The meridians are
unaffected  by removing, by an isotopy, all simple closed curves in
$P
\cap Q$ which are inessential in both $P$ and $Q$. The meridians 
of $A$ and
$B$ which account for the labels must intersect, by \ref{str.irred}, so
they cannot be on opposite sides of $Q$.  If any curve of $P \cap Q$ is
essential in $P$ and inessential  in $Q$ then it is a meridian of $A$,
say, that can be pushed to lie on the opposite side of $Q$ from the
meridian of $B$, a contradiction.  So every curve in $P \cap Q$ is
essential in $Q$.  

Say the meridians of $A$ and $B$ that are disjoint from $Q$ both
lie in $X$.  If any component of $P_X$ is a $\bdd$--parallel
annulus, push it across into $Y$---this has no effect on the
labelling.  If possible, $\bdd$--reduce $X$ in the complement of $P_X$. 
We will assume that no such $\bdd$--reductions are possible, so $X$
remains a genus two handlebody---the argument is easier if
$\bdd$--reductions can be done.  This guarantees that no component of
$P_X$ is a meridian disk of $X$ so every curve in $P \cap Q$ is also
essential in $P$.  

Then the boundary of any meridian disk of $X$ must intersect $\bdd
P$, since $P$ is strongly irreducible (\ref{str.irred}).  In
particular,  no curve of $P \cap Q$ is a meridian curve for $X$, nor
can $P$ lie  entirely inside of $X$.  

Since $P_X$ is compressible yet no boundary component is a meridian of
$X$ it follows that $\chi (P_X) = -2$ and a compression of $P_X$
creates a set of incompressible annuli.  Since all curves of $P
\cap Q$ are essential in both surfaces, one of $Q_A$ or $Q_B$, say the
former, has $\chi (Q_A) = -2$.  Let  $\Aaa$ be the incompressible annuli
in $X$ obtained by compressing $P_X$ into $A$. Dually, $P_X$ is
obtained from $\Aaa$ by attaching a tube along an arc \bbb\ dual to the
compression disk.  It follows from \cite[Theorem 2.1]{JS} that
there is a meridian disk $D'$ for $X$, isotoped to minimize $|D' \cap
\Aaa|$, so that the arc \bbb\ lies in $D'$. 

Consider how a distant meridian disk $D$ of $X$ intersects $\Aaa$ and
how it intersects a compressing disk $E$ for $B$.  First consider an
outermost arc \aaa\ of $\Aaa \cap D$.  Suppose $\bdd E$ is disjoint
from \aaa\ (as we can assume is true if the disk cut off by \aaa\ in
$D$ lies in $B$).  $\bdd$--compress $\Aaa$ to $Q$ via the disk cut off
by \aaa.  This changes an annulus of $\Aaa$ to a disk $\Ddd$.  If the
tube along \bbb\ were attached to
$\Ddd$ it would violate strong irreducibility of $P$ (since $\bdd
\Ddd$ is a meridian disk for $A$ disjoint from $E$), so
$\Ddd
\subset P_X$.  If $\Ddd$ is not $\bdd$ parallel it is parallel to a
meridian disk for $X$ disjoint from $P$ and if it is $\bdd$ parallel
then the original annulus was a
$\bdd$--parallel annulus in $P_X$.  Either is a contradiction.  So we can
assume that each outermost arc of $\Aaa$ in $D$ intersects $\bdd E$
and that the disk in $D$ cut off by the outermost arc lies in $A$.  

This means that there is a disk in $B \cap D$ all but at most one of
whose boundary arcs in $\Aaa$ are outermost arcs, and each of these
intersects $\bdd E$.  It is now easy to argue (see \cite{JS} for
details) that in fact \bbb\ is isotopic in $B \cap X$ to an arc of $E
\cap D$ which connects two adjacent outermost arcs, ie, \bbb\ is
parallel in $B \cap X$ to a spanning arc of one of the annuli of
$Q_B$.  But this implies that there is a meridian disk for $B$ (the
complement of the tube $\eta (\bbb)$ in the annulus of $Q_B$ of which
\bbb\ has been made a spanning arc) that intersects the compressing disk
for
$A$ dual to \bbb\ in two points.  This contradicts \ref{str.irred}.
\qed

\begin{lemma}

There is an unlabelled region.

\end{lemma}

\pf  The argument is a variant of that in \cite{RS}.  Combining
Lemmas \ref{neighbors} and \ref{both} we see that adjacent regions can't
have labels $A$ and $B$ or labels $X$ and $Y$.  So either there is an
unlabelled region or there is a vertex whose four adjacent regions are
each labelled with one label, appearing in order around the vertex $A$,
$X$, $B$, $Y$.  Then no region is $P$--odd and $Q$--even or vice versa,
by Lemma \ref{different}.  By Lemma \ref{odds} the regions labelled $A$
and $B$ must be $P$--even and those labelled $X$ and $Y$ must be
$Q$--even, so in fact all must be both $P$--even and $Q$--even.  But this
would contradict Lemma \ref{evens}.  \qed 

\begin{thm}
\label{positioned} 

Suppose \xuy\ and \aub\ are two genus two Heegaard splittings
of a closed hyperbolike manifold $M$. Then $P$ and $Q$ can be isotoped
in $M$ so that each curve in $P \cap Q$ is essential in both $P$ and
$Q$, so that $\chi(\px) = \chi (\py) = \chi(\qa) = 
\chi (\qb) = -1$ and so that \px\ (resp.\ \py, \qa, \qb) is
incompressible in $X$ (resp.\ $Y$, $A$, $B$).

\end{thm}

\pf  Consider the positioning of $P$ and $Q$ represented by an
unlabelled region.  Curves of intersection that are inessential in
both surfaces can be removed by an isotopy without introducing
meridians in \px, \py, \qa, or \qb, ie, without altering the fact that
the configuration is unlabelled.  Then all curves of intersection must
be essential in both surfaces, for otherwise at least one such
curve would be a meridian.  If the configuration is $P$--odd (hence
also $Q$--odd) then we are done.  

So suppose the configuration is $P$--even (hence $Q$--even).  With no
loss of generality, assume $\chi (\px) = \chi (\qa) = -2$.  Since $X$ is
a handlebody and the region is unlabelled, \px\ is
$\bdd$--compressible.  Do a
$\bdd$--compression.  If the boundary compression is on an annulus
component $\Aaa_P$ of \px\ then the result is a disk.  It can't be a
meridian disk for $X$, by assumption, so $\Aaa_P$ must be
$\bdd$--parallel in $X$.  Push $\Aaa_P$ across the annulus $\Aaa_Q$ to
which it is parallel in $Q$.  Clearly this does not create a meridian
in either \px\ (only an annulus has been removed) or in \py\ (an
annulus has been attached to other annuli).  Similarly, if $\Aaa_Q
\subset \qa$ no meridian is created in \qa\ or
\qb.  

Suppose $\Aaa_Q \subset \qb$.  Then after the annuli are pushed
across each other, \qa\ is enlarged, so one might expect that it could
contain a meridian curve.  But note that if the meridian disk lay in
$X$ then, after compressing along it one would get a solid torus or
two, in which \px\ is incompressible.  But this is impossible since
$\chi (\px) = -2$.  Alternatively, if the meridian disk lay in $Y$ then
note that before the annulus is pushed across, the meridian curve
intersects only one end of each of the annuli in \py\ which have ends
on the ends of $\Aaa_Q$.  But in this case, an easy outermost
argument shows that there is a meridian curve in \qa\ before the
annulus is pushed across, another contradiction.  So we conclude
that nothing is lost by pushing such boundary parallel annuli in \px\
across $Q$. 

Eventually, after these parallel annuli are removed, \px\ is
$\bdd$--compressible along an arc lying in a non-annular component of
\px.   The component of \qa\ or \qb\ to which \px\ can be
$\bdd$--compressed is not an annulus, since \px\ contains no meridian
curves of $P$.  Do the $\bdd$--compression. The result is a positioning
of $P$ and $Q$ which is both $P$--odd and $Q$--odd and, essentially by
\ref{odds}, it remains unlabelled.  This configuration is as
required.  \qed

\section{Alignment of $P$ and $Q$}
\label{aligning}

\begin{lemma}
\label{bdd.compress}

Suppose \xuy\ and \aub\ are two genus two Heegaard splitting
of a closed manifold $M$, the surfaces \px, \py, \qa, and \qb\ are
incompressible in, respectively $X$, $Y$, $A$, and $B$. Then the surface
\px\ $\bdd$--compresses to one of \qa\ or \qb, and \py\
$\bdd$--compresses to the other.

\end{lemma}

\pf Each surface $\bdd$--compresses in the handlebody in which it
lies. With no loss of generality assume that \px\ $\bdd$--compresses to 
\qa.  Suppose it also $\bdd$--compresses to \qb.  Then since \py\ 
$\bdd$--compresses to one of \qa\ or \qb, we are done.  If \px\ fails to
$\bdd$--compress to  \qb\ then, symmetrically, \qb\ fails to
$\bdd$--compress to \px, so it must $\bdd$--compress to \py.  Hence \py\
$\bdd$--compresses to \qb. \qed

\begin{defn}

Suppose $P$ and $Q$ are closed surfaces in a $3$--manifold $M$ and $P$
(resp.\ $Q$) is the union of two subsurfaces $P_0$ and $P_+$ (resp.\
$Q_0$ and $Q_+$) along their common boundary curves. (That is, $P = P_0
\cup_{\bdd} P_+$ and similarly for $Q$).  Suppose finally that $P_0 =
Q_0$ whereas $P_+$ and $Q_+$ are transverse.  Then we say that $P$ and
$Q$ are {\em aligned along} $P_0 = Q_0$.

\end{defn}

\begin{lemma}
\label{align}

Suppose $M = \aub = \xuy$ are two genus
two Heegaard splittings of a hyperbolike closed $3$--manifold.  Then
the surfaces $P$ and $Q$ can be aligned along a subsurface $P_0 = Q_0$
with $\chi(P_0) = -2$ in such a way that each component of $\bdd
P_0 = \bdd Q_0$ is essential in all four handlebodies $A, B, X, Y$.

\end{lemma}

\pf  Following Theorem \ref{positioned}, isotope $P$ and $Q$ so that
each curve in $P \cap Q$ is essential in both $P$ and $Q$, so that
$\chi(\px) = \chi (\py) = \chi(\qa) =  \chi (\qb) = -1$ and so that
\px\ (resp.\ \py, \qa, \qb) is incompressible in $X$, (resp.\ $Y$, $A$,
$B$).  (Since $\chi(\px) = -1$, the last condition is equivalent to
saying that each curve in $\bdd P_X$ is essential in $X$.)  If an
annulus component of \px, say, is parallel to an annulus component of
\qa, say, then one can be pushed across the other without affecting
these hypotheses.  So we can assume that no component of any surface
\px, \py, \qa, and \qb\ is a $\bdd$--parallel annulus in the handlebody
in which it lies.  Then, moreover, any
$\bdd$--compression of \px, if it $\bdd$--compresses an annulus of \px,
would create a compressing disk for either \qa\ or \qb, contradicting
the hypothesis.  So we can assume that any $\bdd$--compression of any
surface is on the unique component whose Euler characteristic is $-1$.

Now apply Lemma \ref{bdd.compress} to find
a disk $D_{a,x}$ that $\bdd$--compresses \px\ to one of \qa\ or \qb, say
\qa, and a disk $D_{b,y}$ that $\bdd$--compresses  \py\ to \qb.  The
boundaries of the disks $D_{a,x} \subset A \cap X$ and $D_{b,y} \subset
B \cap Y$ lie on different surfaces so the disks can be made
disjoint.  

The curve $\bdd D_{a,x}$ is the union of two arcs, $\aaa \subset \px$
and $\bbb \subset \qa$. A collar of $D_{a,x}$ is a $3$--ball whose
boundary is the union of two disks $D_{\pm}$ parallel to $D_{a,x}$ and
a collar of each of \aaa\ and \bbb\ in \px\ and \qa\ respectively.
The $3$--ball can be used to define an isotopy of \px\ that replaces a
collar neighborhood of \aaa\ with the union of the disks $D_{\pm}$
and the collar of \bbb.  After this isotopy (and the kinking of
collars of the curve(s) of $P \cap Q$ on which the ends of \aaa\ lie),
$P$ and $Q$ will be aligned along an essential surface $P_0 = Q_0$
with $\chi (P_0) = -1$.  Repeat the process on $D_{b,y}$.  \qed

\begin{thm}

Suppose $M = \aub = \xuy$ are two genus
two Heegaard splittings of a hyperbolike closed $3$--manifold.   
Then the splittings are both some variation of one of the examples
of Section \ref{examples}.  In particular, \Thp\ and \Thq\ commute. 

\end{thm}

\pf  Following \ref{align}, we assume that 
the surfaces $P$ and $Q$ are aligned along a subsurface $P_0 = Q_0$
with $\chi(P_0) = -2$, and each component of $\bdd P_0 = \bdd Q_0$ is
essential in all four handlebodies $A, B, X, Y$.  We may further assume
that no component of $P_0$ is an annulus, for any such annulus could be
removed by a small isotopy of the surfaces, perhaps creating a curve of
transverse intersection.  We further assume that
$|P_+ \cap Q_+|$ has been minimized by isotopy rel
$\bdd P_0$.  Then $P_X, P_Y, Q_A$ and $Q_B$ all consist of
incompressible annuli.  Any of these annuli that is $\bdd$--parallel in
the handlebody in which it lies could be removed by an isotopy (possibly
adding it to $P_0$ or $Q_0$), so in fact all these annuli are
essential.  

According to \ref{annuli} the ends of $P_X$ in $Q$ can be isotoped to
lie parallel to (at most) two essential non-separating simple
closed curves in $Q$, and similarly for $P_Y$.  Since ends of $P_X$
and $P_Y$ can't cross, there are (at most) three 
simple closed curves $c_1, c_2, c_3$ in $Q$, decomposing $Q$ into two
pairs of pants, so that any component of $\bdd P_X \cup \bdd P_Y$ is
parallel to one of the three curves.  Moreover, if all three curves
$c_1, c_2, c_3$ are ends of annuli of $P - Q$ then the number of annuli
cannot be high, because of the following ``Rule of Three'': 

\begin{lemma}
\label{three}

If three or more ends of $P_X$ are parallel in $Q$ then at least one
must attach to an end of a curve in $P_Y$.  In particular, if $P_X$ has
three or more ends at each of two of the curves $c_i$, then all ends of
$P_Y$ must also be parallel to those two curves. 

\end{lemma}

\pf Immediate. (See Figure 26.)  \qed 


\begin{figure}[ht!] 
\centerline{\relabelbox\small
\epsfxsize=56mm \epsfbox{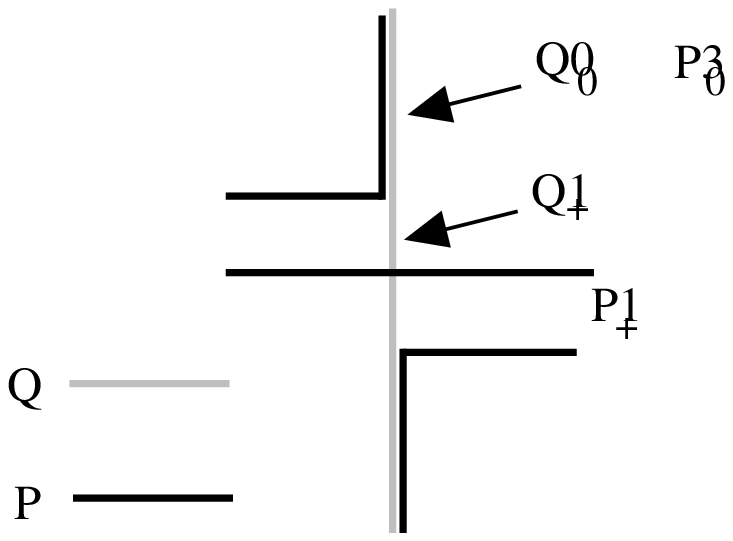}
\relabela <-1pt,0pt> {Q}{$Q$}
\relabela <-1pt,0pt> {Q0}{$Q_0=P_0$}
\relabela <-1pt,0pt> {Q1}{$Q_+$}
\relabela <-1pt,0pt> {P}{$P$}
\relabela <-1pt,0pt> {P1}{$P_+$}
\endrelabelbox}
\caption{}
\end{figure} 

\medskip
We now consider the
possibilities:

\medskip
{\bf Case 1}\qua $P_Y = \emptyset$

Consider the annuli $P_X$ in the context of Lemma \ref{annuli2}.  By the
Rule of Three (Lemma \ref{three}) and the fact that $P_X$ is separating,
$I(\bdd P_X) = (2,0,0,0)$ or $(2,2,0,0)$.  So either $P_X$ is a single
annulus with both ends parallel to the same curve $c_1$ in $Q$, or two
annuli, one with both ends at $c_1$ and the other with both ends at
$c_2$, or two annuli,  each with one end at
$c_1$ and one at
$c_2$.  In the first two cases, since $P_X$ is essential in $X$, $P$ and
$Q$ differ by a cabling into $X$, either on one longitude (example
\ref{symcable}, Variation 1.  See Figure 27) or on two longitudes
(example \ref{doubcable}, Varation 1, with one Dehn surgery done at each
site; see Figure 28.)  


\begin{figure}[ht!] 
\begin{minipage}{.59\textwidth}
\centerline{\relabelbox\small
\epsfxsize=55mm \epsfbox{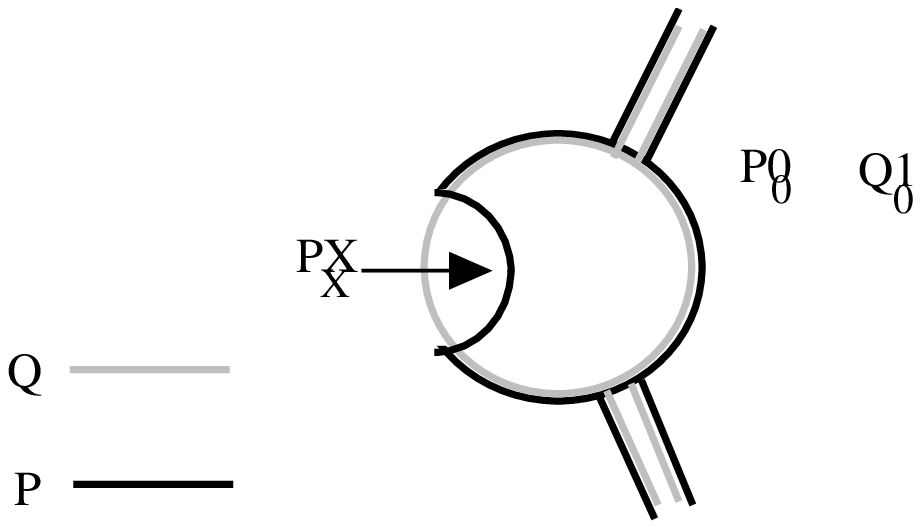}
\relabela <-1pt,0pt> {Q}{$Q$}
\relabela <-1pt,0pt> {P0}{$P_0=\,Q_0$}
\relabela <-1pt,0pt> {P}{$P$}
\relabela <-4pt,-2pt> {PX}{$P_X$}
\endrelabelbox}
\vglue 5mm
\caption{}
\end{minipage} 
\begin{minipage}{.39\textwidth}
\centerline{\relabelbox\small
\epsfxsize=20mm \epsfbox{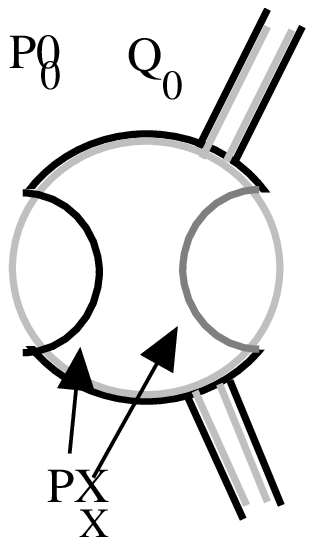}
\relabela <-4pt,0pt> {P0}{$P_0=Q_0$}
\relabela <-4pt,-2pt> {PX}{$P_X$}
\endrelabelbox}
\caption{}
\end{minipage} 
\end{figure}

If $P_X$ is a pair of annuli, each with one end at $c_1$ and one at
$c_2$ then the annuli are non-separating.  The annuli may both be
longitudinal (hence parallel) in $Q$.  Then the splittings
appear as \ref{K4}, Variation 1 (when the $c_i$ are longitudes of
$P$ as well; see Figure 29) or Variation 2, with the Dehn surgery curve
in one of $\mu_{b_{\pm}}$ (when the $c_i$ are twisted in $P$; see
Figure 30).   The annuli $P_X$ could be twisted and not parallel, so
that lying between them is a solid torus on which their cores are torus
knots.  (See Figure 31.)  This is
\ref{K4}, Variation 2, with Dehn surgery in \sss.  Or they could be
twisted and parallel, corresponding to the same Variation but with Dehn
surgery in one of $\mu_{a_{\pm}}$. Note that Variations $3$ and $4$
don't arise, since $M$ is hyperbolike.
 
\begin{figure}[ht!] 
\begin{minipage}{.32\textwidth}
\centerline{\relabelbox\small
\epsfxsize=25mm \epsfbox{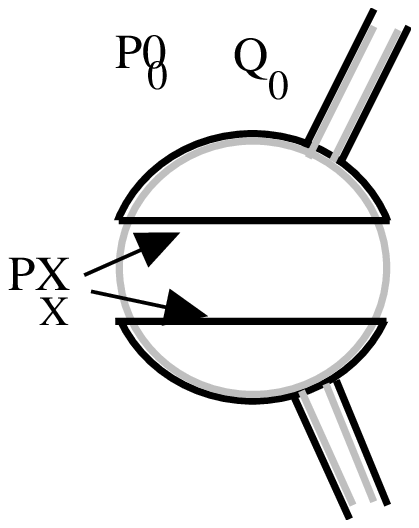}
\relabela <-3pt,0pt> {P0}{$P_0=\,Q_0$}
\relabela <-2pt,-2pt> {PX}{$P_X$}
\endrelabelbox}
\caption{}
\end{minipage} 
\begin{minipage}{.32\textwidth}
\centerline{\relabelbox\small
\epsfxsize=30mm \epsfbox{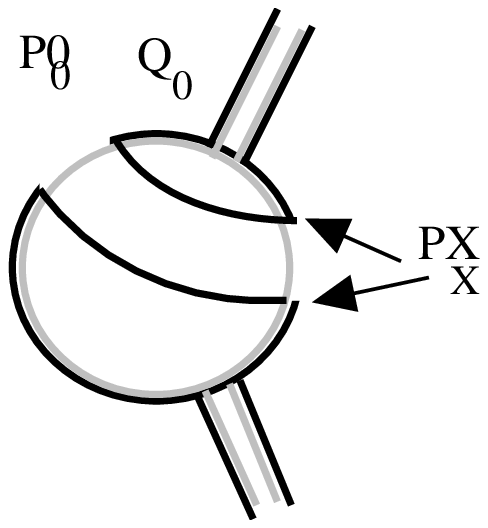}
\relabela <-2pt,0pt> {P0}{$P_0=Q_0$}
\relabela <-2pt,-2pt> {PX}{$P_X$}
\endrelabelbox}
\caption{}
\end{minipage} 
\begin{minipage}{.32\textwidth}
\centerline{\relabelbox\small
\epsfxsize=30mm \epsfbox{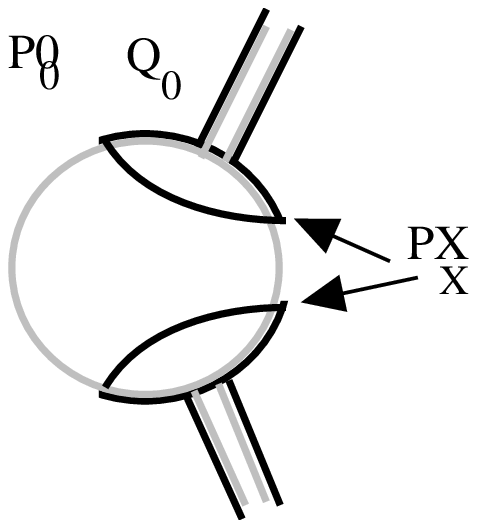}
\relabela <-3pt,0pt> {P0}{$P_0=Q_0$}
\relabela <-2pt,-2pt> {PX}{$P_X$}
\endrelabelbox}
\caption{}
\end{minipage} 
\end{figure}

\medskip
{\bf Case 2}\qua  $P_X$ and $P_Y$ are both non-empty and the end of each
curve in $\bdd P_X \cup \bdd P_Y$ is parallel to one of $c_1$ or $c_2$.

If at least one annulus in each of $P_X$ or $P_Y$ is non-separating,
then together they would give a non-separating, hence essential, torus
in $M$. This contradicts our assumption that $M$ is hyperbolike. So we
may as well assume that each annulus in $P_Y$ is separating.  Hence the
ends of $P_Y$ are twisted in $Y$.  They cannot then be twisted in $X$,
since $M$ is hyperbolike, so $P_X$ is a collection of parallel
non-separating longitudinal annuli in $X$.  

If $P_Y$ has ends at both $c_1$ and $c_2$ (as happens
automatically if $P_X$ has more than two components) then both curves
are twisted in $Y$.  Attach a non-separating annulus in $X$ with ends
at $c_1$ and $c_2$ to the torus (or tori) in $Y$ on which the $c_i$ are
twisted.  The boundary of the thickened result would exhibit a Seifert
manifold in $M$, again contradicting the assumption that $M$ is
hyperbolike.  We conclude that $P_X$ has exactly two
components and that $P_Y$ has ends only at $c_1$, say.  

If there were
more than two annuli in $P_Y$ (hence more than four ends of $\bdd P_Y$)
there would have to be more than two ends of $P_X$ at $c_1$, so we
conclude that $P_Y$ is made up of one or two annuli.  If it's two
annuli, necessarily separating and parallel in $Y$, then the relation
between $P$ and $Q$ can be seen as follows (See Figure 32): In \ref{K4},
Variation 2 let $P$ be the  splitting given there with Dehn surgery
curve in $\mu_{a_+}$ and $Q$ be the same splitting given there but with
Dehn surgery curve in $\mu_{a_-}$. To view these simultaneously as
splittings of the same manifold $M$, of course, the Dehn surgery curve
has to be moved from $\mu_{a_+}$ to $\mu_{a_-}$, dragging some annuli
along, until the splitting surfaces $P$ and $Q$ intersect as described.


\begin{figure}[ht!] 
\centerline{\relabelbox\small
\epsfxsize=45mm \epsfbox{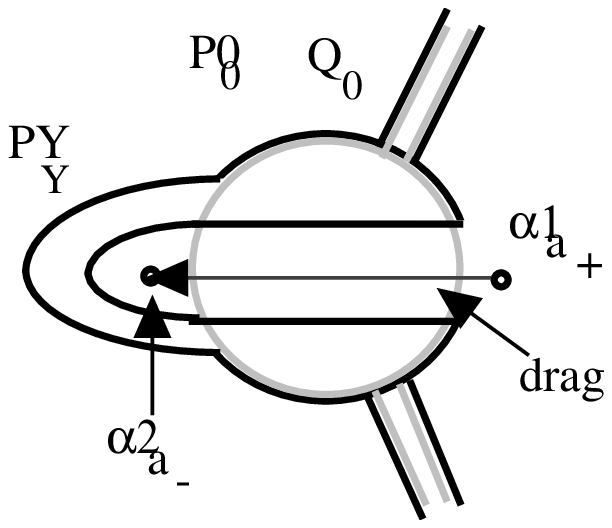}
\relabela <-2pt,0pt> {P0}{$P_0=\,Q_0$}
\relabela <-2pt,-3pt> {PY}{$P_Y$}
\relabela <-2pt,-2pt> {drag}{drag}
\relabela <-2pt,-2pt> {a1}{$\alpha_{a_+}$}
\relabela <-2pt,2pt> {a2}{$\alpha_{a_-}$}
\endrelabelbox}
\caption{}
\end{figure}

If $P_Y$ is a single annulus, it must have one
end on $P_0$ and one end on an end of $P_X$, and the annulus is twisted
in $Y$.  The initial splitting by $Q$ is as in Example \ref{K4}
Variation 1 ($X = A_- \cup \sss$), with a Dehn surgery curve lying in
$\mu_{b_+}$, say.  If the splitting is altered by first putting the
Dehn surgery curve in $\mu_{a_+}$  (yielding the
same manifold $M$), then altering as in Example
\ref{K4} (ie, considering \aub\ where $B = B_- \cup \sss$) and
then dragging the Dehn surgery curve from  $\mu_{a_+}$ to
$\mu_{b_+}$, pushing before it an annulus from the $4$--punctured sphere
along which $A_-$ and $B_-$ are identified, we get the splitting
surface $P$, intersecting $Q$ as required. (See Figure 33.)


\begin{figure}[ht!] 
\centerline{\relabelbox\small
\epsfxsize=45mm \epsfbox{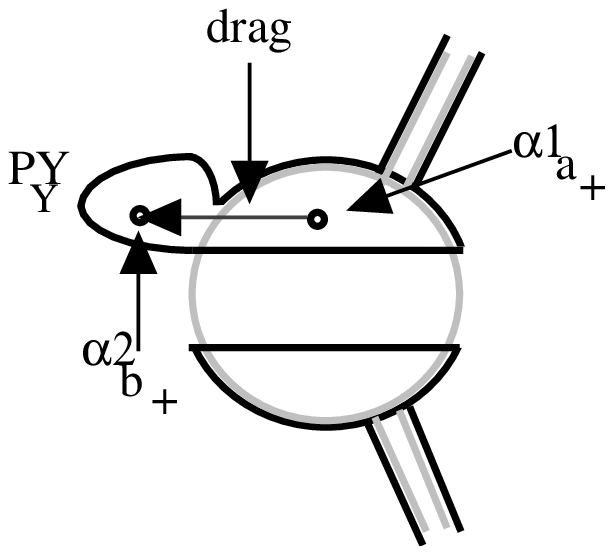}
\relabela <-2pt,-3pt> {PY}{$P_Y$}
\relabela <-2pt,-2pt> {drag}{drag}
\relabela <-2pt,-2pt> {a1}{$\alpha_{a_+}$}
\relabela <-2pt,2pt> {a2}{$\alpha_{a_-}$}
\endrelabelbox}
\caption{}
\end{figure}

\medskip
{\bf Case 3}\qua  Some end component(s) of $P_X$ or $P_Y$ lie parallel to
each of $c_1, c_2, c_3$,

Then one of the $c_i$, say $c_1$, is parallel only to ends of $P_X$,
and at most (hence exactly) two of them.  Another, say $c_3$, is
parallel only to ends of $P_Y$ (again exactly two). 

\medskip
{\bf Subcase a}\qua No ends of $P_Y$, say, are parallel to $c_2$.

Then all ends of $P_Y$ are parallel to $c_3$ and some ends of $P_X$ are
parallel to $c_1$ and some to $c_2$.  So, by the Rule of Three (Lemma
\ref{three}), $P_Y$ is a single separating annulus with ends at $c_3$
and $P_X$ is either a pair of separating annuli, one
each with ends at
$c_1$ and $c_2$, or a pair of non-separating annuli, each having one
end at $c_1$ and one end at $c_2$. (See Figures 34 and 35.) It follows
that when
$X$ is cut open along
$P_X$ the component that contains $c_3$ is a genus two handlebody in
which the cores of the annuli $P_X$ appear as separated curves (at
least one a longitude), and
$c_3$ is a longitude not
parallel to either (since $P$ splits $M$ into handlebodies).  Then the
technical Lemma
\ref{tech0} ($c_3$ here becomes $c_2$ there) precisely places $c_3$ with
respect to the annuli.  In particular, if the $P_X$ are separating, (so
the argument is precisely symmetric moving from $P$ back to $Q$), the
splitting is given in Example
\ref{symcable}, Variant 2.  If the $P_X$ are non-separating, the
splitting is given in Example \ref{K4}, Variation 5 (if $c_1$ and
$c_2$ are not twisted in any of the handlebodies) or Variation 6, with
one Dehn surgery curve appropriately placed (if the curves $c_1$ or
$c_2$ are twisted in a handlebody).   


\begin{figure}[ht!] 
\centerline{\relabelbox\small
\epsfxsize=72mm \epsfbox{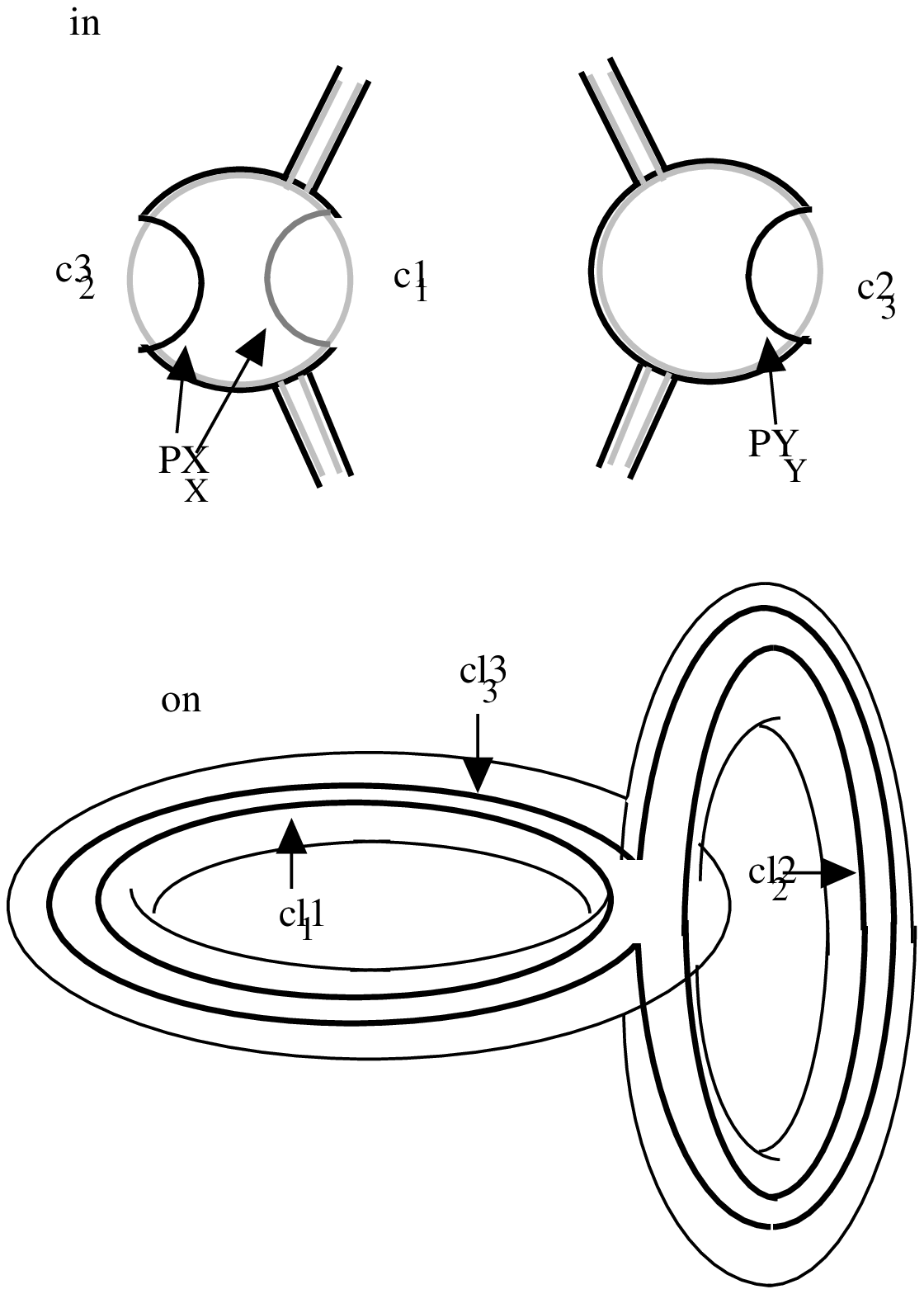}
\relabela <-2pt,-2pt> {PY}{$P_Y$}
\relabela <-2pt,-2pt> {PX}{$P_X$}
\relabela <-2pt,0pt> {c1}{$c_1$} 
\relabela <-2pt,0pt> {c2}{$c_2$} 
\relabela <-2pt,0pt> {c3}{$c_3$} 
\relabela <-2pt,0pt> {cl1}{$c_1$} 
\relabela <-2pt,0pt> {cl2}{$c_2$} 
\relabela <-2pt,-2pt> {cl3}{$c_3$} 
\relabela <-10pt,0pt>  {on}{on $X$ (actual)}
\relabel {in}{in $X$ (schematic)\qquad in $Y$ (schematic)}
\endrelabelbox}
\caption{}
\end{figure}


\begin{figure}[ht!] 
\centerline{\relabelbox\small
\epsfxsize=.9\hsize\epsfbox{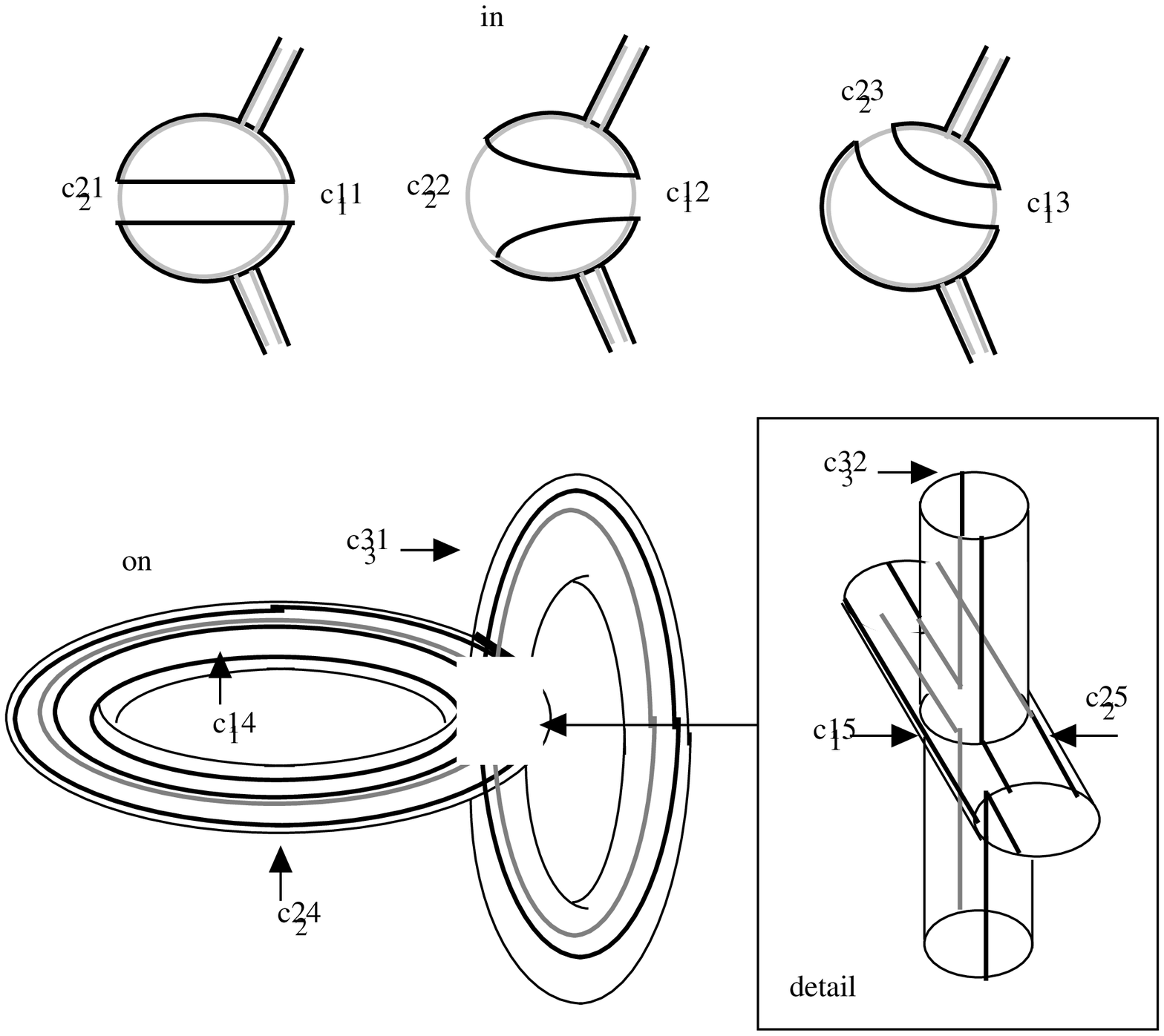}
\relabela <-2pt,0pt> {c11}{$c_1$} 
\relabela <-2pt,0pt> {c21}{$c_2$} 
\relabela <0pt,0pt> {c31}{$c_3$} 
\relabela <-2pt,0pt> {c12}{$c_1$} 
\relabela <-2pt,0pt> {c22}{$c_2$} 
\relabela <0pt,0pt> {c32}{$c_3$} 
\relabela <-2pt,0pt> {c13}{$c_1$} 
\relabela <-2pt,0pt> {c23}{$c_2$} 
\relabela <-2pt,0pt> {c14}{$c_1$} 
\relabela <-2pt,0pt> {c24}{$c_2$} 
\relabela <-2pt,0pt> {c15}{$c_1$} 
\relabela <-2pt,0pt> {c25}{$c_2$} 
\relabela <-10pt,0pt> {on}{on $X$ (actual)}
\relabela <-10pt,0pt> {in}{in $X$ (schematics)}
\relabela <-2pt,0pt> {detail}{detail} 
\endrelabelbox}
\caption{}
\end{figure}

\medskip
{\bf Subcase b}\qua Ends of both $P_X$ and $P_Y$ are parallel to
$c_2$.    

The curve $c_2$ can be twisted in at most one of $X$ and $Y$, so assume
that $c_2$ is not twisted in $X$.  By \ref{annuli2} this means that
$P_X$ is a pair of parallel annuli running between longitudes $c_1$ and
$c_2$ of $X$.  

$P_Y$ has two ends at $c_3$ and, as in Case 2, either $2$ or $4$ ends
at $c_2$. If the annuli are separating and there are $4$
ends of $P_Y$ at $c_2$ then the cores of the two annuli whose ends
are at $c_2$ cobound an annulus in $Y$.  It follows from \ref{tech1}
that any twisted or longitudinal curve in $P$ must be parallel to one
of these cores.  But that would make the core of the annulus in $P_Y$
at $c_3$ parallel to one of these cores, hence $c_3$ parallel to
$c_1$ or $c_2$ in $P$.  Since this is impossible, the case does not
arise.

Suppose the annuli in $P_Y$ are separating and there
are two ends of $P_Y$ at $c_2$. Then $P_Y$ consists of two
separating annuli, $\Aaa_2$ which has both ends at $c_2$ and $\Aaa_3$
which has both ends at $c_3$. The situation is analogous to the last
example in Case 2 above, $K4$ variation 2, with Dehn surgery curve
in  $\mu_{b_+}$ when \sss\ is attached to $A_-$, and in $\mu_{a_+}$
when \sss\ is attached to $B_-$. But the presence of $\Aaa_3$ adds the
additional complexity of Variation 6:  \rrr\ is simultaneously moved
from $B_-$ to $A_-$.

If the annuli in $P_Y$ are not all separating, and $P_Y$ has four ends
at $c_2$ then $P_Y$ consists of a pair of non-separating annuli with
ends at $c_2$ and $c_3$ and a single separating annulus with ends at
$c_2$.  Much as in the other case when there were $4$ ends at $c_2$
this leads to a contradiction, this time with Lemma \ref{tech2}. (See
Figure 36.)


\begin{figure}[ht!] 
\centerline{\relabelbox\small
\epsfxsize=36mm \epsfbox{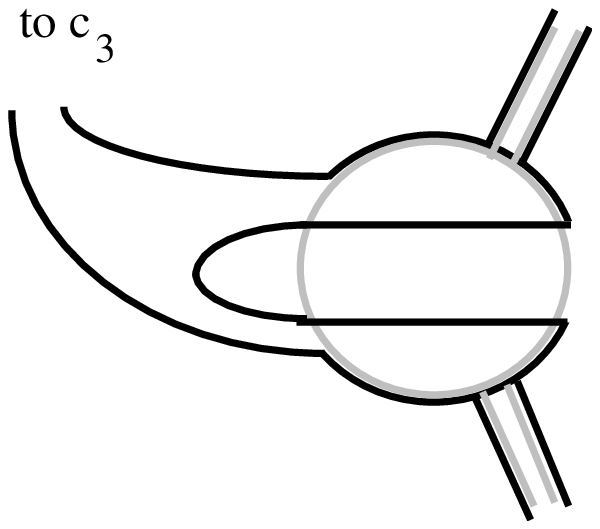}
\relabela <-5pt,0pt> {to c}{to $c_3$} 
\endrelabelbox}
\caption{}
\end{figure} 

The final possibility is that the annuli in $P_Y$ are not all
separating, and
$P_Y$ has two ends at $c_2$.  Then $P_Y$ consists of a pair of
non-separating annuli with ends at $c_2$ and $c_3$.  (See Figure 37).


\begin{figure}[ht!] 
\centerline{\relabelbox\small
\epsfxsize=.8\hsize \epsfbox{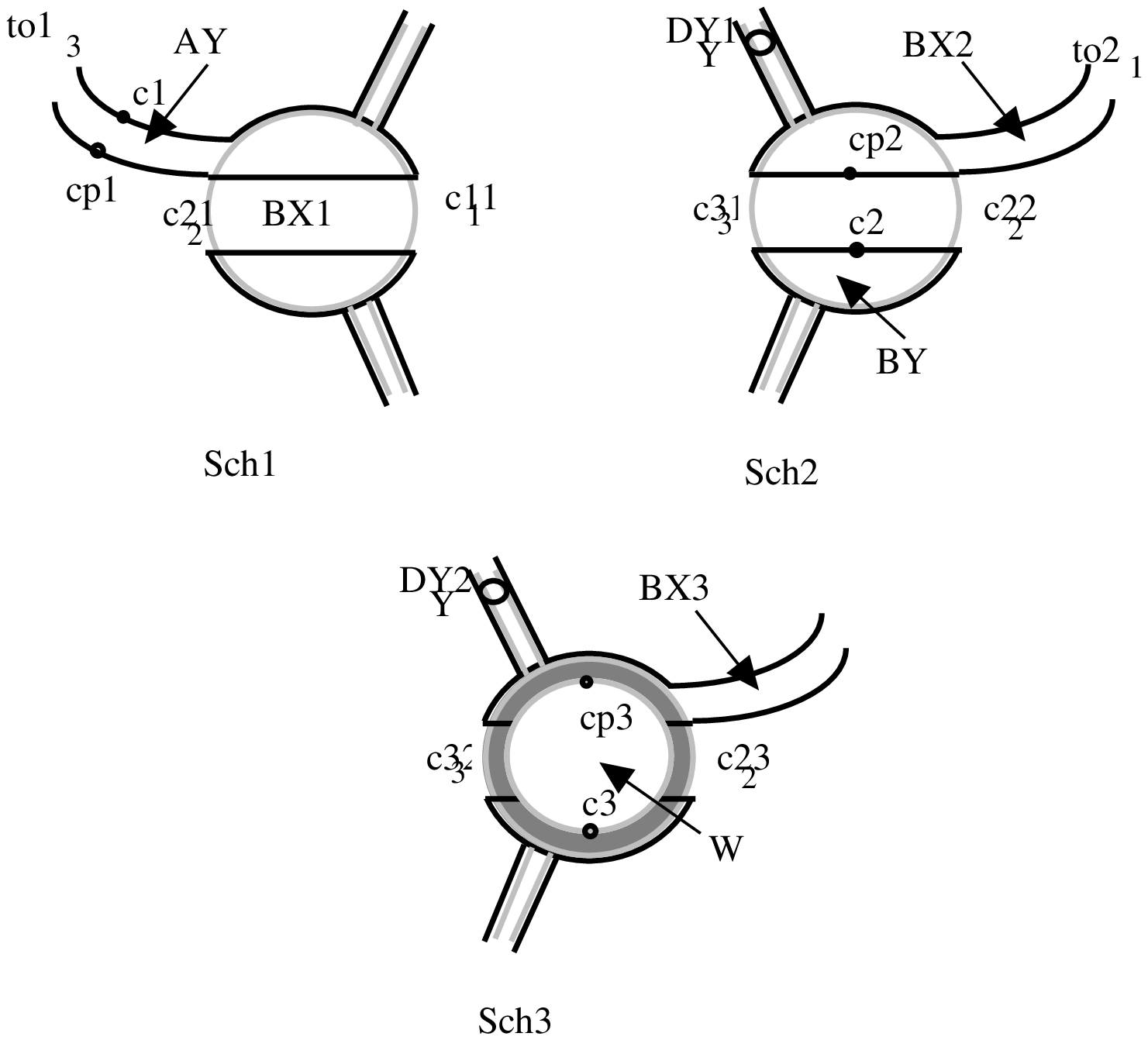}
\relabela <-2pt,0pt> {to1}{to $c_3$} 
\relabela <-2pt,2pt> {to2}{to $c_1$} 
\relabel {BX1}{$B\cap X$}
\relabel {BX2}{$B\cap X$}
\relabel {BX3}{$B\cap X$}
\relabel {BY}{$B\cap Y$}
\relabel {AY}{$A\cap Y$}
\relabel {Sch1}{Schematic in $X$}
\relabel {Sch2}{Schematic in $Y$}
\relabel {Sch3}{Schematic of $Y'$}
\relabela <-2pt,0pt> {DY1}{$D_Y$}
\relabela <-2pt,0pt> {DY2}{$D_Y$}
\relabela <-2pt,0pt> {cp1}{$c'$}
\relabela <-2pt,0pt> {cp2}{$c'$}
\relabela <-2pt,0pt> {cp3}{$c'$}
\relabela <-2pt,0pt> {c1}{$c$}
\relabela <-2pt,0pt> {c2}{$c$}
\relabela <-2pt,0pt> {c3}{$c$}
\relabela <-2pt,0pt> {c11}{$c_1$}
\relabela <-2pt,0pt> {c21}{$c_2$}
\relabela <-2pt,0pt> {c22}{$c_2$}
\relabela <-2pt,0pt> {c23}{$c_2$}
\relabela <-2pt,0pt> {c31}{$c_3$}
\relabela <-2pt,0pt> {c32}{$c_3$}
\relabela <-2pt,0pt> {W}{$W$}
\endrelabelbox}
\caption{}
\end{figure}

We will show that this case, too, cannot arise, for if it did:

\medskip
{\bf Claim}\qua {\sl Then $M$ is a Seifert manifold.}

\medskip
{\bf Proof of claim}\qua  Suppose, with no loss of generality,  that the
(solid torus) region between the annuli $P_Y$ in $Y$ lies in $A$.   Then
$B \cap Y$ is a genus two handlebody $H$ in which the two annuli $P_Y$,
have cores $c, c'$, which are separated curves on the boundary.  $B$ is
obtained from  $H = B \cap Y$ by attaching the region $(B \cap X) \cong
(annulus \times I)$ that lies between the annuli $P_X$. One end of
$annulus \times I$ is attached to a curve which is parallel in $\bdd Y$
to $c_2$ and in $\bdd H$ to $c'$ say.  The other end is attached at
$c_1$.  It follows from \ref{tech} that $c$ is longitudinal in $H$, and
$c_1$ crosses exactly twice a disk $D_Y \subset H$ that separates $c$
and
$c'$. Notice that $Y$ is obtained from $H$ by attaching a thickened
annulus with ends at $c$ and $c'$.  So $D_Y$ is a non-separating
meridian disk for $Y$. 

Cut $Y$ open along $D_Y$ to get a solid torus containing $P_Y$ and
remove from this solid torus all but a collar of the boundary.  That
is, remove a solid torus $W \subset Y$ whose complement $Y_-$ in $Y$
consists of a collar of $\bdd W$ to which a $1$--handle, with cocore
$D_Y$, is attached.  Note that the two annuli
$P_Y$ intersect $Y_-$ in four parallel spanning annuli.  Let \aaa\
denote the slope of their ends on $\bdd W$.  

We will show that $M - W$
is a Seifert manifold. In fact, we will
show that $M - W$ fibers over the circle with fiber the
three--punctured sphere (ie, the pair of pants) and this
suffices. Then a Seifert
manifold structure on $M$ will be obtained by filling in the solid
torus $W$.

Another way to view the compression body $Y_-$ is to begin with $torus
\times I$ and attach to it a genus two handlebody $H' \subset H$ by
attaching collars of two separated longitudes $c, c'' \subset \bdd
H'$ to two annuli in $torus
\times \{ 1 \}$.  The annuli have slope \aaa. Note that whereas $c'$
may be twisted in $H$, we've obtained $H'$ by removing a solid torus
so large it contains any cable space, so the attachment circle $c''$ is
indeed longitudinal in $H'$.   Now when $(B \cap X)
\cong (annulus \times I)$ is attached to $Y_-$, one end is attached to 
$torus \times \{ 1 \}$ along a curve parallel to \aaa\ and the
other end is attached along the curve $c_1 \subset \bdd H$ which crosses
$D_Y$ twice.  Then by Lemma \ref{tech}, $H' \cup (B \cap X)$ is a genus
two handlebody homeomorphic to $V \times I$, where $V$ is a pair of
pants whose boundary is the triple of curves $c, c_1, c''$. The upshot
is that $Y_- \cup (B \cap X)$ can be viewed as $torus \times I$ with $V
\times I$ attached by identifying each of $c \times I, c_1
\times I, c'' \times I$ with different parallel annuli on $torus \times
\{ 1 \}$.  Each annulus has slope \aaa. 

What remains of $M - W$ is the genus two handlebody $X \cap
A$.  When this is attached along its entire boundary, it's easy to see
that the three annuli remaining on $torus \times \{ 1
\}$ are identified with annuli corresponding to two distinct
separated longitudes $d_1, d_2$ in $\bdd (X \cap A)$ (both parallel to
$c_2$ in 
$\bdd X$) and an annulus whose core is $c_3$.  The fact that $A$ is
constructed by attaching $(A \cap Y) \cong (S^1 \times D^2)$ to the
handlebody $A \cap X$ along $c_2$ and $c_3$ means (see
Lemma \ref{tech+}) that
$A \cap X$ can be viewed as $V' \times I$, where $V'$ is a pair of
pants with $\bdd V'$ the three curves $d_1, d_2, c_3$.
The upshot is that adding in $A \cap X$ is the same as attaching $V'
\times I$ by identifying $\bdd V' \times I$ with the remaining annuli of
$torus \times \{ 1 \}$ and then the rest of the boundary, $V' \times
\bdd I$, with $V \times \bdd I$.  But the union of $V \times I$ and
$V' \times I$ fibers over the circle with fiber a pair of pants. It's
easy to show then that this manifold is also Seifert fibered by circles
transverse to the pair of pants, with a generic fiber running through
each of $V \times I$ and $V' \times I$ exactly three times and, in
$torus \times \{ 1 \}$ crossing a curve with slope \aaa\ exactly
once.  The base is a disk and there are two exceptional fibers, each of
type $(3, 1)$.  

Now $M$ is created from the Seifert manifold just described by
filling in a solid torus along its boundary, namely the solid torus
that was removed from $Y$ at the beginning.  The result is either a
Seifert manifold (if the filling slope differs from that of the fiber)
or a reducible manifold, if the filling slope is that of the fiber. 
In any case, $M$ is not hyperbolic.  \qed

\section{Heegaard splittings when tori are present}

Suppose $M$ is a closed orientable irreducible $3$--manifold of Heegaard
genus two, $M$ is not itself a Seifert manifold, and $M$ contains an
essential torus.  Following Section \ref{character}, let \Fff\ denote
the collection of tori that constitute the canonical tori of
$M$. The discussion of Section \ref{character} shows that a genus two
Heegaard splitting $M = \aub$ can be isotoped to intersect \Fff\ so
that there is exactly one component (here to be denoted
$W_P$) of $M - \Fff$ which contains a non-annular part of $P$. 
Moreover, $V_P = M - W_P$ is a Seifert manifold with at most two
components, each of which $P$ intersects in fibered annuli. (More is
shown there about $V_P$).  $W_P$ is atoroidal, but could perhaps be a
Seifert manifold over a disk with two exceptional fibers or over an
annulus or M{\"o}bius band with one exceptional fiber, as long as the
fibering doesn't match a fibering of $V$. As in Section
\ref{character} we let $P_{-} = P \cap W_P$ and let $A_- = A \cap W_P,
B_- = B \cap W_P$ have spines  $\Ss_A$ and $\Ss_B$.

Consider how two different such splittings $M = \aub = \xuy$ compare.
One possibility is 

\begin{example}
\label{split}

$W_P = V_Q$ and $W_Q = V_P$.  

\end{example}  

Then in each of $W_P$ and $W_Q$ there are essential annuli whose slopes
differ from those of the Seifert fibering $\bdd V_P$ and $\bdd V_Q$
respectively. Exploiting Lemma \ref{squar} and Proposition
\ref{splitchar} we can write down explicit and simple descriptions of
all variations possible here and deduce that, for these two
splittings, the commutator $\Thp \Thq \Thp^{-1}
\Thq^{-1}$ can be obtained by Dehn twisting around the unique essential
torus \Fff.

So henceforth we will assume that $W_P = W_Q$ and $V_P = V_Q$ and
revert to $W$ and $V$ as notation.  Then
$W = \awb = \xwy$, where the splitting surface $Q_-$, the
handlebodies $X_-, Y_-$ and the spines $\Ss_X$ and
$\Ss_Y$ are defined analogously to
$P_{-}, A_-, B_-, \Ss_A$ and $\Ss_B$.

\begin{thm}
\label{tori}

Suppose $M = \aub = \xuy$ are two non-isotopic genus two Heegaard
splittings of an irreducible orientable closed $3$--manifold $M$. 
Suppose $M$ contains an essential torus.  Then either the
splittings are isotopic, or the relation between the Heegaard
splittings is described in one of the variations of one of the
examples in Section
\ref{examples} or in Example \ref{split}.  In particular, the commutator
$\Thp \Thq \Thp^{-1} \Thq^{-1}$ can be obtained by Dehn twists around
essential tori in $M$, and the two splittings become equivalent after
a single stabilization.

\end{thm}

\pf 
Isotope $P$ and $Q$ so that they each intersect the canonical tori
\Fff\ of $M$ as described in (possibly different cases of) Section
\ref{character}, and continue with the same notation.  There are two
possibilities:

\noindent

\medskip
{\bf Case 1}\qua $\bdd P_-$ and $\bdd Q_-$ have the same slope on each
component of \Fff.

Then the annuli in $\Ss_A$ and $\Ss_B$ can be chosen to
overlap so that their complements in \Fff\ are disjoint.  This means
that during the sweep-outs of $W - \eta (\Ss_A \cup \Ss_B)$ and $W - 
\eta (\Ss_X \cup
\Ss_Y)$ determined by $P_{-}$ and $Q_-$ respectively,  $\bdd P_{-} \cap
\bdd Q_- = \emptyset$.  (See the discussion preceding \ref{squar} for a
description of the sweepout).  In particular, the generic
intersection of $P_{-}$ and
$Q_-$ during the sweepout consists of closed curves. Apply the argument
of Sections \ref{positioning} and
\ref{aligning} almost verbatim to the two sweep-outs. 
The upshot is a positioning of
$P_{-}$ and $Q_-$ so that they are aligned except along some collection
of subannuli.  That is, $(P_-)_X = closure(P_{-} - Y_-)$ and $(P_-)_Y =
closure(P_{-} - X_-)$ consist of incompressible annuli in $X_-$ and
$Y_-$ respectively and none of these is parallel in $X_-$ or $Y_-$ to a
subannulus of $Q_-$.  

Consider first of all the case in which 
$P_{-}$ and $Q_-$ are both $4$--punctured spheres, so any incompressible
annulus with boundary disjoint from \Fff\ is $\bdd$--parallel.  (This
excludes only the case when $V$ fibers over the circle with two
exceptional fibers and either $P$ or $Q$ intersects $V$ as in the
single annulus case.)  Then
$(P_-)_X$ and $(P_-)_Y$ consist entirely of annuli parallel to one of
the two annuli in $\Fff \cap X_-$ (resp.\ $\Fff \cap Y_-$).  It's easy
to see that these can be removed by an isotopy of $P_{-}$ which slides
$\bdd P_{-}$ around \Fff.  Thus, after an isotopy of $P_{-}$ which may
move the boundary of $P_{-}$, we can make $P_{-}$ and $Q_-$ coincide. 
Such an isotopy is equivalent to Dehn twists around tori in \Fff.  The
fibered annuli $P \cap V$ and $Q \cap V$ may also differ within $V$, but
can be made to coincide by Dehn twists around essential tori in $V$.

Suppose now that $P_{-}$ and $Q_-$ are both twice--punctured tori. 
This can arise when $V$ fibers over the disk with two exceptional
fibers, and both $P$ and $Q$ intersect as in the single annulus case. 
Then more complicated
essential annuli $(P_-)_X$ and
$(P_-)_Y$ can occur.  In any of $A_-, B_-, X_-, Y_-$, say $X_-$,
essential annuli with boundaries disjoint from
$\Fff$ can be of two types: parallel annuli non-separating in $Q_-$,
each with one end parallel to $\bdd Q_-$ and the other parallel to a
curve $c' \subset Q_-$; or parallel separating annuli, both ends
parallel to the same single twisted curve $c' \subset Q_-$. (See
Figure 38.)


\begin{figure}[ht!] 
\centerline{\relabelbox\small
\epsfxsize=85mm \epsfbox{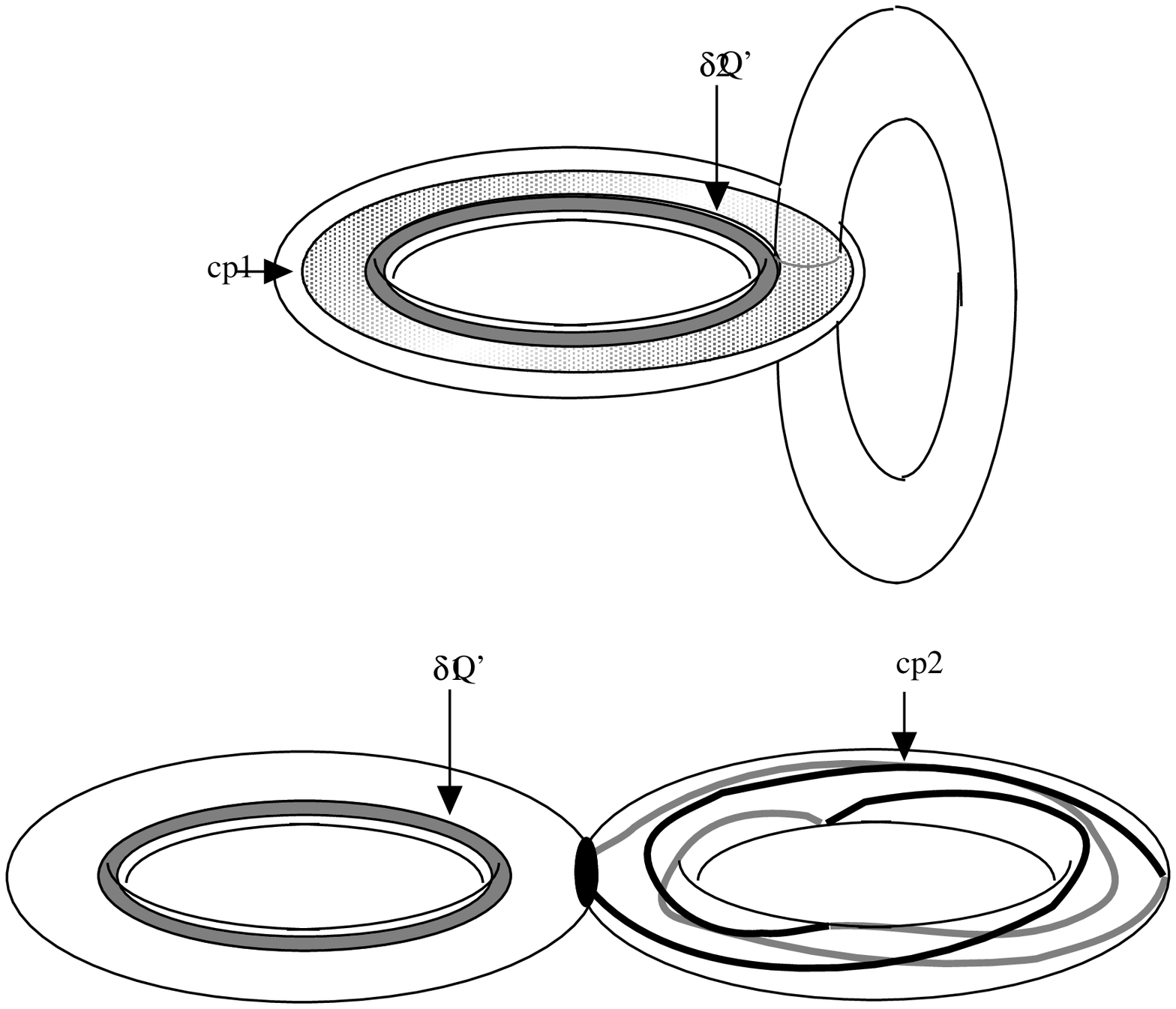}
\relabela <-2pt,0pt> {cp1}{$c'$}
\relabela <-2pt,0pt> {cp2}{$c'$}
\relabela <-2pt,0pt> {d1}{$\partial Q'$}
\relabela <-2pt,0pt> {d2}{$\partial Q'$}
\endrelabelbox}
\caption{}
\end{figure}

Since \Fff\ is the set of canonical tori of $M$, there is no
essential torus in $W$, nor is there an essential annulus
with end having the slope of $\bdd Q_-$ ($= \bdd P_-$).  It follows that
one of $(P_-)_X$ or $(P_-)_Y$ is empty.  For if both $(P_-)_X$ and
$(P_-)_Y$  were non-empty and separating, then copies of each could be
matched up along their boundaries in $Q_-$ to create an essential
separating torus in $W$.  If both were non-separating, then their ends
could be matched up to create an essential non-separating torus in
$W$.   If $P_Y$ were separating and $P_X$ were non-separating (or vice
versa), then an annulus of $P_Y$ attached to two annuli in $P_X$, each
with their other end on \Fff, would give an essential annulus in $W$
with slope the fiber of $V$.

So we may assume $(P_-)_Y$ is empty, and so $(P_-)_X$ consists of one
separating or two parallel and non-separating annuli.  If
$(P_-)_X$ is a single separating annulus then the splittings are
completely described by Example \ref{doubcable} Variation 1, with three
Dehn surgery curves.  One pair produces $V$ the other cables $A$ into
$B$ to produce $\xuy$ and vice versa.  

Suppose $(P_-)_X$ consists of two parallel non-separating annuli and 
suppose that the region between the parallel annuli of $(P_-)_X$ lies in
$A_-$, say.  Then $A_- \cap X_-$ can be viewed as $\sss = square \times
S^1$, where $(P_-)_X$ comprises two opposite annuli in the boundary
and the other pair of opposite annuli is $(Q_-)_A$.  Since the annuli
$(P_-)_X$ are $\bdd$--compressible in $X$ it follows that, viewed in
$P$, one of the annuli is longitudinal in $B$.  Similarly, one of the
annuli in $(Q_-)_A$ is longitudinal in $Y$.  This suffices to
characterise the two splittings as those of Example $\ref{K4}$ Variation
3, with one of the Dehn
surgery curves placed in one of $\mu_{a_{\pm}}$ and the other in one of
$\mu_{b_{\pm}}$.

Finally, suppose that  say, $P_{-}$ is a $4$--punctured
sphere, and $Q_-$ is a twice--punctured torus.  This means that  $V$
fibers over the disk with two exceptional fibers; that $P$ 
intersects $V$ as in the parallel annuli case; and that $Q$ intersects
$V$ as in the single annulus case.  Then the
argument is a mix of earlier ideas:  Once again, after an isotopy of
$\bdd Q_-$ on $\Fff$ we can ensure that the annuli in $Q_-$ which are
not aligned with $P_-$ are not contained in collars of $\bdd Q_-$ and,
in $A_-$ and $B_-$, these annuli are parallel to the annuli $\Fff \cap
A$ and/or $\Fff \cap B$ respectively. It follows that $P_{- 0}$ (that
is, the part of $P_-$ that is aligned with $Q_-$) is a single
$4$--punctured sphere.  It follows then that the complement of $Q_{- 0}$
is a single collection of parallel annuli. Since $P_-$ has two more
boundary components than $Q_-$ the annuli of $P_-$ not aligned with $Q$
include collars of exactly two components of $\bdd P_-$.  Since no
component of $P_{- 0}$ is an annulus, it follows that in fact $Q_-$ is
aligned with $P_-$ except for a single annulus, lying in $A_-$ say. 
That annulus cuts off collars of two components of $P_-$, which are the
only parts of $P_-$ not aligned with $Q_-$.  Put another way:  $Q_-$ is
obtained from $P_-$ by attaching a copy of an annulus component of
$\Fff \cap A$.  Then the setting is exactly as in Lemma \ref{switch} and
preceding.  In particular, both splittings are described in \ref{K4}
Variation 3.

\noindent \medskip
{\bf Case 2}\qua $\bdd P_-$ and $\bdd Q_-$ have different slopes
on some component of \Fff.

Then $V$ is the neighborhood of either a one-sided Klein bottle or a
non-separating torus, for otherwise the slope of $P_-$ and $Q_-$ must
be that of the unique Seifert fibering of $V$. We will concentrate on
the latter, for the proof in the former, more specialized, case is
similar but easier:  A combinatorial proof comparing $P_-$ and $Q_-$,
much as in the non-separating torus case below, shows that there is an
essential annulus in $W$, so $W$ is in fact a Seifert piece attached to
$V$, but the fibers do not match.  This is Example \ref{split}.

The argument when $V$ is the neighborhood of a non-separating torus
will eventually bear a striking resemblance to the hyperbolic case,
Section \ref{positioning}.  

$W$ is the manifold obtained by cutting open along
the non-separating torus \Fff. $\bdd W$ consists of two copies
of \Fff, which we denote $\bdd^{\pm} W$. We will denote $\bdd P_- \cap
\bdd^{\pm} W$ by $\bdd^{\pm} P$ (and similarly for $\bdd^{\pm} Q$).

\noindent \medskip
{\bf Subcase 2a}\qua  $W \cong T^2 \times I$.

Then $M$ is the mapping
cylinder of a torus.  It is shown in \cite{TO} (more detailed argument
relevant here can be found also in \cite{CS}) that the only such
mapping cylinders allowing a genus two splittings are those with
monodromy of the form
$$L = \left(\begin{array}{cc} \pm m & -1 \\ 1 & 0 \end{array}
\right).$$ If, for example, $\aub$ is the genus two
splitting, then with respect to the coordinates for which $L$ is
the monodromy, the slope of
$P_-$ is $\pm \left( \begin{array}{c} 0 \\ 1 \end{array} \right)$. 
Hence, when we consider both splittings $\aub$ and $\xuy$, it
follows that there is an automorphism of the torus
$\Fff$ that carries the slope of $P_-$ to that of $Q_-$ and this
automorphism must commute with $L$.   If $|m| \leq 2$ it is easy to
check that the matrix of such an automorphism is a power of $L$. It
follows that the original splitting
$P$ can be ``spun'' around the mapping cylinder until the slopes of
$P_-$ and
$Q_-$ coincide, and so, as above, $P$ and $Q$ are isotopic.  (See \cite{CS} for more detailed
explanation.) 

If $|m| \geq 3$ then $M$ is a solvmanifold, whose Heegaard splittings
are described in \cite{CS}:  With precisely two
exceptions, each solvmanifold has exactly one isotopy class of
irreducible Heegaard splittings (sometimes genus two, sometimes genus
three).  The two exceptions, corresponding to the case $|m| = 3$, each
have exactly two genus two splittings, for which the associated
standard involutions commute, as desired.

\noindent \medskip
{\bf Subcase 2b}\qua $W$ contains an essential spanning annulus

That is, $W$ contains an essential annulus \Aaa\ with one end on each
of  $\bdd^{\pm} W$.  These ends are denoted $\bdd_{\pm} \Aaa$.  Note
that if $W$ contains two essential spanning annuli of different slopes
then, since $M$ is irreducible,  $W \cong T^2 \times I$ and we are
done by the previous subcase.  So we may as well assume that $W$
contains a unique (up to proper isotopy) essential spanning annulus
\Aaa. 

\begin{lemma}
\label{spanner}

Neither end of \Aaa\ is parallel to $\bdd P_-$ (or $\bdd Q_-$).

\end{lemma}

\pf  If both ends are parallel to $\bdd P_-$, then, arguing as in
\ref{squar}, we can arrange that $P_- \cap \Aaa$ consists of
essential closed curves, parallel to the core of $\Aaa$.  Then
isotope $\Aaa$ to minimize the number of curves; the result is that
$\Aaa$ can be made disjoint from $P_-$ and so lies entirely in $A_-$ or
$B_-$.  But this would imply that $A$ or $B$ contained an essential
torus, obtained by isotoping the two curves $\bdd_{\pm} \Aaa$ so that
they coincide in
$\Fff$.  But a handlebody does not contain an essential torus.  

If one end of \Aaa\ is parallel to $\bdd P_-$ and the other
end is not, then the involution
$\Thp | W$, which interchanges $\bdd^+ W$ and $\bdd^- W$, carries \Aaa\
to a second spanning annulus in $W$ whose slope at each end differs
from that of \Aaa, contradicting our hypothesis that $W$
contains a unique essential spanning annulus \qed

\begin{lemma}
\label{disjarcs}

$P_-$ does not contain two disjoint arcs \aaa\ and \bbb, the former
boundary compressing via a disk in $A_-$ and the latter via a disk in
$B_-$.  (Similarly for $Q_-$.)

\end{lemma}

\pf Suppose such curves existed.  The ends of \aaa\ lie on $\bdd^+
W$. A $\bdd$--compression of $P_-$ along \aaa\ changes it to a pair of
pants with one boundary component an inessential circle in $\bdd^+ W$. 
It follows that any essential arc in $P_-$ that is disjoint from
\aaa\ and which has both ends on $\bdd^- W$ will $\bdd$--compress via a
disk in $A_-$.  In particular, \bbb\ also has both ends on $\bdd^+ W$. 
(See Figure 39.)  Then simultaneous $\bdd$--compressions on both \aaa\
and \bbb\ give two parallel spanning annuli in $W$.  Their ends on
$\bdd^+ W$ have slope perpendicular to that of $\bdd^+ P_-$ and on
$\bdd^- W$ they have slope parallel to $\bdd^- P_-$. The result then
follows from Lemma \ref{spanner}.  \qed


\begin{figure}[ht!] 
\centerline{\relabelbox\small
\epsfxsize=85mm \epsfbox{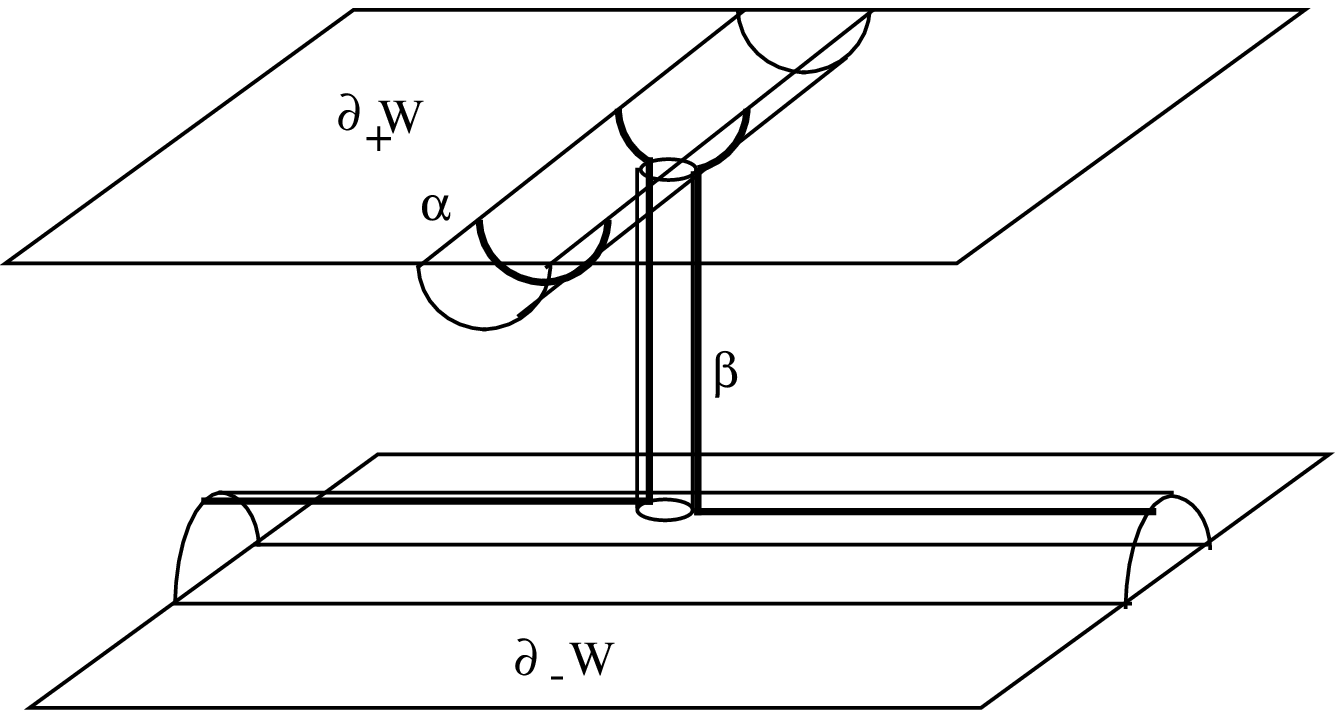}
\relabela <-2pt,0pt> {d2}{$\partial_+W$}
\relabela <-2pt,0pt> {d1}{$\partial_-W$}
\relabela <-2pt,0pt> {a}{$\alpha$}
\relabela <-2pt,0pt> {b}{$\beta$}
\endrelabelbox}
\caption{}
\end{figure} 

\medskip
Following Lemmas \ref{squar} and
\ref{spanner}, we can assume that $P_-$  and $Q_-$ each intersect \Aaa\
in pairs of spanning arcs of \Aaa.  These pairs of arcs can be made
disjoint by a proper isotopy of, say, $P_- \subset W$.  As noted
in the remarks following Proposition
\ref{splitchar}, since we only seek to understand the involutions up
to Dehn twists along $\bdd W$, we may allow proper isotopies in
$W$ that are not fixed on $\bdd P_-$.  Then, following Proposition
\ref{splitchar}, the involutions
$\Thp|W$ and $\Thq|W$ preserve \Aaa\ as well as \Fff, and induce the
standard involution on the solid torus $U = W - \Aaa$. Of the three
possibilities for such an annulus preserving involution of $W$ (see
the proof of Proposition \ref{splitchar}) only one sends each boundary
component of $W$ to itself, as do $\Thp|W \Thq|W$ and $\Thp|W
\Thq|W$.  Hence these products coincide with that involution.   This
implies that the involutions
$\Thp|W$ and $\Thq|W$ commute.  This implies that $\Thp$ and $\Thq$
commute, up to Dehn twists along
$\bdd W$. \qed

\medskip
\noindent {\bf Subcase 2c}\qua $W$ contains no essential spanning annulus.

This case closely parallels that of
Section \ref{positioning}, so some of it will just be sketched.  We
consider the square $I \times I$ parameterizing the sweep-outs by $P_-$
and $Q_-$.  We label any region in which the two surfaces are transverse
with label $A$ if there is a meridian disk of $A_-$ whose boundary
lies entirely in $P_- - Q_-$ or if there is an arc component of $P_-
\cap Q_-$ which $\bdd$--compresses to \Fff\ via a disk in $A_-$. 
Similarly apply labels $B$, $X$, and $Y$.  Labels $A$ and $B$
(or $X$ and $Y$) can't appear on the same or adjacent regions, in part
by Lemma \ref{disjarcs}.  So there will be regions with no labels at
all.  

Consider how $P_- \cap Q_-$ appears in an unlabelled region.  We can
think of the intersection arcs in $P_-$ as a graph $\Ggg_P \subset S^2$
whose edges are the arcs of intersection and whose fat vertices are
disks filling in the four boundary components of $P_-$.  Two of these
vertices, $u^{+}_e$ and $u^{+}_w$ lie on $\bdd^+ P_-$ and two of
them $u^{-}_e$ and $u^{-}_w$ lie on $\bdd^- P_-$. (See Figure 40.)
Similar remarks hold for the graph $\Ggg_Q \subset S^2$ which describes
the arcs of intersection in $Q_-$.  Label the vertices in this graph by
$v^{+}_e$, $v^{+}_w$, $v^{-}_e$ and $v^{-}_w$ in a similar fashion.  


\begin{figure}[ht!] 
\centerline{\relabelbox\small
\epsfxsize=45mm \epsfbox{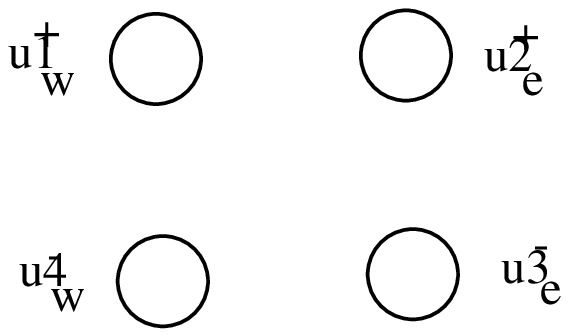}
\relabela <-2pt,0pt> {u1}{$u^+_w$}
\relabela <-2pt,0pt> {u2}{$u^+_e$}
\relabela <-2pt,0pt> {u4}{$u^-_w$}
\relabela <-2pt,0pt> {u3}{$u^-_e$}
\endrelabelbox}
\caption{}
\end{figure}

The valence of each vertex is $2 p \cdot q$, where $p$ and $q$ are the
slopes of $\bdd P_-$ and $Q_-$ in \Fff\ respectively.   Since the
region has no labels, it follows that there are no trivial loops in
$\Ggg_P$ or $\Ggg_Q$, hence no loops at all.  No loops in $\Ggg_P$
means that any edge in $\Ggg_Q$ has one end in one of $v^{+}_e$,
$v^{+}_w$ and the other end in one of $v^{-}_e$ and $v^{-}_w$. (An
orientation parity argument is used here.)  That is, each edge has
one end on a $+$ vertex and one end on a $-$ vertex, in fact in both
graphs.  If three or more edges are parallel in $\Ggg_P$ say, then the
bigons lying between them can be assembled to give a spanning annulus
in $W$, contradicting our hypothesis. So we may as well assume
that $p \cdot q = 1$, so each vertex has valence 2.

Now restrict attention to those regions of $I \times I$ which are
unlabelled.  In positionings corresponding to these regions, $\Ggg_P$
and $\Ggg_Q$ are bipartite graphs, so each face has an even number of
edges.  For each face $F$ in $\Ggg_P$ or $\Ggg_Q$, define the index
to be $$J(F) = \frac{|edges \, in \, \bdd F|}{2} - \chi (F).$$ 
The sum of the indices of all faces in $\Ggg_P$ or $\Ggg_Q$ is $\chi
(P_-) = 2$.  If the sum of the indices of all faces in $P \cap X$
(hence also $P \cap Y$) is odd (resp.\ even) we say the positioning is
$P$--odd (resp.\ $P$--even), and similarly for $Q \cap A$.  Since there
are no loops in either graph, and the valence of each vertex is $2$, to
say a position is $P$--odd is equivalent to saying that both $P \cap X$
and $P \cap Y$ are disks with four edges. (See Figure 41.) Examination
of the few combinatorial possibilities shows that $P$--odd is equivalent
to $Q$--odd, so we will refer to unlabelled regions as either {\em odd}
or {\em even}. Regions which already have labels are neither even nor
odd.  Note that a bigon in $\Ggg_Q$ lying in $A_-$ corresponds to a
properly imbedded square $I \times I \subset (A_- \cap X)$ so that $I
\times \{ 0 \}$ (resp.\ $I \times \{ 1 \}$) is an edge of $\Ggg_P$
running between
$u^{\pm}_e$ (resp.\ $u^{\pm}_w$), and  $\{ 0 \} \times I$ and $\{ 1 \}
\times I$ are spanning arcs of the  annuli $A_- \cap \bdd^{\pm} W$, one
arc in each. Such a square (with two sides spanning the annuli
$\bdd_{\pm}$ and the other two essential arcs in $P_-$) is called a
{\em spanning square} in $A_-$.  No side of a spanning square in
$A_-$ can be isotopic to a side of a spanning square in $B_-$, for
otherwise the two squares could be assembled to give a spanning
annulus, contradicting our hypothesis.


\begin{figure}[ht!] 
\centerline{\relabelbox\small
\epsfxsize=45mm \epsfbox{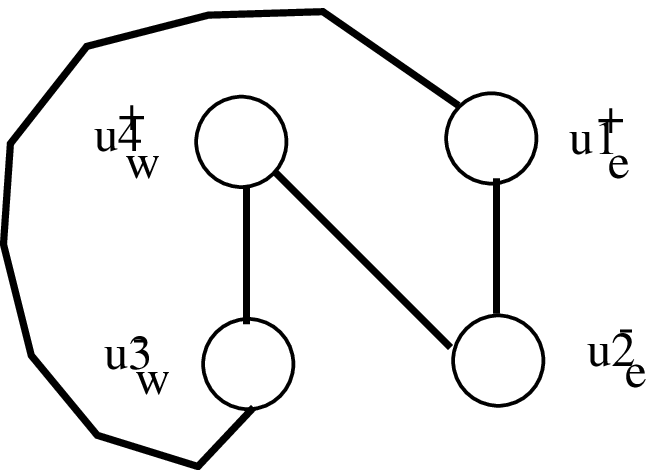}
\relabela <-2pt,0pt> {u4}{$u^+_w$}
\relabela <-2pt,0pt> {u1}{$u^+_e$}
\relabela <-2pt,0pt> {u3}{$u^-_w$}
\relabela <-2pt,0pt> {u2}{$u^-_e$}
\endrelabelbox}
\caption{}
\end{figure}

Expand the rules for labelling, much as
in Section \ref{positioning}, to include the label $A'$ if there is an
arc in $Q - P$ which $\bdd$--compresses to \Fff\ via a $\bdd$--compressing
disk in $A_-$, and similarly for the other three labels $B', X', Y'$. 
(The difference between $A'$ and $A$ is that for the label $A$ the
$\bdd$--compressing arc needs to be an arc of 
$Q \cap P$ whereas for label $A'$ it only needs to lie in $Q - P$.)  The
labels
$A, B, A', B'$ (resp.\
$X, Y, X', Y'$ will be called $P$--labels (resp.\ $Q$--labels.)

\begin{lemma}
\label{persist}

A previously unlabelled region adjacent to a region that has label $A$
now has label $A'$.  {\em Any} region adjacent to a region that has {\em
only} label $A$ either itself has label $A$ or it is even and has
label $A'$.  Similarly for labels $B, X, Y$.

\end{lemma}

\pf The move from the region with label $A$ to the adjacent region
corresponds to a band  move.  The band itself is in a face, and hence
is disjoint from the edge of
$\Ggg_P$ that $\bdd$--compresses in $A$. (See Figure 42.) So the
$\bdd$--compressing disk persists even after the band move, though its
edge in $P_-$ is no longer in the graph. Thus if the adjacent
region had previously been unlabelled it now gets label $A'$.


\begin{figure}[ht!] 
\centerline{\relabelbox\small
\epsfxsize=45mm \epsfbox{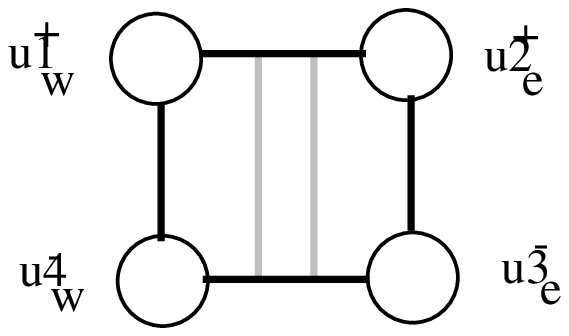}
\relabela <-2pt,0pt> {u1}{$u^+_w$}
\relabela <-2pt,0pt> {u2}{$u^+_e$}
\relabela <-2pt,0pt> {u4}{$u^-_w$}
\relabela <-2pt,0pt> {u3}{$u^-_e$}
\endrelabelbox}
\caption{}
\end{figure}

If the original region has only label $A$ (and not label $X$
or $Y$) then $A_-$ must contain two $\bdd$--compressing disks, one with
edge in $P_-$ running from $u^{+}_e$ to $u^{+}_w$ and the other with
edge running from $u^{-}_e$ to $u^{-}_w$. A band move that destroys
both edges would result in an even region and so one labelled $A'$. 
Otherwise one of the edges persists and the label $A$ remains.
\qed

\begin{lemma}
\label{oddsarc}

Any even region has two primed labels.  Adjacent regions cannot both be
even.  If a region is odd then its labels are a subset (possibly with
primes removed) of the labels of any adjacent region.  

\end{lemma}

\pf A region that is even corresponds to a positioning
where in $\Ggg_P$ there are exactly two edges running between
$u^{\pm}_e$ and two between $u^{\pm}_w$.  The resulting
bigons lie either both in $X$ or both in $Y$, say the former.  (See
Figure 43.)  Then
$P_-$ intersects $X_-$ only in two parallel spanning squares, so $Q_-$
$\bdd$--compresses to $\bdd W$ in the complement of $P_-$, forcing the
label $X'$.  A dual argument works from $\Ggg_Q$ to give a label $A'$
or $B'$.


\begin{figure}[ht!] 
\centerline{\relabelbox\small
\epsfxsize=45mm \epsfbox{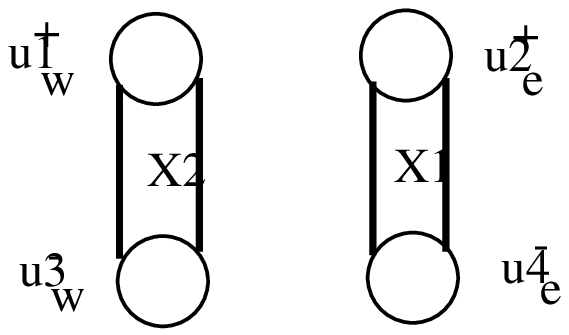}
\relabela <-2pt,0pt> {u1}{$u^+_w$}
\relabela <-2pt,0pt> {u2}{$u^+_e$}
\relabela <-2pt,0pt> {u3}{$u^-_w$}
\relabela <-2pt,0pt> {u4}{$u^-_e$}
\relabela <-2pt,0pt> {X1}{$X$}
\relabela <-2pt,0pt> {X2}{$X$}
\endrelabelbox}
\caption{}
\end{figure}

The band move in $\Ggg_P$ corresponding to a move to an adjacent
region creates either a loop, an edge corresponding to a loop in
$\Ggg_Q$ (eg an edge with one end on each of $u^{+}_e$ and
$u^{+}_w$) or a positioning that is odd.  The former two possibilities
would have given the corresponding region unprimed labels, so it
couldn't be even.

Consider a positioning corresponding to an odd region and, say, $A'$
is a label.  That is, suppose an arc in $P \cap X$, say,
$\bdd$--compresses through
$A$ to $\Fff$.  Then the arc has an end on each of $u^{+}_e$ and
$u^{+}_w$, say. (See Figure 44.) If one performed this
$\bdd$--compression one would see that there is also an arc in $P \cap
Y$ with ends on
$u^{-}_e$ and $u^{-}_w$ that $\bdd$--compresses through $A$.  One of
these two arcs will persist in any adjacent region of the graphic,
since the corresponding change of positioning of $P_-$ with respect to
$Q_-$ is via a band move in either $P_- \cap X$ (so the second one
persists) or
$P_- \cap Y$ (so the first persists).  \qed


\begin{figure}[ht!] 
\centerline{
\epsfxsize=45mm \epsfbox{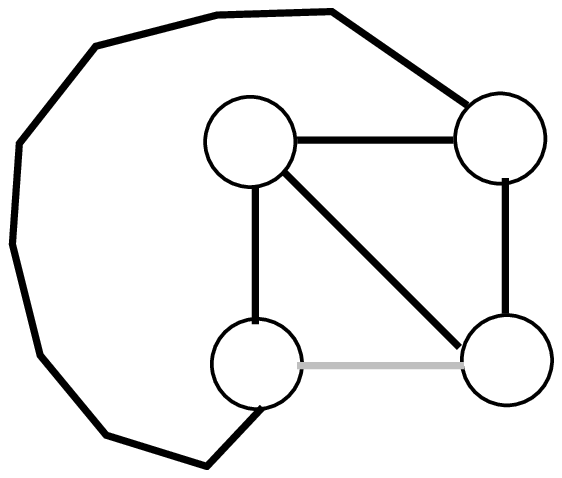}}
\caption{}
\end{figure}

\begin{lemma}
\label{squares}

In a positioning corresponding to an even region, with label $A'$, there
is a properly imbedded square $I \times I \subset (A_- \cap X)$ so that
$I \times \{ 0 \}$ (resp.\ $I \times \{ 1 \}$) is parallel to an edge of
$\Ggg_P$ running between $u^{\pm}_e$ (resp.\ $u^{\pm}_w$), and 
$\{ 0 \} \times I$ and $\{ 1 \} \times I$ are spanning arcs of the 
annuli $A_- \cap \bdd^{\pm} W$, one arc in each.   Similarly for labels
$B', X', Y'$.

\end{lemma}

\pf  Since the label is $A'$, there is a $\bdd$--compressing disk
$D^+$ for $P_-$, lying in $A_-$, one of whose sides is a spanning arc of
the annulus $A_- \cap \bdd^{+} W$, say, and the other side is in $P_-$
but disjoint from the arcs $P_- \cap Q_-$. If one performed this 
$\bdd$--compression one would see that there is also a 
$\bdd$--compressing disk $D^-$ for $P_-$, lying in $A_-$, one of whose
sides is a spanning arc of the annulus $A_- \cap \bdd^{-} W$ and
the other side is also in $P_-$ but disjoint from the arcs $P_- \cap
Q_-$.  Piping these disks together in $P_- - Q_-$ gives the required
square. (See Figure 45.) \qed


\begin{figure}[ht!] 
\centerline{\relabelbox\small
\epsfxsize=45mm \epsfbox{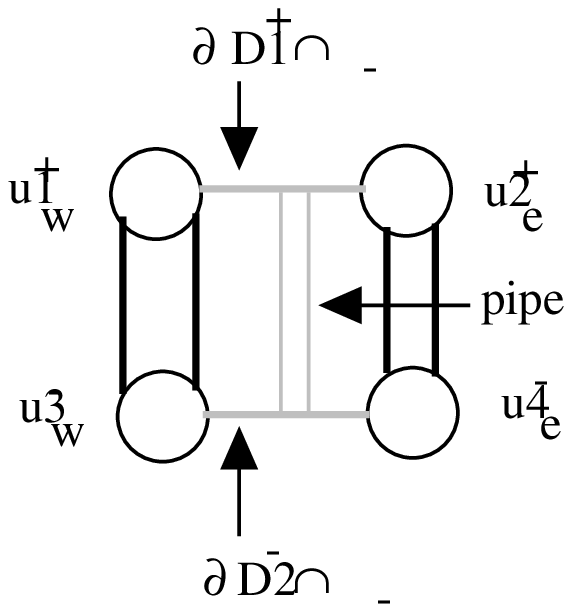}
\relabela <-2pt,0pt> {u1}{$u^+_w$}
\relabela <-2pt,0pt> {u2}{$u^+_e$}
\relabela <-2pt,0pt> {u3}{$u^-_w$}
\relabela <-2pt,0pt> {u4}{$u^-_e$}
\relabela <-8pt,0pt> {D1}{$\partial D^+\cap P_-$}
\relabela <-8pt,0pt> {D2}{$\partial D^-\cap P_-$}
\relabela <0pt,0pt> {pipe}{pipe}
\endrelabelbox}
\caption{}
\end{figure}

\begin{lemma}
\label{notboth}

No region can be labelled both $A$ and $B$ or both $A'$ and
$B'$.  Similarly for labels $X, X', Y, Y'$.

\end{lemma}

\pf  If both labels $A$ and $B$ occur then there would be a
spanning annulus.  If both labels $A'$ and $B'$ occur and the
region is even, then \ref{squares} shows how to construct squares in
both $A_-$ and $B_-$ which assemble to give a spanning annulus.   If the
region is odd, then note that in each of the two faces of $\Ggg_P$
there is only one isotopy class of arcs with one end point on
each of  $u^{+}_e$ and $u^{+}_w$. (See Figure 46.) It follows that the
boundary compressing disks in $A_-$ and $B_-$ either are disjoint or
would assemble to make a compressing disk for $\bdd^+ W$.  The latter
violates the hypothesis and the former would create a spanning
annulus, contradicting our hypothesis. \qed 


\begin{figure}[ht!] 
\centerline{
\epsfxsize=50mm \epsfbox{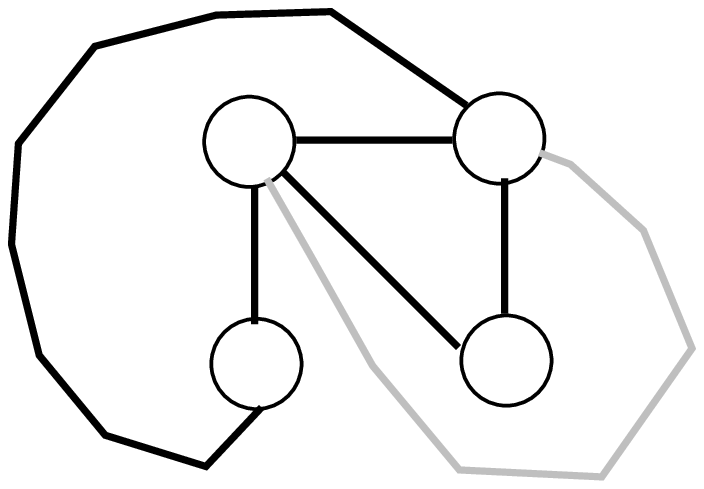}}
\caption{}
\end{figure}

\begin{lemma}
\label{notadj}

No two adjacent regions can be labelled so that one has label $A$ or
$A'$ and the other has label one of $B$ or $B'$.  Similarly for labels
$X, X', Y, Y'$.

\end{lemma}

\pf If adjacent regions are labelled $A$ and $B$ then in fact one
can find disjoint $\bdd$--compressing disks for $P_-$, one of them in
$A_-$ and the other in $B_-$, contradicting \ref{disjarcs}. 

If adjacent regions are labelled $A'$ and $B'$, then by \ref{oddsarc}
one is odd and so one has both labels.  This contradicts \ref{notboth}.

If a region labelled $A'$ is adjacent to one labelled $B$ then by
\ref{oddsarc} and \ref{notboth}, the region labelled $A'$ must be
even. The label $B$ on the adjacent region forces, by \ref{persist}, the
label $B'$ onto the region labelled $A'$. This again contradicts
\ref{notboth}.  \qed

\begin{lemma}

There is an unlabelled region.

\end{lemma}

\pf  Following \ref{notboth} and \ref{notadj}, the alternative is that
there is a vertex whose four adjacent regions are each labelled with one
label, appearing in order around the vertex: $A$ or $A'$, $X$ or $X'$,
$B$ or $B'$, $Y$ or $Y'$.  It follows immediately from \ref{oddsarc}
that no region is odd and no two adjacent regions are even, so at least
one of the labels is not primed.  The labelling then contradicts
\ref{persist}. \qed

\medskip
To complete the proof of Theorem \ref{tori} begin with the positioning
of $P_-$ and $Q_-$ that corresponds to an unlabelled, necessarily odd,
region.  As in \ref{align} one can align $P_-$ and $Q_-$, first pushing
arcs in the quadrilaterals $(P_-)_X$ and $(Q_-)_A$ (say) together and
arcs in the quadrilaterals $(P_-)_Y$ and $(Q_-)_B$ together.  The
result is that $P_-$ and $Q_-$ are aligned except along a set of
bigons, since (essentially by \ref{persist}) no loops can be formed in
either graph by a band move of $P_-$ across $Q_-$.  If two bigons were
parallel there would be a spanning annulus, contradicting
our hypothesis. So, after the alignment, there are exactly four bigons
in
$P_- - Q_-$, one for each possible way of connecting a vertex $u^{+}_e$
or $u^{+}_w$ with $u^{-}_e$ or $u^{-}_w$.  Similarly, there are exactly
four bigons in $Q_- - P_-$, one for each possible way of connecting a
vertex $v^{+}_e$ or $v^{+}_w$ with $v^{-}_e$ or $v^{-}_w$.  Each bigon
corresponds to a spanning square. (See Figure 47.) The picture is now so
explicit that
$P$ and $Q$ can be recognized as Variation 3 of Example \ref{nonsep}.
\qed


\begin{figure}[ht!] 
\centerline{\relabelbox\small
\epsfxsize=120mm \epsfbox{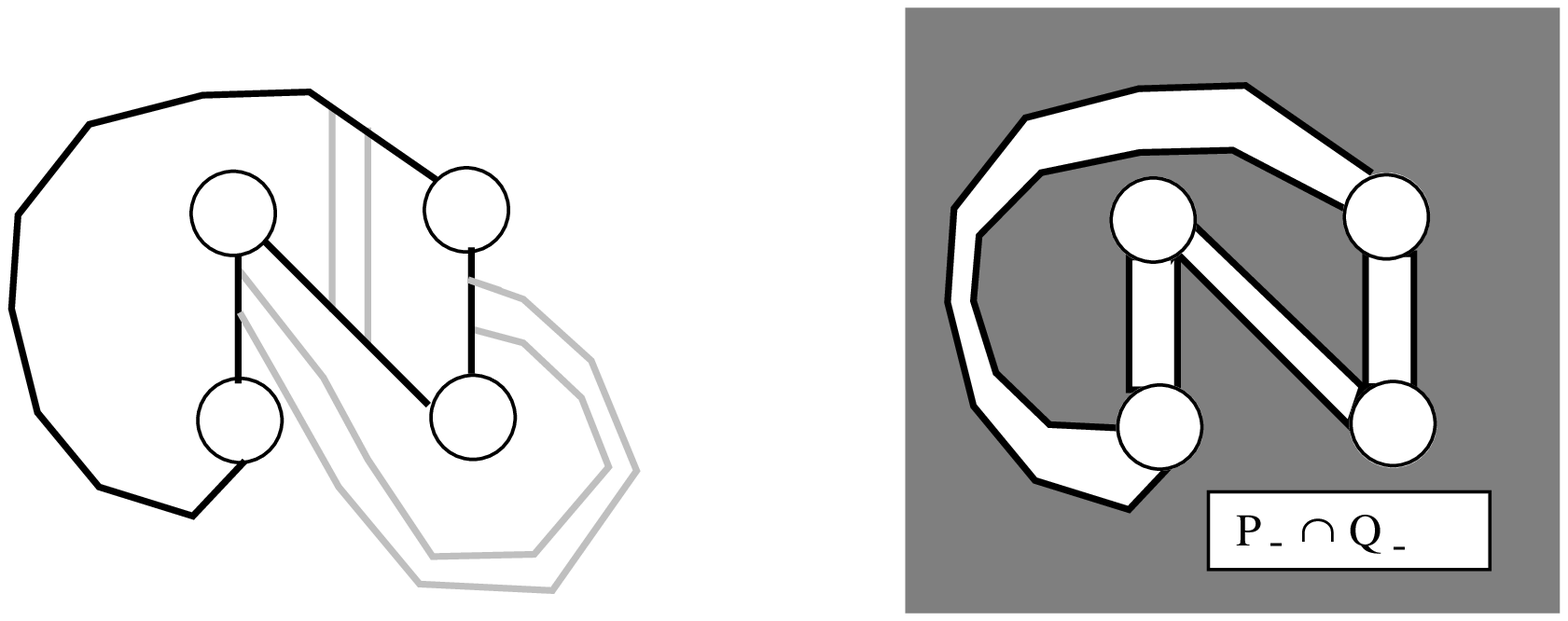}
\relabela <5pt,0pt> {P}{$P_-\cap Q_-$}
\endrelabelbox}
\caption{}
\end{figure}

%
%

%
\Addresses\recd

\begin{thebibliography}



\bibitem{BGM} {\bf J~Birman}, {\bf F~Gonzalez-Acuna}, {\bf J\,M~Montesinos},
{\it Heegaard splittings of prime $3$--manifolds are not unique}, Mich.
Math. J. 23 (1976) 97--103

\bibitem{BH} {\bf  J~Birman}, {\bf M~Hilden}, {\it Heegaard splittings of
branched coverings of $S^3$}, Trans. Amer. Math. Soc. 213
(1975) 315--352

\bibitem{BCZ} {\bf M~Boileau}, {\bf D\,J~Collins}, {\bf H~Zieschang},
{\it Genus 2 Heegaard decompositions of small Seifert manifolds},
Ann. Inst. Fourier, 41 (1991) 1005--1024.

\bibitem{BO} {\bf M~Boileau}, {\bf J-P~Otal}, {\it Groupes des diff{\'e}otopies
de certaines vari{\'e}t{\'e}s de Seifert},  C. R. Acad. Sci.
Paris,  303-I (1986) 19--22

\bibitem{Bo} {\bf F~Bonahon}, {\it Diffeotopies
des espaces lenticulaires}, Topology,
22 (1983) 305--314

\bibitem{BoO} {\bf F Bonahon}, {\bf J-P~Otal}, {\it Scindements de Heegaard des
espaces lenticulaires}, Ann. scient. {\'E}c. Norm. Sup. 16 (1983) 451--466

\bibitem{CG} {\bf A~Casson}, {\bf C\,McA~Gordon}, {\it Reducing Heegaard
splittings}, Topology  and its applications, 27 (1987) 275--283

\bibitem{CS} {\bf D~Cooper}, {\bf M~Scharlemann}, {\it The structure of
a solvmanifold's Heegaard splittings}, to appear

\bibitem{HR} {\bf C~Hodgson}, {\bf H~Rubinstein}, {\it Involutions and
isotopies of lens spaces}, from: ``Knot theory and manifolds'', Lecture Notes in
Mathematics 1144, Springer (1985) 60--96

\bibitem{JS} {\bf M~Jones}, {\bf M~Scharlemann}, {\it How strongly
irreducible Heegaard splittings intersect handlebodies}, to appear

\bibitem{Ja} {\bf W~Jaco}, {\it Lectures on three--manifold topology},
Regional Conference series in Mathematics, no. 43, AMS (1980)

\bibitem{Ko} {\bf T~Kobayashi}, {\it Structures of Haken manifolds with
Heegaard splittings of genus two}, Osaka J. of Math. 
21 (1984) 437--455

\bibitem{Mon} {\bf J\,M~Montesinos}, {\it Sobre la conjectura de Poincare y
los recubridores ramificados sobre un nodo}, Tesis, Facultad de
Ciencias, Universidad Complutense de Madrid (1972)

\bibitem{Mo} {\bf Y~Moriah}, {\it Heegaard splittings of Seifert fibered
spaces}, Invent. Math. 91 (1988) 465--481

\bibitem{MS} {\bf Y~Moriah}, {\bf J~Schultens}, {\it Irreducible Heegaard
splittings of Seifert fibered spaces are either vertical or
horizontal}, Topology, 37 (1998) 1089--1112

\bibitem{RS} {\bf H~Rubinstein}, {\bf M~Scharlemann}, {\it Comparing Heegaard
splittings of non-Haken 3--manifolds}, Topology, 35 (1997) 1005--1026

\bibitem{Sc} {\bf M~Scharlemann}, {\it Heegaard splittings of compact
$3$--manifolds}, from: ``Handbook of Geometric Topology'',
(R~Daverman and R~Sher, editors), Elsevier (to appear)

\bibitem{TO} {\bf M~Takahashi}, {\bf M~Ochiai}, {\it Heegaard diagrams of
torus bundles over $S^1$}, Commentarii Math. Univ. Sancti Pauli, 
31 (1982) 63--69

\bibitem{Vi} {\bf O~Viro}, {\it Linkings, 2--sheeted coverings and braids},
Mat. Sb. (N. S.),  87(129) (1972), 216--228,  English
translation; Math. USSR-Sb. 16 (1972) 223--226

\end{thebibliography}
\end{document}